\documentclass[ijoc,nonblindrev]{informs3_blank}

\OneAndAHalfSpacedXI 

\usepackage[hidelinks]{hyperref}
\hypersetup{
	colorlinks,
	allcolors=black
}
\usepackage{multirow}
\usepackage{subcaption}

\usepackage{natbib}
 \bibpunct[, ]{(}{)}{,}{a}{}{,}%
 \def\bibfont{\small}%

\usepackage{caption}
\usepackage{relsize}
\usepackage{algorithm}
\usepackage{algpseudocode}
\usepackage{hhline}
\usepackage{float}
\usepackage{graphics}
\usepackage{amsmath}
\usepackage{xcolor}
\usepackage{xfrac}
\usepackage{bbm}
\usepackage{booktabs}
\usepackage{placeins}
\usepackage{chngcntr}

\usepackage[font=small]{subcaption}
\usepackage{lmodern}

\newcommand{\cref}[1]{Chapter~\ref{#1}}

\newcommand{\bs}[1]{\boldsymbol{#1}}

\TheoremsNumberedThrough     

\EquationsNumberedThrough    

\usepackage{xcolor}

\begin{document}


\RUNAUTHOR{Bertsimas et al.}

\RUNTITLE{Stochastic Benders Decomposition for Network Design}

 \TITLE{A Stochastic Benders Decomposition Scheme for Large-Scale Stochastic Network Design}

\ARTICLEAUTHORS{%
\AUTHOR{Dimitris Bertsimas}
\AFF{Sloan School of Management and Operations Research Center, Massachusetts Institute of Technology, Cambridge, MA, US, 
	\\ ORCID: \href{https://orcid.org/0000-0002-1985-1003}{0000-0002-1985-1003}
	\\ \EMAIL{dbertsim@mit.edu}}
\AUTHOR{Ryan Cory-Wright}
\AFF{Department of Analytics, Marketing and Operations, Imperial College Business School, London, UK\\
        ORCID: 0000-0002-4485-0619 \\
        \EMAIL{r.cory-wright@imperial.ac.uk}}
\AUTHOR{Jean Pauphilet}
\AFF{Management Science and Operations, London Business School, London, UK
	\\ ORCID: \href{https://orcid.org/0000-0001-6352-0984}{0000-0001-6352-0984}
	\\ \EMAIL{jpauphilet@london.edu}}
\AUTHOR{Periklis Petridis}
\AFF{Operations Research Center, Massachusetts Institute of Technology, Cambridge, MA, USA
	\\ ORCID: \href{https://orcid.org/0000-0002-1019-0763}{0000-0002-1019-0763}
	\\ \EMAIL{periklis@mit.edu}}
} 

\ABSTRACT{%
Network design problems involve constructing edges in a transportation or supply chain network to minimize construction and daily operational costs. We study a stochastic version where operational costs are uncertain due to fluctuating demand and estimated as a sample average from historical data. 
This problem is computationally challenging, and instances with as few as 100 nodes often cannot be solved to optimality using current decomposition techniques. We propose a stochastic variant of Benders decomposition that mitigates the high computational cost of generating each cut by sampling a subset of the data at each iteration and nonetheless generates deterministically valid cuts, rather than the probabilistically valid cuts frequently proposed in the stochastic optimization literature, via a dual averaging technique. We implement both single-cut and multi-cut variants of this Benders decomposition, as well as a variant that uses clustering of the historical scenarios. To our knowledge, this is the first single-tree implementation of Benders decomposition that facilitates sampling.
On instances with 100--200 nodes and relatively complete recourse, our algorithm achieves 5--7\% optimality gaps, compared with 16--27\% for deterministic Benders schemes, and scales to instances with 700 nodes and 50 commodities within hours.
Beyond network design, our strategy could be adapted to generic two-stage stochastic mixed-integer optimization problems where second-stage costs are estimated via a sample average.
}%


\KEYWORDS{Generalized Benders Decomposition; Network Design; Stochastic Integer Optimization}

\maketitle

\linespread{1.2}
\section{Introduction}\label{section:introduction}
Network design is one of the most famous and frequently studied problems in the Operations Research literature, with widespread applications in logistics, air transportation \citep{barnhart2003applications}, supply chains \citep{santoso_stochastic_2005, pishvaee_accelerated_2014}, telecommunications \citep{balakrishnan_models_1991}, and energy markets \citep{binato_new_2001} among other domains. These problems are large-scale and involve uncertain parameters which reflect deviations between the forecast and realized utilization of a network, e.g., uncertain consumer demand in an air traffic control problem or uncertain renewable generation output in a capacity expansion problem. Moreover, we often have data on past realizations of the uncertain parameters. Unfortunately, despite the rapid advances in the scalability of branch-and-bound solvers over the past 25 years, stochastic network design problems with as few as 100 nodes are, to our knowledge, currently regarded as intractable and instead are solved via domain-specific approximation algorithms or heuristics \citep{crainic_network_2021}. 

To scale to network design problems with up to 50 nodes, the mixed-integer optimization (MIO) community has developed a suite of algorithms for mixed-integer nonlinear problems over the past 25 years, originating with the works of \cite{ceria1999convex, stubbs1999branch} and refined by \cite{gunluk_perspective_2009, crainic2016partial} among others. These methods tackle
mixed-integer problems with logical constraints and a partially separable objective function, and 
enforce logical constraints
implicitly via perspective functions, thus tightening the Boolean relaxation. 
Indeed, mixed-integer decomposition schemes that exploit 
perspective reformulations often solve problems to optimality at sizes an order of magnitude larger than was previously possible; see \cite{fischetti_redesigning_2017, bertsimas_unified_2021} for related decomposition schemes.

In a different direction, the machine-learning community has enjoyed considerable success over the past 25 years in improving the scalability of unconstrained stochastic optimization. A common meta-approach is to modify a classical optimization algorithm to sample from an observed dataset at each iteration of the algorithm, and not consider the entire dataset as part of each iterate. Remarkably, each sample often conveys the same essential information as the entire dataset but can be processed multiple orders of magnitude faster. This sampling approach routinely produces a multiple-order-of-magnitude scalability improvement on classical optimization algorithms. Stochastic variants of first-order methods 
such as Stochastic Gradient Descent \citep[SGD,][]{davis2020stochastic}, the Stochastic Average Gradient method \citep{schmidt2017minimizing}, or Adam \citep{kingma2014adam} are currently considered to be state-of-the-art for unconstrained problems. 

In this paper, we propose to embed a sampling technique within a Benders decomposition \citep[]{geoffrion1972generalized} scheme run on the perspective reformulation \citep[]{gunluk_perspective_2009} of a network design problem. To our knowledge, this is the first single-tree implementation of Benders decomposition that facilitates sampling scenarios while maintaining deterministic optimality guarantees. We demonstrate that this approach obtains bound gaps of $5$--$7\%$ on instances with $100$--$200$ nodes, three times smaller than the bound gaps obtained by deterministic Benders decomposition schemes in a comparable amount of time.
Moreover, our approach successfully scales to obtain bound gaps of $10$--$40\%$ on instances with $700$ nodes and $50$ commodities. At this scale, deterministic Benders schemes obtain optimality gaps of $25$--$55\%$.
Our numerical success can be explained by the fact that sampling allows us to generate significantly more Benders cuts within a given time budget than is possible via a deterministic Benders approach, while conveying most of the essential information stored in each deterministic cut. 
Although developed for the special case of stochastic multi-commodity capacitated fixed-charge network design problems, we believe our approach could be applied to two-stage stochastic optimization problems where the first-stage variables are discrete, and the second-stage cost is evaluated via a sample average approximation.

\subsection{Problem Formulation and Main Contributions}
\paragraph{Problem Formulation: }We propose a new approach for solving stochastic Multi-commodity Capacitated Fixed-charge Network Design (MCFND) problems to certifiable optimality, which we formally define in the next paragraph. Similar models appear in \citet[][]{magnanti_network_1984, costa_survey_2005, crainic2016partial, rahmaniani2018accelerating, ramirez2023benders} among other works.

In MCFND problems, there is an index set of commodities $\mathcal{K}$ to be shipped over a capacitated directed network $(\mathcal{N}, \mathcal{E})$, 
where $\mathcal{N}$ denotes a set of nodes and $\mathcal{E}$ denotes a set of edges. Our overall objective is to perform this transshipment in a manner that minimizes the construction plus flow transportation cost. Let $\boldsymbol{A}$ denote this network's corresponding flow conservation matrix. The capacity of arc $(i,j) \in \mathcal{E}$ is given by $u_{i,j}$ and each node $n \in \mathcal{N}$ supplies or demands an amount $d_{n}^{k,r}$ of each commodity $k \in \mathcal{K}$ in each scenario $r \in \mathcal{R}$.
There is a fixed cost $c_{ij}$ of activating each edge $(i,j) \in \mathcal{E}$, and given this problem data, we introduce binary design variables $z_{i,j} \in \{0, 1\}$ to denote whether the $(i,j)$th edge is activated. In addition to taking activation cost into account in the objective, some applications can also involve a fixed limit on the number of edges to be activated, $c_0$. The flow variable $x_{ij}^{k,r}$ then denotes the quantity of commodity $k$ routed on edge $(i,j)$ in scenario $r$, and $f_{ij}^{k}$ denotes the marginal transportation cost, i.e., the per unit cost of transporting the $k$th commodity through edge $(i,j)$. 
Moreover, we follow the standard Sample Average Approximation (SAA) paradigm \citep[see][for a general theory]{shapiro2021lectures} in placing equal weight on each observation of historical data $r$ in our objective.

The complete optimization formulation for MCFND can then be written as:
\begin{equation}\label{eqn:prob_description}
	\begin{aligned}
		\min \quad & \sum_{(i,j) \in \mathcal{E}} c_{i,j}z_{i,j} + 
		\dfrac{1}{|\mathcal{R}|} \sum_{r \in \mathcal{R}} \sum_{(i,j) \in \mathcal{E}} \left( \sum_{k \in \mathcal{K}} f_{i,j}^{k}x_{i,j}^{k,r} + \frac{1}{2 \gamma} \left( \sum_{k \in \mathcal{K}} x_{i,j}^{k,r}\right)^2 \right) &
		\\
		\text{s.t.} 
		\quad & 
		\boldsymbol{A}\boldsymbol{x}^{k,r} = \boldsymbol{d}^{k,r}, \quad \forall k \in \mathcal{K}, r \in \mathcal{R},\\
		\quad & 
		\sum_{k \in \mathcal{K}} x^{k,r}_{i,j} \leq u_{i,j},  
		\quad \forall (i,j) \in \mathcal{E}, r \in \mathcal{R}, \\
		& 
		\boldsymbol{x}^{k,r} \geq 0 ,\ x^{k,r}_{i,j} = 0 \textit{\rm\ if } z_{i,j} = 0, \quad \forall (i,j) \in \mathcal{E},\\
        \quad & \sum_{(i,j) \in \mathcal{E}} z_{i,j} \leq c_0, \ z_{i,j} \in \{0, 1\} \quad \forall (i,j) \in \mathcal{E},
	\end{aligned}
\end{equation}
where $\gamma > 0$ controls a strongly quadratic regularization term in the objective, which can be seen as a penalization of the hard constraint on each edge's capacity, $\sum_{k \in \mathcal{K}} x^{k,r}_{i,j} \leq u_{i,j}$  (see also \citet{atamturk2018note} for a discussion of capacity constraints in network design problems). We refer to this term as a ``regularization'' term throughout the paper, and justify its use from both a theoretical and a practical perspective in Appendix \ref{sec:a.regularization.theory}. 

Observe that in Problem \eqref{eqn:prob_description}, we link the discrete and continuous decisions in \eqref{eqn:prob_description} with a logical `if' statement. 
In the network design literature, these logical constraints are typically replaced with big-M constraints of the form $\sum_{k \in \mathcal{K}}x_{i,j}^{k,r}\leq u_{i,j}z_{i,j}$ by default \citep[]{glover1975improved}. However, there are alternative ways to convexify logical constraints, which sometimes lead to tighter formulations, e.g., by leveraging the presence of the strongly quadratic term in the objective ---leading to the so-called perspective formulation, with second-order cone constraints \citep{ceria1999convex,gunluk_perspective_2009}--- or by leveraging both the quadratic term and the capacity constraints (as we do in this paper). 
Accordingly, we formulate network design with logical constraints to facilitate tighter convexifications and stronger Benders cuts than are achievable via the big-$M$ technique alone; see also \cite{wei2022ideal} for a detailed study of conic formulations that give tighter relaxations of logically constrained problems than big-$M$ relaxations in other contexts.

\paragraph{Main Contributions: }In this paper, we provide two main contributions. 

First, we propose a new 
decomposition method that combines sampling-based methods from the stochastic optimization and machine learning literature with a Generalized Benders Decomposition approach in the spirit of \cite{geoffrion1972generalized}. 
Our approach can tackle large-scale mixed-integer problems 
by leveraging weak duality to obtain valid dual variables for scenarios we do not explicitly sample. To our knowledge, this is the first single-tree implementation of a Generalized Benders Decomposition scheme that facilitates sampling.

Second, we implement and benchmark our approach across a wide variety of large-scale network design instances, and explore the performance benefits of various design and implementation choices. 
Our approach allows us to solve network design problems with $200$ nodes and relatively complete recourse to within $7\%$ of optimality in hours, and obtain high-quality feasible solutions on instances with relatively complete recourse and up to $700$ nodes. 

\subsection{Background and Literature Review}
Our work is built on two intertwined literatures. 
First, decomposition schemes for large-scale deterministic problems with logical constraints developed by the MIO community. Second, sampling algorithms for problems with exogenous uncertainty developed by the stochastic optimization community. 
We further remark that, owing to Problem \eqref{eqn:prob_description}'s significant computational difficulty, a wide variety of approximation algorithms \citep{agrawal1991trees, goemans1993survivable, bertsimas1998valid} and heuristic methods have also been proposed for solving Problem \eqref{eqn:prob_description}; see \citet{rodriguez2010local, gendron2018matheuristics} for reviews.

\paragraph{Cutting-Plane Schemes for Mixed-Integer Optimization: }
Problem \eqref{eqn:prob_description} is a computationally challenging mixed-integer problem that encompasses hard combinatorial problems such as Steiner tree optimization \citep[][]{garey1977rectilinear} and possesses extremely poor Boolean relaxations \citep{gendron1999multicommodity}. Indeed, generic branch-and-bound solvers cannot currently solve network design (ND) problems at even moderate problem sizes with tens of nodes \citep[see][Section 6.1 for an investigation of \texttt{CPLEX} version $12.8$'s performance on synthetic ND instances with ten nodes]{crainic_partial_2021}. Accordingly, and due to its cardinal importance in practice, ND has emerged as one of the most frequently studied problems in the MIO literature over the past 50 years. 

Throughout the first 30 years of the field of Operations Research, there was a spirited debate regarding the most efficient technique for solving ND problems, with many proposals, including branch-and-bound \citep{boyce1973optimal}, Lagrangian methods \citep{cornuejols1980canonical}, and dynamic programming \citep{erickson1987send}. The idea of solving ND problems via Generalized Benders decomposition \citep[][]{geoffrion1972generalized} was moved front-and-center by \cite{magnanti1981accelerating, magnanti_network_1984}. Building upon several influential prior works, including \cite{geoffrion1974multicommodity, florian1976engine, richardson1976optimization}, they found that an accelerated Benders decomposition was a viable and often more scalable alternative for ND problems than several other optimization approaches, including the three aforementioned ones. 
Ever since, Benders decomposition has been widely recognized as one of the most competitive methods for solving ND problems; we refer to \cite{fischetti_redesigning_2017, crainic_partial_2021} for modern reviews of Benders decomposition for ND problems. 

In a related direction, a significant line of work has developed a suite of cutting planes that iteratively strengthen Problem \eqref{eqn:prob_description}'s Boolean relaxation upon their imposition; see, e.g., \cite{van1985valid, magnanti1993convex, magnanti1995modeling, bienstock1998minimum, gunluk1999branch, atamturk2021multicommodity} and references therein. Remarkably, these approaches are so numerically successful and easy to implement that they are usually incorporated within commercial branch-and-cut solvers within several years of their proposal \citep[][]{bixby2012brief}. 
As a result, some of the decomposition schemes reviewed above may even be considered ``sleeping beauties'' in the sense of \cite{ke2015defining}, i.e., were not originally considered numerically successful but would be if proposed today, implicitly in conjunction with these valid inequalities.

\paragraph{Decomposition Schemes for Large-Scale Optimization Under Uncertainty: } Cotemporally, a considerable amount of attention has been devoted by the stochastic optimization community to solving large-scale convex optimization problems with uncertain parameters for which we have access to either a joint probability distribution or observations from historical data. Initiated by the independent works of \cite{dantzig1955linear,beale1955minimizing}, and subsequently refined by \cite{wets1966programming, van1969shaped}, contemporary optimizers for large-scale stochastic problems typically invoke 
the Minkowski-Weyl theorem \citep[c.f.][Chapter 4]{bertsimas1997introduction} to solve their deterministic equivalents via Benders decomposition \citep[which was termed the L-shaped method by][]{van1969shaped}. Alternatively, 
works like \citet{zakeri_inexact_2000, fabian2000bundle, rei2009accelerating, guigues2020inexact} propose generating Benders cuts without solving each subproblem to optimality.

The two main variants of Benders decomposition invoked for two-stage stochastic integer optimization problems such as Problem \eqref{eqn:prob_description} are called single-cut and multi-cut Benders. Single-cut schemes maintain a single epigraph variable that upper bounds the expected transshipment cost and generates a single cut at each iteration of Benders decomposition. Multi-cut schemes associate a separate epigraph variable with the cost incurred in each scenario and generate a separate cut for each epigraph variable in each iteration \citep{birge_multicut_1988}. Therefore, single-cut schemes typically require more iterations to converge but require less time to perform each iteration \citep[see][for comparisons]{birge_introduction_2011,de2008benders,you2013multicut}. 
Problems with fewer scenarios are typically solved faster via multi-cut approaches. However, the relative performance of each variant is highly problem-dependent.

More recently, considerable attention has been devoted to designing variants of Benders decomposition that avoid solving a subproblem for each scenario at each iteration by sampling. \cite{higle1991stochastic, pereira1991multi, dantzig_multi-stage_1993, infanger_monte_1992} initiated this line of inquiry by proposing stochastic cutting-plane schemes that converge almost surely 
\citep[see also][]{bertsimas_stochastic_2021}. 
Determining convergence of these schemes is technically challenging. Various statistical tests exist \citep[see, e.g.,][]{higle1996duality, morton1998stopping, mak1999monte} that provide confidence intervals on the duality gap. Yet, to avoid multiple-testing problems, practitioners typically run stochastic cutting-plane methods for a prespecified number of iterations and then perform a statistical test on termination \citep[][]{de2015improving}. 

\subsection{Structure}
We propose a stochastic Benders decomposition scheme that combines the perspective reformulation technique from the MIO literature with sampling ideas from the stochastic optimization literature to, for the first time, successfully solve data-driven capacitated network design problems with hundreds of nodes to certifiable (near) optimality. 
The rest of this paper is laid out as follows:
\begin{itemize}
\item In Section \ref{section:cp_method}, we propose stochastic variants of the single- and multi-cut versions of Benders decomposition to solve a perspective reformulation of \eqref{eqn:prob_description}. Our algorithms randomly sample a subset of scenarios $\mathcal{R}_t \subseteq \mathcal{R}$ at each iteration and use a dual averaging technique to generate cuts that are deterministically valid for all $r \in \mathcal{R}$, while previous stochastic approaches generate cuts that are only valid on average or with high probability. We prove high probability bounds on the approximation error stemming from our dual averaging technique. 
\item In Section \ref{section:bounds}, we propose rigorous convergence criteria to terminate our stochastic decomposition schemes at a certifiable optimal solution. Since our master optimization problem is an MIO problem, we also discuss the specific termination challenges arising when Benders Decomposition is implemented via branch-and-cut (or lazy constraints). We also review techniques for accelerating the convergence of our methods, by warm-starting their upper and lower bounds.
\item In Section \ref{section:experiments}, we apply our decomposition schemes to a collection of network design instances that are synthetically generated or obtained from the literature \citep{crainic2016partial, crainic_network_2021}. 
On the synthetic instances (which exhibit relatively complete recourse),
our best stochastic cutting-plane strategy achieves 7--11\% (resp. 20--30\%) optimality gaps within two hours for instances with 70--300 nodes (resp. 500-700 nodes) compared with 12--26\% (resp. 50--55\%) for its deterministic counterpart.
On the {\bf R} instances introduced by \citet{crainic_simplex-based_2000} and frequently benchmarked against in the literature, we find that our approach provides a noticeable reduction in optimality gap (by 5--10 percentage points) as the number of scenarios increases.
\end{itemize}

\subsection*{Notation}
We let non-boldface characters such as $b$ denote scalars, lowercase bold-faced characters ($\bs{x}$) denote vectors, uppercase bold-faced characters ($\bs{A}$) denote matrices, and calligraphic uppercase characters ($\mathcal{Z}$) denote sets. We let $[n]$ denote the running set of indices $\{1, \dots, n\}$. We let $\bs{e}$ denote the vector of ones, and $\bs{0}$ denote the vector of all zeros. 
Finally, we let $\langle \cdot, \cdot \rangle$ denote the Euclidean inner product between two vectors of the same size, i.e., $\langle \bs{x}, \bs{y}\rangle:=\sum_{i=1}^n x_i y_i$ for any $\bs{x}, \bs{y} \in \mathbb{R}^n$. 

\section{Deterministic and Stochastic Cutting-Plane Methods}\label{section:cp_method}
This section proposes an efficient numerical strategy for solving Problem \eqref{eqn:prob_description} to certifiable optimality. The backbone of our approach is a Generalized Benders Decomposition scheme run on a perspective reformulation of Problem \eqref{eqn:prob_description}, which uses sampling techniques to avoid explicitly solving each scenario at each iteration of the method. Instead, we use dual-optimal solutions from the sampled subproblems to construct dual-feasible solutions to the remaining subproblems and thereby construct valid cuts. We further discuss the convergence properties of our method. 

\subsection{A Two-Stage Reformulation} \label{ssec:cp.twostage}
We observe that the flow minimization problem with respect to each $\boldsymbol{x}^{\cdot,r}$ in \eqref{eqn:prob_description} is decomposable across scenarios $r \in \mathcal{R}$. Therefore, consider a set of demand vectors $\boldsymbol{d}^k \in \mathbb{R}^{\mathcal{N}}$ for $k \in \mathcal{K}$ and define 
\begin{align}\label{eqn:f_z_definition_y}
f(\bs{z}; \bs{d}) := \min_{\boldsymbol{x}^k \in \mathbb{R}_+^{\mathcal{E}}, k \in \mathcal{K}} \:
 \sum_{k \in \mathcal{K}} \langle \bs{f}^{k}, \bs{x}^{k}\rangle  +\frac{1}{2 \gamma} \sum_{{(i,j)} \in\mathcal{E} } \left(  \sum_{k \in \mathcal{K}}  {x}_{i,j}^{k} \right)^2 \mbox{ s.t. } & 
\bs{Ax}^{k} = \bs{d}^{k}, \, \forall k \in \mathcal{K}, \\
&  \sum_{k \in \mathcal{K}} x^{k}_{i,j} \leq u_{i,j}, \, \forall (i,j) \in \mathcal{E},  \nonumber
\\
&
x^{k}_{i,j} = 0 \textit{ if } z_{i,j} = 0, \, \forall (i,j) \in \mathcal{E},\nonumber
\end{align}
to be the operational cost of serving demand $\bs{d}$ on network $(\mathcal{N}, \mathcal{E})$ with design variables $\bs{z}$.
Observe that the minimization problem defining $f(\bs{z}; \bs{d})$ is not decomposable across commodities because of shared capacity constraints. With this notation, Problem \eqref{eqn:prob_description} is equivalent to
\begin{equation}\label{eqn:prob_outer_approx}
	\min _{\boldsymbol{z} \in \mathcal{Z}} \: \langle \bs{c}, \boldsymbol{z} \rangle + \dfrac{1}{|\mathcal{R}|}\sum_{r \in \mathcal{R}} f(\bs{z}; \bs{d}^{\cdot,r}),
\end{equation}
where $\bs{d}^{\cdot,r}$ denotes the collection of demand vectors $\{ \bs{d}^{k,r}$, $k \in \mathcal{K}\}$ and 
$\mathcal{Z} = \{ \bs{z} \in \{0,1\}^{\mathcal{E}} : \sum_{(i,j)} z_{i,j} \leq c_0 \}$ denotes the set of feasible edges. 
The network design formulation \eqref{eqn:prob_outer_approx} separates the discrete design variables $\bs{z}$ from the continuous second-stage routing variables $\bs{x}^k$, thus giving a pure integer optimization formulation that is readily amenable to outer-approximation techniques.

\subsection{A Linear Lower Approximation of the Second-Stage Cost Function} \label{ssec:cp.duality}
In this section, we derive a family of Benders cuts that successfully outer-approximate a perspective reformulation of \eqref{eqn:f_z_definition_y}. 

Since the objective function in \eqref{eqn:prob_outer_approx} involves the average of the function $f(\bs{z},\bs{d})$ over $|\mathcal{R}|$ realizations of $\bs{d}$, we start by analyzing properties of the function $f(\bs{z},\bs{d})$ in isolation, with a view to establish that $f(\bs{z},\bs{d})$ is convex in $\bs{z}$ and a valid subgradient can be obtained by solving a dual problem, as has already been done in the literature, e.g., in \citet{bertsimas_unified_2021}.

\begin{proposition}\label{prop:f.dual} For any $\bs{z} \in \{0,1\}^{\mathcal{E}}$ and demand vectors $\bs{d}^k$, $k \in \mathcal{K}$ such that Problem \eqref{eqn:f_z_definition_y} admits a feasible solution, we have: 
\begin{equation}\label{eqn:f.dual}
\begin{aligned} 
f(\bs{z}; \bs{d}) = & \max_{\substack{\boldsymbol{\alpha} \in \mathbb{R}^{\mathcal{E}}, \boldsymbol{\beta} \in \mathbb{R}_+^{\mathcal{E}} \\ \boldsymbol{p}^k \in \mathbb{R}^{\mathcal{N}}, k \in \mathcal{K}}} \quad
		\sum_{k \in \mathcal{K}} 
		\left\langle \bs{p}^{k}, \bs{d}^{k} \right\rangle 
		- \sum_{{(i,j)} \in\mathcal{E} } 
		 z_{i,j}  \left[ u_{i,j} \beta_{i,j} + \frac{\gamma}{2} (\alpha_{i,j} + \beta_{i,j})^2 \right]
  \mbox{ s.t. } \bs{A}^\top \bs{p}^{k} \leq \bs{f}^{k} - \bs{\alpha}. 
\end{aligned}
\end{equation}

\end{proposition}
The proof of Proposition \ref{prop:f.dual} follows analogously to \citet[Theorem 2.5]{bertsimas_unified_2021} and relies on deriving the dual of the minimization problem defining $f(\bs{z}; \bs{d})$ by using a variable decomposition {\it \`a la Fenchel}; for completeness, we provide a formal proof in Appendix \ref{ssec:a.proof.fdual}. Observe that the optimization problem \eqref{eqn:f.dual} remains well defined if there are no hard constraints on edge capacity (i.e., if $u_{i,j} = +\infty$, we set $\beta_{i,j}=0$) or if there is no quadratic term in the objective (i.e., if $\gamma = +\infty$ we set $\alpha_{i,j}+\beta_{i,j}=0$).

Proposition \ref{prop:f.dual} calls for a few observations. First, according to the dual reformulation, $f(\bs{z}; \bs{d})$ can be expressed as the point-wise maximum of affine functions in $\bs{z}$, hence $f(\bs{z}; \bs{d})$ is convex in $\bs{z}$. Second, any feasible dual solution $\bs{\alpha}  \in \mathbb{R}^{\mathcal{E}}$, $\bs{\beta}\in \mathbb{R}_+^{\mathcal{E}}$, $\bs{p}^k \in \mathbb{R}^{\mathcal{N}}$ such that $\bs{A}^\top \bs{p}^{k} \leq \bs{f}^{k} - \bs{\alpha}$ provides a valid linear lower approximation of $f(\bs{z}; \bs{d})$. Namely, for any $\bs{z}$,
\begin{align*}
f(\bs{z}; \bs{d}) \geq \sum_{k \in \mathcal{K}} 
		\langle \bs{p}^{k}, \bs{d}^{k} \rangle 
        - \sum_{{(i,j)} \in\mathcal{E} } 
		 z_{i,j}  \left[ u_{i,j} \beta_{i,j} + \frac{\gamma}{2} (\alpha_{i,j} + \beta_{i,j})^2 \right].
\end{align*}
When the dual variables are optimal for a particular vector $\bs{z}^0$, the resulting offset and slope in the above linear approximation are exactly the value of $f(\bs{z}^0;\bs{d})$ and a subgradient of $f$ at $\bs{z}^0$, i.e., 
\begin{align*}
f(\bs{z}; \bs{d}) \geq f(\bs{z}^0; \bs{d}) + \left\langle \nabla f(\bs{z}^0; \bs{d}), \bs{z} - \bs{z}^0 \right\rangle.
\end{align*}
Third, Proposition \ref{prop:f.dual} applies if Problem \eqref{eqn:f_z_definition_y} is feasible for the current design vector $\bs{z} = \bs{z}^0$. On the other hand, if \eqref{eqn:f_z_definition_y} is not feasible, 
then the following feasibility problem does not admit a solution:
\begin{align*}
\exists \, \bs{x} \in \mathbb{R}_+^{\mathcal{E} \times \mathcal{K}}: \ \bs{A}\bs{x}^k=\bs{d}^k \ \forall k \in \mathcal{K}, \ \sum_{k \in \mathcal{K}} x_{i,j}^k \leq u_{i,j} z^0_{i,j}, \ \forall (i,j) \in \mathcal{E}.
\end{align*}
Hence, by Farkas's lemma \citep[see, e.g.,][Theorem 4.6]{bertsimas1997introduction}, we can find a certificate of infeasibility, i.e., we can find $\bs{\beta} \in \mathbb{R}_+^{\mathcal{E}}, \bs{p}^k \in \mathbb{R}^{\mathcal{N}}, k \in \mathcal{K}$ such that
$\bs{A}^\top \bs{p}^k - \bs{\beta} \leq \bs{0}$ and 
$\displaystyle \sum_{k \in \mathcal{K}} \langle \bs{p}^{k}, \bs{d}^{k} \rangle - \sum_{(i,j) \in \mathcal{E}} z^{0}_{i,j} u_{i,j} \beta_{i,j} > 0$. In particular, the existence of such vectors $\bs{\beta}, \, \{\bs{p}^k\}_{k \in \mathcal{K}}$ implies that Problem \eqref{eqn:f.dual} is unbounded.
Therefore, we can separate the infeasible incumbent solution $\bs{z}^0$ 
by imposing the feasibility cut
\begin{align} \label{eqn:feasibility.cut}
\sum_{k \in \mathcal{K}} 
		\langle \bs{p}^{k}, \bs{d}^{k} \rangle -  
		 \sum_{(i,j) \in \mathcal{E}} z_{i,j} u_{i,j} \beta_{i,j} \leq 0
\end{align}
on the first-stage variable $\bs{z}$.

Finally, as has already been observed in the literature \citep{xie2020scalable, bertsimas_unified_2021}, our reformulation can alternatively be achieved by performing a perspective reformulation on \eqref{eqn:f_z_definition_y} to rewrite it as a mixed-integer second-order cone problem \citep[c.f.][]{gunluk_perspective_2009} and taking the dual of this perspective reformulation with respect to the continuous variables. 

\subsection{Epigraph Formulations: Modeling Choice and Algorithmic Implications} \label{ssec:cp.benders}
In this section, we exploit our previously developed characterization of $f(\bs{z}, \bs{d})$ as the pointwise maximum of functions linear in $\bs{z}$ to revisit three deterministic outer-approximation methods that solve Problem \eqref{eqn:prob_description} to certifiable optimality. For simplicity, we focus our description on 
optimality cuts in this section; feasibility cuts follow in much the same way. 

Outer-approximation methods such as generalized Benders decomposition solve \eqref{eqn:prob_outer_approx} by constructing a lower approximation of the second-stage operational cost $\tfrac{1}{|\mathcal{R}|} \sum_{r \in \mathcal{R}} f(\bs{z}; \bs{d}^{\cdot,r})$ and refining this approximation at each step. However, since the second-stage cost is the average operational cost over $|\mathcal{R}|$ scenarios, one can either approximate each term $f(\bs{z}; \bs{d}^{\cdot,r})$ separately or their sum, 
which we refer to as multi-cut and single-cut approaches respectively.

In a multi-cut approach, we consider the following epigraph formulation of Problem \eqref{eqn:prob_outer_approx}, as originally proposed by \cite{birge_multicut_1988} for two-stage stochastic linear optimization: 
\begin{equation*}
	\min _{\substack{\boldsymbol{z} \in \mathcal{Z} \\ \eta_r \in \mathbb{R}, \forall r\in \mathcal{R}}} \langle \bs{c}, \boldsymbol{z} \rangle + \dfrac{1}{|\mathcal{R}|}\sum_{r \in \mathcal{R}} \eta_r \mbox{ s.t. } \eta_r \geq f(\bs{z}; \bs{d}^{\cdot,r}), \, \forall r \in \mathcal{R},
\end{equation*}
and iteratively refine a piecewise linear lower approximation of $f(\bs{z}; \bs{d}^{\cdot,r})$ for each epigraph constraints until convergence. Specifically, at each iteration $T$, the multi-cut cutting-plane algorithm solves the MIO problem
\begin{equation} \label{eqn:cp.multi}
	\min _{\substack{\boldsymbol{z} \in \mathcal{Z} \\ \eta_r \in \mathbb{R}, \forall r\in \mathcal{R}}} \langle \bs{c}, \boldsymbol{z} \rangle + \dfrac{1}{|\mathcal{R}|}\sum_{r \in \mathcal{R}} \eta_r \mbox{ s.t. } \eta_r \geq f(\bs{z}^t; \bs{d}^{\cdot,r}) + \left\langle \nabla f(\bs{z}^t; \bs{d}^{\cdot,r}), \bs{z} - \bs{z}^t\right\rangle, \, \forall t \in [T], \, \forall r \in \mathcal{R}.
\end{equation}
Observe that, in this implementation, each of the $|\mathcal{R}|$ functions $f(\bs{z}; \bs{d}^{\cdot,r})$ is linearized at $T$ points $\bs{z}^t$, so \eqref{eqn:cp.multi} comprises $|\mathcal{R}| \times T$ linear constraints. The solution of \eqref{eqn:cp.multi}, $\bs{z}^{T+1}$, then serves as a linearization point to further improve the approximations of the functions $f(\bs{z}; \bs{d}^{\cdot,r})$ at the next iteration.

Alternatively, the single-cut approach, as originally proposed for two-stage stochastic linear optimization by \cite{van1969shaped}, considers a more compact epigraph formulation: 
\begin{equation*}
	\min _{\substack{\boldsymbol{z} \in \mathcal{Z} \\ \eta \in \mathbb{R}}} \langle \bs{c}, \boldsymbol{z} \rangle + \dfrac{1}{|\mathcal{R}|} \eta \mbox{ s.t. } \eta \geq \sum_{r \in \mathcal{R}} f(\bs{z}; \bs{d}^{\cdot,r}),
\end{equation*}
and constructs a piece-wise linear lower-approximation of $\sum_{r \in \mathcal{R}} f(\bs{z}; \bs{d}^{\cdot,r})$ directly. In a single-cut cutting-plane algorithm, at a given iteration $T$, the epigraph constraint is replaced by linear constraints of the form 
\begin{equation}\label{eqn:cp.single.constraint}
	\eta \geq \sum_{r \in \mathcal{R}} f(\bs{z}^t; \bs{d}^{\cdot,r}) + \left\langle \sum_{r \in \mathcal{R}} \nabla f(\bs{z}^t; \bs{d}^{\cdot,r}), \bs{z} - \bs{z}^t\right\rangle. 
\end{equation}
The single-cut approach involves only one epigraph variable $\eta$ (compared with $|\mathcal{R}|$ in the multi-cut implementation) and adds one linear constraint at each iteration (vs. $|\mathcal{R}|$). As a result, the MIO problems involved in the single-cut approach are smaller and usually more tractable than those solved by the multi-cut approach. Yet, multi-cut methods approximate the second-stage cost function more accurately and might require fewer iterations to converge. Various studies, including \citet[][]{birge_introduction_2011, de2008benders, you2013multicut} have reported mixed results on the relative merits of single and multi-cut methods, and which method works best appears to depend on the underlying problem and the number of scenarios. 

Regarding feasibility, if there exists a scenario $\bs{d}^{\cdot,r}$ for which the incumbent solution $\bs{z}^T$ is not feasible (i.e., $f(\bs{z}^T,\bs{d}^{\cdot,r}) = +\infty$), then a feasibility cut of the form \eqref{eqn:feasibility.cut} is imposed. In the single-cut approach, the feasibility cut is imposed instead of an optimality cut \eqref{eqn:cp.single.constraint}. However, in the multi-cut approach, optimality cuts on the other epigraph variables $\eta_{r'}, r' \neq r$, which correspond to feasible scenarios, can still be added.

We remark that all these methods converge in a finite but possibly exponential number of iterations by the finiteness of $\{0, 1\}^\mathcal{E}$ and since no method visits a binary vector $\bs{z}$ twice \citep[see also][Theorem 2.4]{geoffrion1972generalized}.
A common thread between these approaches is that evaluating values of functions of the form $f(\bs{z},\bs{d})$ (and their subgradients) ---an operation referred to as the separation oracle--- is the main computational bottleneck, and the number of function evaluations is the same, $|\mathcal{R}|$, which can be prohibitive, especially when the number of past scenarios $|\mathcal{R}| $ increases. Accordingly, we propose stochastic versions of these approaches with improved per-iteration complexity in the next section.

\begin{remark} \OneAndAHalfSpacedXI
To successfully combine the best aspects of single and multi-cut approaches,
\cite{trukhanov2010adaptive, contreras2011benders} proposed to partition the scenarios into subsets of similar scenarios and introduce one epigraph variable per cluster. For conciseness, we discuss this approach (and propose a stochastic variant), which we refer to as a $k$-cut approach, in Appendix \ref{sec:a.kcut}.
\end{remark}

\subsection{A Stochastic Cutting-Plane Algorithm} 
\label{ssec:cp.stochastic}
In this section, we propose stochastic variants of the cutting-plane methods proposed in the previous section, which obtain high-quality deterministically valid lower bounds without explicitly solving an optimization problem in each scenario and each commodity at each iteration of the method. We also discuss the convergence of these methods. As these methods do not provide deterministically valid upper bounds from a single sample, we defer a detailed discussion of their upper bounds, the corresponding termination criteria, and their single-tree implementation to Section \ref{section:bounds}, and assume for ease of exposition that all cutting-plane methods are multi-tree throughout the section.

First, a stochastic variant of the multi-cut algorithm can be developed in a straightforward manner. Indeed, in its deterministic implementation, at each iteration $t$ of the multi-cut cutting-plane algorithm, we add one linear constraint for each epigraph variable $\eta_r$, for $r \in \mathcal{R}$. Instead, we can sample a subset $\mathcal{R}_t \subseteq \mathcal{R}$ of scenarios and only add linear constraints for these scenarios. Formally, at iteration $T$ of the algorithm, we solve 
\begin{equation*} 
	\min _{\substack{\boldsymbol{z} \in \mathcal{Z}\\ \eta_r \in \mathbb{R}, r\in \mathcal{R}}} \langle \bs{c}, \boldsymbol{z} \rangle + \dfrac{1}{|\mathcal{R}|}\sum_{r \in \mathcal{R}} \eta_r \mbox{ s.t. } \eta_r \geq f(\bs{z}^t; \bs{d}^{\cdot,r}) + \left\langle \nabla f(\bs{z}^t; \bs{d}^{\cdot,r}), \bs{z} - \bs{z}^t\right\rangle, \, \forall t \in [T], \, \forall r \in \mathcal{R}_t,
\end{equation*}
instead of \eqref{eqn:cp.multi}, as sketched in the multi-tree case in Algorithm \ref{alg:scp-multicut} (we defer a detailed discussion of its single-tree implementation and termination criteria to Section \ref{section:bounds}). Consequently, each iteration only requires solving $|\mathcal{R}_{T+1}|$ optimization problems that define $f(\bs{z}^{T+1};\bs{d}^{\cdot,r})$, 
which can be significantly faster. Moreover, it is not too hard to see that this algorithm converges almost surely under any reasonable sampling scheme (e.g., sampling subsets of $\mathcal{R}$ of fixed cardinality uniformly) since we almost surely sample each subset $\mathcal{R}_t$ infinitely often and there are finitely many binaries. Note that, in the pseudo-code, Algorithm \ref{alg:scp-multicut} is initialized with a set of valid constraints generated from scenarios $r \in \mathcal{R}_0$. However, in practice, these constraints do not have to be generated at $\bs{z}_0$, nor be binding. We can initialize the algorithm with any set of valid (linear) constraints on $(\bs{z}, \bs{\eta})$.
\begin{algorithm}[!ht]
	\caption{A Multi-Cut Sample-Based Cutting Plane Method}\label{alg:scp-multicut}
	\begin{algorithmic}[1]
		\State  $\textbf{initialize}\ \bs{z}_0; \, f(\bs{z}_{0}; \bs{d}^{\cdot,r}), \nabla f(\bs{z}_{0}; \bs{d}^{\cdot,r}), \forall r \in \mathcal{{R}}_0$.
		\State  $\textbf{set}\ T \leftarrow 0$ 
		\Repeat
		\State $\begin{array}{crl} \textbf{compute}\ \bs{z}^{T+1}, \boldsymbol{\eta}^{T+1} \leftarrow 
			&  \displaystyle \underset{\bs{z, \eta_r}}{\argmin} & \langle \bs{c}, \bs{z}\rangle + \dfrac{1}{|\mathcal{R}|}\mathlarger{\sum_{r \in \mathcal{R}}} \eta_{r} \\[1em]
			& \text{s.t.} & \eta_{r} \geqslant f(\bs{z}^t; \bs{d}^{\cdot,r})+ \langle \nabla f(\bs{z}^t; \bs{d}^{\cdot,r}), \bs{z}-\bs{z}^{t}\rangle, \forall t \in [T], \forall r \in \mathcal{R}_{t},  
		\end{array}$
		\State $\textbf{sample}\ \mathcal{{R}}_{T+1} \subseteq \mathcal{R}$
		\State $\textbf{calculate}\ f(\bs{z}^{T+1}; \bs{d}^{\cdot,r}), \nabla f(\bs{z}^{T+1}; \bs{d}^{\cdot,r})$ for $r \in \mathcal{R}_{T+1}$
		\State $\textbf{set}\ T \leftarrow T+1$
		\Until{Termination Criterion Met}
	\end{algorithmic}
\end{algorithm}

On the other hand, developing a stochastic version of the single-cut method is technically challenging because constraint \eqref{eqn:cp.single.constraint} aggregates information across scenarios. To address this issue, \cite{infanger_monte_1992} 
propose generating probabilistic cuts by sampling a subset of scenarios $\mathcal{R}_t \subseteq \mathcal{R}$ at each iteration and imposing the constraint
\begin{equation}\label{eqn:scp.single.bertli}
	\frac{|\mathcal{R}_t|}{|\mathcal{R}|} \times \eta \geq \sum_{r \in \mathcal{R}_t} f(\bs{z}^t; \bs{d}^{\cdot,r}) + \left\langle \sum_{r \in \mathcal{R}_t} \nabla f(\bs{z}^t; \bs{d}^{\cdot,r}), \bs{z} - \bs{z}^t\right\rangle,
\end{equation}
instead of \eqref{eqn:cp.single.constraint}, where the quantities $\tfrac{|\mathcal{R}|}{|\mathcal{R}_t|} \sum_{r \in \mathcal{R}_t} f(\bs{z}^t; \bs{d}^{\cdot,r})$ and $\tfrac{|\mathcal{R}|}{|\mathcal{R}_t|} \sum_{r \in \mathcal{R}_t} \nabla f(\bs{z}^t; \bs{d}^{\cdot,r})$ are unbiased estimates of the original offset and slope terms, $\sum_{r \in \mathcal{R}} f(\bs{z}^t; \bs{d}^{\cdot,r})$ and $\sum_{r \in \mathcal{R}} \nabla f(\bs{z}^t; \bs{d}^{\cdot,r})$ respectively, so that \eqref{eqn:scp.single.bertli} is a reasonable approximation of the original constraint \eqref{eqn:cp.single.constraint}. This intuition is similar to that of SGD in unconstrained continuous optimization. Unfortunately, these cuts are only valid probabilistically and may cut off part of the feasible region when combined. Moreover, while the sampled cuts are unbiased estimates of the slope, optimizing these estimates via Benders decomposition yields solutions that suffer from 
the so-called optimizer's curse \citep{smith2006optimizer}. SGD shares the same drawbacks but mitigates them by performing only one gradient step at each iteration and forgetting estimation errors between iterations. Conversely, in a cutting-plane algorithm, cuts added at one iteration are imposed in subsequent iterations, until termination. 

We reconcile the computational benefits of sampling with the aforementioned drawbacks of the stochastic single-cut approach by leveraging the dual formulation of $f(\bs{z};\bs{d})$ in Proposition \ref{prop:f.dual} to derive deterministically valid lower-approximations for scenarios $r$ that are not sampled. Further, we argue that provided the sampled scenarios are sufficiently representative of the remaining scenarios, this approximation is sufficiently accurate that we eventually obtain a near-optimal solution with high probability; see also \cite{zakeri_inexact_2000} for an ``inexact'' Benders decomposition method. 

Specifically, recall that any feasible dual solution $(\bs{\alpha}, \bs{\beta}, \bs{p})$ provides a valid lower bound:
\begin{align*}
f(\bs{z}; \bs{d}^{\cdot,r}) \geq q(\bs{z}^t, \bs{\alpha}, \bs{\beta}, \bs{p}; \bs{d}^{\cdot,r}) + \left\langle \nabla_{\bs{z}} q(\bs{z}^t, \bs{\alpha}, \bs{\beta}, \bs{p}; \bs{d}^{\cdot,r}), \bs{z} - \bs{z}^t \right\rangle,
\end{align*}
with $q(\bs{z}, \bs{\alpha}, \bs{\beta}, \bs{p}; \bs{d}) := \sum_{k \in \mathcal{K}} \langle \bs{p}^{k}, \bs{d}^{k} \rangle 
- \sum_{{(i,j)} \in\mathcal{E} } 
		 z_{i,j}  \left[ u_{i,j} \beta_{i,j} + \frac{\gamma}{2} (\alpha_{i,j} + \beta_{i,j})^2 \right] $. 
Hence, we replace \eqref{eqn:cp.single.constraint} by a constraint of the form
\begin{align} \label{eqn:scp.single.dual}
\eta \geq \sum_{r \in \mathcal{R}} q(\bs{z}^t, \bs{\alpha}^r, \bs{\beta}^r, \bs{p}^r; \bs{d}^{\cdot,r}) + \sum_{r \in \mathcal{R}} \left\langle \nabla_{\bs{z}} q(\bs{z}^t, \bs{\alpha}^r, \bs{\beta}^r, \bs{p}^r; \bs{d}^{\cdot,r}), \bs{z} - \bs{z}^t \right\rangle,
\end{align}
for some feasible dual solutions $(\bs{\alpha}^r, \bs{\beta}^r, \bs{p}^r)$. Observe that, unlike \eqref{eqn:scp.single.bertli}, the constraint \eqref{eqn:scp.single.dual} is a deterministically valid (although not necessarily tight) lower bound on the true operational cost. 

Collecting these observations yields our overall stochastic single-cut approach: First, to reduce the computational burden of solving an optimization problem for each scenario, at each iteration, we only solve a random subset of scenarios $r \in \mathcal{R}_t \subseteq \mathcal{R}$ --hence effectively computing $f(\bs{z}^t; \bs{d}^{\cdot,r})$ and $\nabla f(\bs{z}^t; \bs{d}^{\cdot,r})$. Second, for the remaining scenarios $r \notin \mathcal{R}_t$, we refrain from solving \eqref{eqn:f.dual} and instead use the cheap to compute and feasible dual average solution
$(\bs{\alpha}^r, \bs{\beta}^r,\bs{p}^r) = (\bar{\bs{\alpha}}^{\mathcal{R}_t}, \bar{\bs{\beta}}^{\mathcal{R}_t}, \bar{\bs{p}}^{\mathcal{R}_t}) := \dfrac{1}{|\mathcal{R}_t|} \sum_{r' \in \mathcal{R}_t} (\bs{\alpha}^{r'}, \bs{\beta}^{r'},\bs{p}^{r'})$ instead. This gives a stochastic cutting-plane method with a sequence of deterministically valid non-decreasing lower bounds, which we formalize in Algorithm \ref{alg:scp-singlecut} (we defer a detailed discussion of its single-tree implementation and termination criterion to Section \ref{section:bounds}). 

However, whether this method converges towards an optimal solution (e.g., in a limit) or generates a never-ending sequence of deterministically valid but not tight cuts is not obvious. We now provide some reassurance in this direction, by showing that for the incumbent solution $\bs{z}^t$, the approximation error 
of cuts obtained via dual averaging can be decomposed, with high probability, as the sum of two terms: one term that depends on the variance of the optimal dual variables and that captures the heterogeneity in the demand scenarios, and one estimation error term that vanishes as $|\mathcal{R}_t|$ grows (proof deferred to Appendix \ref{ssec:a.proof.approx}):
\begin{proposition}\label{prop:avg.approx} Fix $\bs{z}^t$. For any $r \in \mathcal{R}$, denote  $(\bs{\alpha}^r, \bs{\beta}^r, \bs{p}^r)$ the optimal dual solutions of \eqref{eqn:f.dual} for $\bs{z} = \bs{z}^t$ and $\bs{d} = \bs{d}^{\cdot,r}$. Denote $\nu^2$ the variance in optimal dual variables, defined as 
\begin{align*}
    \nu^2 = \dfrac{1}{|\mathcal{R}|} \sum_{r \in \mathcal{R}} \left\| (\bs{\alpha}^r, \bs{\beta}^r, \bs{p}^r) - (\bar{\bs{\alpha}}^{\mathcal{R}}, \bar{\bs{\beta}}^{\mathcal{R}}, \bar{\bs{p}}^{\mathcal{R}}) \right\|^2 \mbox{ with }  (\bar{\bs{\alpha}}^{\mathcal{R}}, \bar{\bs{\beta}}^{\mathcal{R}}, \bar{\bs{p}}^{\mathcal{R}}) = \dfrac{1}{|\mathcal{R}|} \sum_{r \in \mathcal{R}} (\bs{\alpha}^r, \bs{\beta}^r, \bs{p}^r).
\end{align*}
Then, there exist universal constants $L, M > 0$ 
such that, for any $\delta \in (0,e^{-1})$, when $\mathcal{R}_t$ is sampled without replacement from $\mathcal{R}$ with a fixed size $|\mathcal{R}_t|$, we have with probability $1-3\delta$:
\begin{align}
\sum_{r \notin \mathcal{R}_t} \left|  q(\bs{z}^t, \bar{\bs{\alpha}}^{\mathcal{R}_t}, \bar{\bs{\beta}}^{\mathcal{R}_t}, \bar{\bs{p}}^{\mathcal{R}_t}; \bs{d}^r) - f(\bs{z}^t; \bs{d}^r) \right| 
&\leq  
L \sqrt{|\mathcal{R} \setminus \mathcal{R}_t |} \nu 
+ D \sqrt{|\mathcal{R} \setminus \mathcal{R}_t |\, \log(1/\delta)},
\end{align}
with 
\begin{align*}
    D := L \, M \sqrt{2|\mathcal{E}| + |\mathcal{N}| \times |\mathcal{K}|} \left[ \sqrt{|\mathcal{R}|} \left(\dfrac{1}{|\mathcal{R}_t|} - \dfrac{1}{|\mathcal{R}|} \right)^{1/2} + \left(\dfrac{1}{|\mathcal{R} \setminus \mathcal{R}_t|} - \dfrac{1}{|\mathcal{R}|} \right)^{1/4} \right].
\end{align*}
\end{proposition}

\begin{algorithm}[!ht]
	\caption{A Single-Cut Sample-Based Cutting Plane Method}\label{alg:scp-singlecut}
	\begin{algorithmic}[1]
		\State  $\textbf{initialize}\ \bs{z}_1;\, f(\bs{z}_{0}; \bs{d}^{\cdot,r}), \nabla f(\bs{z}_{0}; \bs{d}^{\cdot,r}), \forall r \in \mathcal{{R}}_0$.
		\State  $\textbf{set}\ T \leftarrow 1$ 
		\Repeat
		\State $\begin{array}{crl} \textbf{compute}\ \bs{z}^{T+1}, {\eta}^{T+1} \leftarrow 
			&  \displaystyle \underset{\bs{z}, \eta}{\argmin} & \langle \bs{c}, \bs{z}\rangle + \tfrac{1}{|\mathcal{R}|}\, \eta \\
			& \text{s.t.} & \eta \geqslant \sum_{r \in \mathcal{R}} q(\bs{z}^t; \bs{d}^{\cdot,r})+ \langle \nabla q(\bs{z}^t; \bs{d}^{\cdot,r}), \bs{z}-\bs{z}^{t}\rangle, \forall t \in [T],
		\end{array}$
		\State $\textbf{sample}\ \mathcal{{R}}_{T+1} \subseteq \mathcal{R}$
		\State $\textbf{calculate}\ f(\bs{z}^{T+1}; \bs{d}^{\cdot,r}), \nabla f(\bs{z}^{T+1}; \bs{d}^{\cdot,r})$ for $r \in \mathcal{R}_{T+1}$
		\State $\textbf{set}\ T \leftarrow T+1$
		\Until{Termination Criterion Met}
	\end{algorithmic}
\end{algorithm}

Proposition \ref{prop:avg.approx} provides a probabilistic guarantee on the quality of each cut in terms of the sample size $\vert \mathcal{R}_t\vert$. Observe that the approximation error is proportional to $\sqrt{|\mathcal{R} \setminus \mathcal{R}_t|}$, which means that the approximation error is zero in the limit where $\mathcal{R}_t \rightarrow \mathcal{R}$ (as expected) but which also means that the approximation error grows sub-linearly in the number of scenarios to approximate $|\mathcal{R} \setminus \mathcal{R}_t|$. 

Moreover, by the probabilistic method \citep[see, e.g.,][]{grimmett2020probability}, Proposition \ref{prop:avg.approx} reveals that, for any $\bs{z}^t$ and sufficiently small $\delta$, there exists some $\mathcal{R}_t$ such that this guarantee holds deterministically. Indeed, setting 
$\delta < (1-{{\vert \mathcal{R}\vert} \choose {\vert \mathcal{R}_t\vert}}^{-1})/3$ reveals that, with the notations of Proposition \ref{prop:avg.approx}, repeatedly sampling $\mathcal{R}_t$ for a given $\bs{z}^t$ eventually gives a cut which is an underestimator of $f(\bs{z}^t)$ by at most $\rho$, where
\begin{align}\label{eqn:definerho}
\rho := L \sqrt{|\mathcal{R} \setminus \mathcal{R}_t|} \nu + D \sqrt{|\mathcal{R} \setminus \mathcal{R}_t |\, \log(1/\delta)}.
\end{align}
The above observation implies that running Algorithm \ref{alg:scp-singlecut} without termination and selecting a $\bs{z}^t$, which minimizes our underestimator in the limit, almost surely returns a $\rho$-optimal solution to Problem \eqref{eqn:prob_description}, where $\rho$ is defined by Equation \eqref{eqn:definerho}. Therefore, in practice, when Algorithm \ref{alg:scp-singlecut}'s lower bound stabilizes, we can either increase the number of scenarios sampled (and thus reduce $\rho$), or terminate with confidence if, according to a statistical test, the gap between our stochastic upper bound (see Section \ref{section:bounds}) and our deterministic lower bound is sufficiently small. As we observe in our numerical results (see Section \ref{section:experiments}), the optimality gap from single-cut at termination with a sample rate of around $10\%$ is usually quite small in practice.

We conclude this section with two remarks that contrast our approach with the recent work of \citet{ramirez2023benders} and incorporate dual averaging within our multi-cut method respectively:

\begin{remark} \OneAndAHalfSpacedXI 
Recently, \citep{ramirez2023benders} proposed an adaptive scenario aggregation scheme that applies to stochastic network design problems. Their scheme clusters scenarios into groups, and generates a lower bound on the average cost within each group by the cost associated with the average scenario for that group (via Jensen's inequality). Their approach shares some commonalities with this work, chiefly applying a separation oracle to a subset of scenarios in a stochastic network design problem, rather than all scenarios. However, it differs from our approach in two important aspects: First, \citet{ramirez2023benders} aggregate demand vectors $\bs{d}^r$ within each group, while we aggregate optimal dual variables. Second, their clustering of scenarios into groups (hence, the scenarios passed to the separation oracle) is fixed throughout Benders algorithm and only refined after termination, while we sample a new subset for each incumbent 
explored through the 
algorithm. 
\end{remark}

\begin{remark} \OneAndAHalfSpacedXI 
Although the dual averaging technique is not needed to develop a stochastic multi-cut cutting-plane algorithm, it can be used to improve its convergence. In Algorithm \ref{alg:scp-multicut}, instead of only imposing a new cut for the epigraph variables $\eta_r$ with $r \in \mathcal{R}_t$, we can also use dual averaging to impose one additional constraint on the variables $\eta_r, r \notin \mathcal{R}_t$. 
Formally, 
\begin{align}
\sum_{r \notin \mathcal{R}_t} \eta_r & \geq \sum_{r \notin \mathcal{R}_t} q(\bs{z}^t, \bar{\bs{\alpha}}^{\mathcal{R}_t}, \bar{\bs{\alpha}}^{\mathcal{R}_t}, \bar{\bs{p}}^{\mathcal{R}_t}; \bs{d}^{\cdot,r}) + \sum_{r \notin \mathcal{R}_t} \left\langle \nabla_{\bs{z}} q(\bs{z}^t, \bar{\bs{\alpha}}^{\mathcal{R}_t}, \bar{\bs{\beta}}^{\mathcal{R}_t}, \bar{\bs{p}}^{\mathcal{R}_t}; \bs{d}^{\cdot,r}), \bs{z} - \bs{z}^t \right\rangle,
\end{align}
In our experiments, we 
refer to this implementation as the \emph{accelerated multi-cut} approach.
\end{remark}

\section{Upper Bounds in Stochastic Cutting Planes with Binary Variables}\label{section:bounds}
In this section, we analyze the upper bounds obtained at each iteration of our cutting-plane methods and design convergence criteria that allow us to terminate our methods with confidence. 

The primary motivation for this section is that while the lower bounds for the three stochastic cutting plane methods introduced in Section \ref{section:cp_method} are deterministic, their per iteration estimates of the cost associated with each incumbent solution $\bs{z}^t$, 
\begin{align}\langle \bs{c}, \bs{z}^t\rangle+\dfrac{1}{|\mathcal{R}_t|} \sum_{r \in \mathcal{R}_t} f(\bs{z}^t; \bs{d}^{\cdot,r})\end{align}
are stochastic estimates that depend on the sample $\mathcal{R}_t$. Accordingly, we cannot simply use these stochastic estimates in the same way as in a deterministic method and terminate when the deterministic lower bound, say
$ \displaystyle
\langle \bs{c}, \bs{z}^t\rangle+\tfrac{1}{\vert \mathcal{R}\vert}\sum_{r \in \mathcal{R}}\eta_r
$ in the multi-cut case, is within $\epsilon$ of our stochastic upper bound, or we may terminate because $\bs{z}^t$ is a high variance solution and we picked an optimistic sample set $\mathcal{R}_t$, rather than because $\bs{z}^t$ is an optimal solution; see also \cite{smith2006optimizer}. 

In addition, another salient characteristic of our problem is that the decision variables $\bs{z}$ are binary. Hence, as described in pseudo-code in Algorithm \ref{alg:scp-multicut} and \ref{alg:scp-singlecut}, a MIO problem needs to be solved at each iteration by constructing a branch-and-bound tree (multi-tree implementation). Nowadays, efficient implementations of these schemes exist that simultaneously construct the branch-and-bound tree and generate cutting planes (single-tree implementation). We also discuss the extent to which the stochastic cutting-plane algorithms we developed in the previous section can be implemented with a single-tree instead of multi-tree approach.

\subsection{Convergence Criteria}
In this section, we define a convergence criterion by using an asymptotically normal estimator of the upper bound and using a related upper confidence bound. 
Suppose that one of our stochastic cutting-plane methods finds a solution $\bs{z}$, and that we would like to evaluate its quality. Then, we can use a sample $\mathcal{W}$ to estimate the true cost of this solution 
$$\bar{c}=\langle \bs{c}, \bs{z}\rangle+\dfrac{1}{|\mathcal{R}|} \sum_{r \in \mathcal{R}} f(\bs{z}; \bs{d}^{\cdot,r})$$
by its estimate on $\mathcal{W}$:
$$\hat{c}^{\mathcal{W}}=\langle \bs{c}, \bs{z}\rangle+\dfrac{1}{|\mathcal{W}|} \sum_{r \in \mathcal{W}} f(\bs{z}; \bs{d}^{\cdot,r}).$$
In this section, for simplicity, we omit the dependency of $\hat{c}^{\mathcal{W}}$, $\bar{c}$, and the following quantities, on the solution $\bs{z}$. We also denote $\mathcal{W}$ the random sample used for termination since it could be a new independent draw from the sample $\mathcal{R}_t$ used in the algorithm (and should be, for our estimation procedure to be unbiased). 

As noted by \cite{morton1998stopping, mak1999monte}, under some mild assumptions on the distribution of $\bs{d}^k$ (e.g., finite variance), for an infinite number of scenarios $|\mathcal{R}|$, this estimator obeys a central limit theorem:
$$\sqrt{\vert \mathcal{W}\vert} \left[\hat{c}^{\mathcal{W}}-\bar{c}\right] \xrightarrow{d} \mathcal{N}(0, \sigma^2_c) \: \mbox{ as } \: \vert \mathcal{W}\vert \rightarrow \infty,$$
where $\sigma^2_c=\mathrm{Var}(f(\bs{z}, d^{,r}))$ can be estimated via the sample variance estimator
\begin{align*}
\hat{\sigma}_c^2:=\frac{1}{\vert \mathcal{W}\vert-1} \sum_{r \in \mathcal{W}}\left(f(\bs{z}, \bs{d}^{,r})-\tfrac{1}{|\mathcal{W}|} \sum_{s \in \mathcal{W}} f(\bs{z}; \bs{d}^{\cdot,s})\right)^2.
\end{align*}
In reality, however, we only have finitely many observations $\mathcal{R}$. Yet, provided $\vert \mathcal{R}\vert $ is large relative to $\vert \mathcal{W}\vert$, we can still apply the CLT to estimate the cost of $\bs{z}$. Consequently, letting ${q}_\alpha$ be such that $\mathbb{P}(\mathcal{N}(0,1) \leq {q}_\alpha)=1-\alpha$, we can construct an asymptotically valid confidence interval for this estimator at level $\alpha$ of the form 
\begin{align*}
\left[ \hat{c}^{\mathcal{W}}-\frac{{q}_{\alpha/2}}{\sqrt{\vert \mathcal{W}\vert}} \hat{\sigma}_c, \hat{c}^{\mathcal{W}}+\frac{{q}_{\alpha/2}}{\sqrt{\vert \mathcal{W}\vert}} \hat{\sigma}_c \right].
\end{align*}

We terminate our method using a modified version of the convergence criteria proposed by \citet{morton1998stopping}. Namely, letting $$\bar{c}_{\alpha,t}:=\hat{c}^{\mathcal{W}}+\frac{{q}_{\alpha/2}}{\sqrt{\vert \mathcal{W}\vert}} \hat{\sigma}$$ denote {an} upper confidence bound at level $\alpha$ on the cost of $\bs{z}^t$, the solution generated at the $t$th iterate of one of our cutting-plane methods,  
we terminate as soon the conservative bound gap falls below a predefined threshold $\epsilon$,
i.e., for the multi-cut method
\begin{align}
\dfrac {\bar{c}_{\alpha,t} - \left( \langle \bs{c}, \bs{z}^t\rangle + \frac{1}{\vert \mathcal{R}\vert} \sum_{r \in \mathcal{R}}\eta_{r,t} \right) }{\bar{c}_{\alpha,t}}  \leq \epsilon \label{eqn:gap.sto}
\end{align}
and the termination criteria for the two remaining methods are similar. Alternatively, we terminate if we exceed a time limit, as discussed in our numerical results. In the latter case, we evaluate the true cost of $\bs{z}^t$ by computing its cost across each scenario in $\mathcal{R}$.

We remark that for some adversarial instances of Problem \eqref{eqn:prob_description}, using the same sample size at each iteration in conjunction with this termination criterion could lead to unattractive results where we terminate at a highly suboptimal solution with high probability \citep[c.f.][Example 1]{morton1998stopping}. To address this issue and provide a confidence bound on our overall solution (accounting for multiple testing problems), we can increase the sample size at each iteration of the method in accordance with \citet[][Theorem 2]{morton1998stopping} or use another sampling rule discussed therein \citep[see also][]{bayraksan2011sequential}. 
However, owing to the single-tree implementation of our cutting-plane methods, as discussed in the next section, we do not test every candidate solution we generate when deciding to terminate. Therefore, as we observe in our numerical results, using the same sample size at each iteration is usually adequate. This is particularly true for the single-cut and $k$-cut methods, which, as discussed previously, often generate conservative lower bounds in practice, meaning that we often terminate at a computational time limit. 

Finally, in circumstances where the total number of scenarios is relatively small, we can evaluate the true upper bound directly, rather than a stochastic estimate of the bound. Accordingly, we take this approach whenever the number of scenarios is sufficiently small.

\subsection{Integrating Optimality Cuts Within a Branch-and-Cut Framework}
Once our cut-generation and termination criterion schemes have been designed, they need to be embedded within a branch-and-cut framework to solve Problem \eqref{eqn:prob_description} to certifiable optimality. 
Indeed, in the naive implementation of our algorithms described in the pseudocode of the previous section,
we need to solve a mixed-integer problem at each iteration. 
For further scalability benefits, we can integrate our stochastic cut generation procedure within a state-of-the-art commercial mixed-integer solver (namely, \verb|Gurobi| version $9.1.2$) using \verb|lazy constraint callbacks|, which 
accelerate cutting-plane methods by constructing a single branch-and-bound tree. 
For example, they have been used to implement deterministic cutting-plane algorithms in a highly efficient and relatively standard way; see, e.g., \citet[Section 4]{fischetti_redesigning_2017}. 

Mixed-integer solvers assume that \verb|lazy constraints| are binding at the point they are generated. Accordingly, they do not visit and do not generate \verb|lazy constraints| twice at the same solution. Our stochastic cuts, however, are not binding, they provide a valid yet not necessarily tight lower bound. Therefore, when we implement our method with \verb|lazy constraints|, the MIO solver can terminate with a highly suboptimal solution it deems optimal, because it (mistakenly) assumes the value of the cut generated at $\bs{z}^t$ and evaluated at that point is precisely the cost of $\bs{z}^t$. 
To avoid this issue, we take a hybrid approach between single- and multi-tree branch-and-cut, which, to our knowledge, has not yet been described in the literature. 

Namely, we maintain an outer loop where, at each iteration, we run a single-tree implementation of branch-and-cut with stochastic cutting-planes. We save all the cuts generated and imposed as \verb|lazy constraints| within a separate cut pool during the branch-and-cut algorithm. 
After the branch-and-cut algorithm, we randomly sample a subset of scenarios $\mathcal{W}$ and compute the termination criterion described in the previous section to determine whether the solution returned by the branch-and-cut algorithm is indeed $\epsilon$-optimal with high probability ($\alpha=0.90$). By computing this convergence criterion at each iteration of the outer loop only, we mitigate the issue of multiple hypothesis testing that would arise when testing the quality of a solution at each iteration of Algorithm \ref{alg:scp-multicut}--\ref{alg:scp-singlecut} (inner loop).
If the convergence criterion is met, we terminate the algorithm. Otherwise, we rerun the branch-and-cut algorithm and ensure the MIO solver no longer considers the previously generated \verb|lazy constraints| as binding: We apply the constraints generated in the lazy cut pool as regular linear constraints, purge the lazy cut pool, and rerun the branch-and-cut algorithm. 
In addition to an optimality gap criterion, we terminate the algorithm when the total computational time exceeds a predefined \texttt{TimeLimit}.
We summarize this procedure in Algorithm \ref{alg:scp-multicut2}. We remark that this approach is related to the notion of restarting a single-tree decomposition in a classical deterministic Benders scheme \citep[see, e.g.,][Section 4.4]{fischetti_Benders_2016}. 

\begin{algorithm}[!ht]
	\caption{Outer Loop for Stochastic Branch-and-Cut}\label{alg:scp-multicut2}
	\begin{algorithmic}[1]
		\State  $\textbf{initialize}\ CutPool = \emptyset,\,  t = 0$ 
		\Repeat
			\State Increment $t \leftarrow t+1$
			\State Initialize Algorithm \ref{alg:scp-multicut}/\ref{alg:scp-singlecut} with constraints in $CutPool$.
			\State Run lazy-constraint implementation of Algorithm \ref{alg:scp-multicut}/\ref{alg:scp-singlecut} 
			\State Save all lazy constraints generated in $CutPool$.
            	\State Obtain candidate optimal solution $\bs{z}^t$.
           	\State Obtain valid lower bound from the MIO branch-and-cut solver.
           	\State Sample $\mathcal{W}$ and compute $\bar{c}_{\alpha,t}$
		\Until{ \eqref{eqn:gap.sto} or \texttt{TimeLimit}}
            \State Return $\bs{z}^t$, stochastic upper bound, and deterministic lower bound
	\end{algorithmic}
\end{algorithm}

Finally, in addition to the hybrid scheme described in this section, one could also consider a pure multi-tree implementation of our stochastic cutting-plane methods, as suggested in Section \ref{section:cp_method} and the classical network design literature \citep[]{geoffrion1974multicommodity}. However, in preliminary numerical experiments, we found that such an approach is significantly slower because it involves solving a different MIO to generate each cut. Accordingly, we do not consider such an approach as part of our numerical experiments.

\subsection{Accelerating the Convergence of our Approach} \label{ssec:bounds.acceleration}
We now describe practical enhancements to our stochastic cutting-plane approaches that improve their convergence, sometimes substantially; see also \cite{fischetti_Benders_2016,fischetti_redesigning_2017, bertsimas_unified_2021} for related discussions on accelerating the convergence of decomposition schemes. To facilitate 
a fair comparison, we implement these strategies 
for all Benders-type 
methods in our experiments.

\paragraph{Warm-Starting the Lower Bound: Cuts at the Root Node}
First, we can warm-start our lower bound by applying cutting planes at the root node obtained after solving a Boolean relaxation of \eqref{eqn:prob_outer_approx} using a continuous analog of our discrete cutting-plane method. 
This strategy is referred to a two-phase Benders approach \citep{mcdaniel1977modified} and has been successfully been applied in network design \citep{crainic2016partial} and other contexts \citep[e.g.,][]{fischetti_redesigning_2017, bertsimas_unified_2021}.
Note that the continuous cutting plane algorithm can also be implemented in a multi- or single-cut fashion and in a deterministic or stochastic version. To balance the tightness of the formulation at the root node against the overall computation cost, we impose a hard constraint on the total number of root node cuts applied (typically $10$ or $20$).

\paragraph{Warm-Starting the Upper Bound: }
We supply the initial network (without any new construction) and the network obtained by constructing all the edges as warm-starts.
However, we do not implement a more sophisticated warm-starting strategy for any of our methods, to better isolate the numerical benefit of our decomposition schemes.
We remark that, in practice, the Boolean relaxation could be randomly rounded to generate provably high-quality feasible solutions \citep[c.f.][Section 3.2]{bertsimas_unified_2021}, and other heuristics specific to network design, as reviewed in the introduction, could also be applied.

\paragraph{Warm-Starting Feasibility Constraints: } Problem-specific inequalities can be added to provide more structure to the master problem \citep{rahmaniani2018accelerating}. Based on the numerical evidence of \citet{rahmaniani2018accelerating}, we implement two types of valid inequalities (origin and destination node inequalities and network connectivity cuts). We also implement partial optimality cuts: Namely, when the incumbent solution $\bs{z}^T$ is infeasible, we not only impose feasibility cuts of the form \eqref{eqn:feasibility.cut} for scenarios $r \in \mathcal{R}_{\texttt{inf}}$, we also derive optimality cuts for scenarios $r \in \mathcal{R} \setminus \mathcal{R}_{\texttt{inf}}$ (or $\mathcal{R}_T \setminus \mathcal{R}_{\texttt{inf}}$ in the stochastic version). In a multi-cut implementation, we can impose these constraints for each $\eta_r$. In the single-cut implementation, we use our dual averaging technique to derive a valid linear inequality on the single epigraph variable $\eta$.

We remark that, on preliminary experiments (Table \ref{tab:overall.deterministic.variants}), we found no clear benefit from using Pareto-optimal cuts \citep{magnanti1981accelerating}. Accordingly, we did not consider them in our implementation. Note that this finding is consistent with prior literature on Pareto-optimal cuts, which finds that they do more harm than good \citep[c.f.][]{papadakos2008practical, fischetti_redesigning_2017}

\section{Numerical Experiments}\label{section:experiments}
In this section, we numerically benchmark our stochastic 
Benders decomposition schemes on data-driven MCFND problems. We also compare their performance with their deterministic counterparts and \verb|Gurobi| on a perspective reformulation of the original MIO formulation \eqref{eqn:prob_description}. 

\subsection{Implementation Details}\label{ssec:exp.implementation}
All experiments were conducted on MIT’s Supercloud Cluster \citep{reuther2018interactive}, which hosts Intel Xeon Platinum 8260 processors. All algorithms were implemented in Julia v1.7.3 \citep{bezanson2017julia} using JuMP v0.21.10 \citep{DunningHuchetteLubin2017Jump} and Gurobi {v9.5.1} \citep{gurobi}. The RAM allocated
varies from 4GB to 176GB for the largest instances, see \ref{ssec:a.exp.memories} for a detailed breakdown.

In Section \ref{ssec:exp.warmstart} and \ref{ssec:exp.scalability}, we consider synthetic instances generated according to a methodology from \citet{gunluk_perspective_2009} and \citet{bertsimas_unified_2021}. In particular, in these instances, the network flow problem for each commodity corresponds to an all-to-one shortest path, and feasibility is not an issue (the pre-existing edges are sufficient to guarantee feasibility).
All in all, we generate instances with varying numbers of nodes $|\mathcal{N}|$, commodities $|\mathcal{K}|$, and scenarios $|\mathcal{{R}}|$, as described in Table \ref{tab:problem_sizes}. We later refer to these instances as small-, medium-, and large-scale instances based on the number of nodes $|\mathcal{N}|$.
\begin{table}
\centering
\caption{Dimensions of the MCFND problems generated, by scale (small-, medium-, and large-scale).}
\label{tab:problem_sizes}
    \begin{tabular}{cccccc}
    \toprule
    Scale & $|\mathcal{N}|$ & & $|\mathcal{K}|$ & & $|\mathcal{R}|$ \\
    \midrule
    Small & \{10,30,50,70\} & $\times$ & \{5,10,25,50\} & $\times$ & \{10,30,50,70,100\} \\
    Medium & \{100,150,200\} & $\times$ & \{5,10,25,50\} & $\times$ & \{10,30,50,70,100\} \\
    Large & \{300,500,700\} & $\times$ & \{5,10,25,50\} & $\times$ & \{10,30,50,70,100\} \\
    \bottomrule
    \end{tabular}
\end{table}
In Section \ref{ssec:exp.rinstances}, we evaluate the experiments on the {\bf R} instances from \citet{crainic2016partial}, with demand scenarios generated by \citet{rahmaniani2018accelerating}.
Details on generating the synthetic and {\bf R} instances are provided in Appendix \ref{ssec:a.instance}.

For our algorithms, we use two termination criteria: a time limit (7,200 seconds) and an optimality gap target $\epsilon = 1\%$ (with $\alpha = 0.90$  for our stochastic algorithms). Note that the time limit applies to the full outer-loop presented in Algorithm \ref{alg:scp-multicut2} (and not on each run of the branch-and-cut algorithm only). 
{For all stochastic methods, we use a sampling rate, $|\mathcal{R}_t| / |\mathcal{R}|$, of 10\%.}  
We also fix the regularization parameter $\gamma$ to $1$ ---we discuss its impact on our algorithms in Appendix \ref{sec:a.gamma}. We warm-start all methods with the original connected graph as an initial solution.

\subsection{Comparison of Different Stochastic Cutting-Plane Algorithms}\label{ssec:exp.warmstart}
In this section, we benchmark the variants of the stochastic cutting-plane algorithm proposed in Section \ref{section:cp_method}, namely the multi-, single-, accelerated multi-cut, 
in terms of their ability to obtain a certifiably near-optimal solution with high confidence. We also measure the impact of warm-starting these methods with cuts obtained from solving the perspective relaxation with a multi- or single-cut stochastic cutting-plane algorithm, and applying these cuts at the root node in our branch-and-cut scheme (which we refer to as multi-cut or single-cut root node cuts respectively). 
We report average computational time (capped at  7,200 seconds) for solving our small and medium-scale synthetic instances in Table \ref{tab:root_cuts_runtime}. To augment these results, Table \ref{tab:root_cuts_gap} reports the average optimality gap at termination, and Table \ref{tab:root_cuts_conv} (see Appendix \ref{ssec:a.exp.warmstart}) reports the fraction of instances solved within the time limit. Note that the optimality gaps reported in Table \ref{tab:root_cuts_gap} are computed using the true cost of the incumbent solution, using all scenarios in $\mathcal{R}$, and that the time required to calculate this true cost is not included in the computational time of any cutting-plane algorithm.

\begin{table}
	\centering
	\caption{Computational time (in seconds) of the multi-, single-, and accelerated multi-cut stochastic cutting plane algorithm, with different warm-start strategies at the root node (none, multi-cut, and single-cut root node cuts). Metrics are averaged across instances with the same number of nodes $|\mathcal{N}|$.}
	\label{tab:root_cuts_runtime}
    \begin{tabular}{l|ccc|ccc|ccc}
        \toprule
         & \multicolumn{3}{c}{Multi-Cut}& \multicolumn{3}{c}{Single-Cut} & \multicolumn{3}{c}{Accelerated Multi-Cut}\\
       $|\mathcal{N}|$  & None     &            Multi &            Single &           None &              Multi &                      Single &          None &             Multi &             Single  \\ 
       \midrule
         10             &  384.87  &           107.11  &          89.25   &        142.96  &              76.89 &                     85.27   &       542.76  &            124.92 &             94.08   \\
         30             &  5760.18 &           5763.62 &          5512.78 &        5166.12 &            4018.97 &                     4414.57 &       6262.14 &           5951.10 &             6409.78 \\
         50             &  7200.00 &           6588.54 &          7200.00 &        5154.68 &            4514.63 &                     4426.90 &       6625.89 &           6865.11 &             7200.00 \\
         70             &  6788.22 &           7200.00 &          7200.00 &        5861.50 &            5248.96 &                     5132.34 &       7200.00 &           7200.00 &             7200.00 \\
         \midrule
        100             &  7105.08 &           7200.00 &          7157.26 &        6117.74 &            5490.28 &                     4775.17 &       7200.00 &           7200.00 &             7103.43 \\
        150             &  7200.00 &           7200.00 &          7200.00 &        6531.33 &            5819.10 &                     5992.10 &       7200.00 &           7200.00 &             6865.95 \\
        200             &  7200.00 &           7200.00 &          7200.00 &        6571.25 &            6430.63 &                     5625.82 &       7200.00 &           7194.36 &             7200.00 \\
        \bottomrule
    \end{tabular}
\end{table}

\begin{table}
	\centering
	\caption{Relative optimality gap (in \%) at termination for multi-, single-, and accelerated multi-cut algorithms, with different warm-start strategies at the root node. Metrics are averaged across instances with same number of nodes $|\mathcal{N}|$.}
	\label{tab:root_cuts_gap}
    \begin{tabular}{l|ccc|ccc|ccc}
        \toprule
         & \multicolumn{3}{c}{Multi-Cut}& \multicolumn{3}{c}{Single-Cut} & \multicolumn{3}{c}{Accelerated Multi-Cut}\\
       $|\mathcal{N}|$  &    None &        Multi &       Single &                None &         Multi &            Single &                   None &                  Multi &         Single \\ 
       \midrule
         10 &               0.26  &        0.02  &        0.93  &               0.10 &           0.23 &            0.23  &                 0.12  &                     0.07 &          0.93  \\
         30 &               10.36 &        4.21  &        3.85  &              14.25 &           5.27 &            4.30  &                 18.48 &                     4.47 &          3.88  \\
         50 &               15.19 &        2.69  &        2.21  &              10.56 &           4.33 &            3.72  &                 26.32 &                     2.37 &          2.32  \\
         70 &               48.23 &        11.78 &        8.17  &              24.78 &          16.58 &            12.12 &                 52.00 &                     7.28 &          8.46  \\
         \midrule
        100 &               50.55 &        8.90  &        4.44  &              29.73 &          18.51 &            11.35 &                 49.75 &                     6.99 &          4.56  \\
        150 &               61.71 &        11.35 &        9.19  &              53.03 &          22.68 &            12.80 &                 58.06 &                     5.51 &          9.93  \\
        200 &               63.49 &        12.49 &        10.43 &              47.34 &          20.70 &            12.80 &                 60.51 &                     7.15 &         12.01  \\
        \bottomrule
    \end{tabular}
\end{table}

We observe that the multi-cut and single-cut warm-start strategies both effectively reduce the relative optimality gap at termination. Indeed, our root node strategies more than halve the optimality gap at termination compared to not applying cuts at the root node. 
For the multi- and single-cut approach, a single-cut strategy at the root node appears to outperform a multi-cut root node strategy in terms of the relative gap at termination.
For the accelerated multi-cut approach, however, both root node strategies are comparable, with a small edge for multi-cut.
All in all, applying a single-cut approach warm-started with a single-cut method at the root node performs best in terms of computational time, while the accelerated multi-cut approach warm-started with a multi-cut method at the root node achieves the lowest average gap at termination. For this reason, we only report results for these two variants in the following two sections.

\begin{figure}
    \centering
    \includegraphics[width=0.95\textwidth]{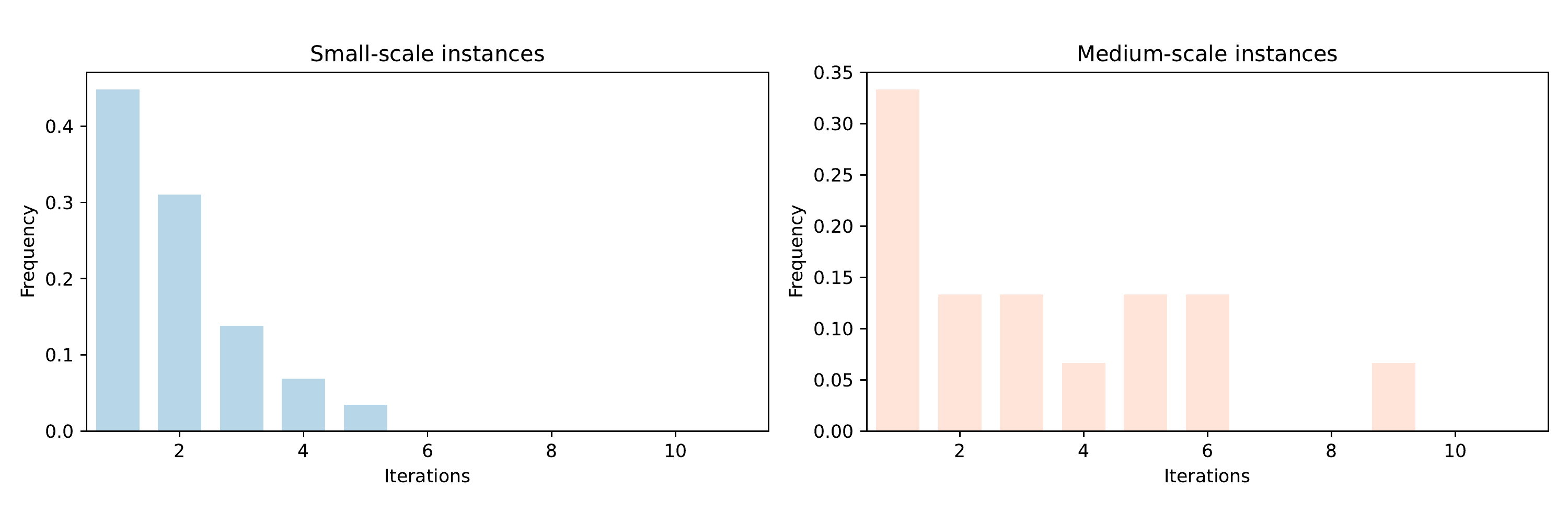}
    \caption{Distribution of the number of outer-loop iterations required by the single-stochastic cutting-plane algorithms with single-cut root node cuts on small-scale (left panel) and medium-scale (right panel) instances; see Table \ref{tab:problem_sizes} for definitions of small and medium-scale instances.}
    \label{fig:outer_loop_iter}
    \scriptsize
\end{figure}

Next, we investigate the number of iterations of the outer loop performed by our methods; recall that in Section \ref{section:bounds}, we propose an outer loop procedure that allows our sampling approach to be safely integrated within a branch-and-cut procedure, without requiring a new branch-and-bound tree each time we generate a cut. To this end, Figure \ref{fig:outer_loop_iter} depicts the number of outer-loop iterations performed by our single-cut algorithm on the small- and medium-scale instances; recall that these instances are defined in Table \ref{tab:problem_sizes}, and comprise instances with $10$--$70$ and $100$--$200$ nodes respectively. We observe that only one outer loop iteration is performed in many cases ($50\%$ for small-scale and $70\%$ for medium-scale instances). In the remaining cases, the first iteration of branch-and-cut with stochastic cuts terminates with a solution that is not $\epsilon$-optimal but 
Algorithm \ref{alg:scp-multicut2} is very efficient, requiring a limited number of additional iterations to identify an optimal solution. This verifies that a single outer loop iteration often wrongly terminates at a solution that is not optimal. On the other hand, only a small number of iterations of the outer loop are usually needed to achieve optimality. Therefore, the tractability of our approach is not compromised by the outer loop. We remind the reader that we impose a global time limit of two hours and Algorithm \ref{alg:scp-multicut2} terminates when it either converges or reaches this time limit. Accordingly, the results on small-scale instances may be less right-censored.

\subsection{Benchmarking Scalability on Synthetic Instances}\label{ssec:exp.scalability}
We now compare the performance of our stochastic cutting plane methods (single- and accelerated multi-cut) against two benchmarks: (a) solving Problem \eqref{eqn:prob_description}'s perspective reformulation directly with \texttt{Gurobi}, (b) a deterministic single-cut method with single-cut root node cuts (we also report the performance of the deterministic method with several acceleration strategies from the literature in Table \ref{tab:overall.deterministic.variants} in Appendix \ref{subsection:a.exp_results}). 
For our stochastic approaches, we use a sampling rate of 10\%. We impose a time limit of 7,200 seconds for all methods. To calibrate our approaches and verify their correctness, we use the smallest instances to verify that all methods terminate with the same optimal solution (see Table \ref{tab:small_scale} in Appendix \ref{subsection:a.exp_results}).

We report the average computational time and optimality gap of all methods, on the small-, medium-, and large-scale instances, in Table \ref{tab:overall_nodes}, with metrics averaged over instances with the same number of nodes $|\mathcal{N}|$. 

\begin{table}
\centering
\caption{Runtime (in seconds) and final optimality gap (in \%) for each algorithm, averaged over instances with the same number of nodes $|\mathcal{N}|$.}
\label{tab:overall_nodes}
\begin{tabular}{lrcrcrcrc}
\toprule
& \multicolumn{2}{c}{\multirow{2}{*}{\texttt{Gurobi} with \eqref{eqn:prob_description}}} & \multicolumn{2}{c}{\multirow{2}{*}{Deterministic}} & \multicolumn{4}{c}{Stochastic} \\
& & & & & \multicolumn{2}{c}{Single} & \multicolumn{2}{c}{Accerated Multi} \\
\cmidrule(lr){2-3} \cmidrule(lr){4-5} \cmidrule(lr){6-7} \cmidrule(lr){8-9} $|\mathcal{N}|$ & Runtime & Gap & Runtime & Gap & Runtime & Gap & Runtime & Gap \\
\midrule
10 & 223.60 & 0.00 & 247.79 & 0.02 & 85.27 & 0.23 & 124.92 & 0.07 \\
30 & 7200.00 & 42.68 & 7163.94 & 6.22 & 4414.57 & 4.30 & 5951.10 & 4.47 \\
50 & 7200.00 & 67.71 & 7200.00 & 4.87 & 4426.90 & 3.72 & 6865.11 & 2.37 \\
70 & 7200.00 & 77.56 & 7200.00 & 11.85 & 5132.34 & 12.12 & 7200.00 & 7.28 \\
\midrule
100 & 7200.00 & 85.15 & 7165.78 & 16.37 & 4775.17 & 11.35 & 7200.00 & 6.99 \\
150 & 7200.00 & 95.97 & 7186.61 & 23.49 & 5992.10 & 12.80 & 7200.00 & 5.51 \\
200 & 7196.63 & 92.87 & 6853.71 & 26.68 & 5625.82 & 12.80 & 7194.36 & 7.15 \\
\midrule
300 & - & - & 6237.87 & 23.11 & 6017.72 & 11.04 & 7200.00 & 11.04 \\
500 & - & - & 6441.49 & 49.09 & 6321.34 & 26.77 & 7200.00 & 21.90 \\
700 & - & - & 6499.08 & 53.39 & 6295.69 & 39.53 & 7200.00 & 30.72 \\
\bottomrule
\end{tabular}
\end{table}

We observe that a perspective reformulation of the original formulation \eqref{eqn:prob_description} cannot be solved by \texttt{Gurobi} with $100$ or more nodes within the time (2 hours) and memory ($>72$GB) limits. Indeed, while this approach converges within minutes for instances with ten nodes, it fails to identify an optimal solution within the two-hour time limit for instances with 20-70 nodes and terminates with large optimality gaps ($> 30\%$) on average. On the other hand, a deterministic Benders decomposition scheme reaches optimality gaps that are an order of magnitude smaller on instances with 20-70 nodes, scales to instances with up to $200$ nodes, but fails to recover a {solution with a meaningful optimality gap within the time limit for larger problems}. 

Our stochastic cutting plane algorithms significantly improve upon their deterministic counterpart. On small- and medium-scale instances, the single-cut stochastic cutting-plane algorithm reduces the average computational time by 40--90\% on the small instances and 20--50\% on the medium ones. A comparison in terms of average computational times might be misleading, however, because of the time limit, and because many of these instances are not solved to $\epsilon$-optimality. Accordingly, we also compare in terms of the optimality gap. We observe that our single-cut stochastic cutting-plane algorithm terminate with gaps half the size of deterministic algorithms on medium to large instances (i.e., around 5\% for the instances with 10--50 nodes, 12\% for the instances with 70--300 nodes, and 30--40\% for the largest instances compared with 5\%, 12--23\%, and 50\% for the deterministic approach respectively). 
Finally, we observe that our accelerated multi-cut stochastic cutting-plane algorithm is generally slower than its single-cut counterpart (probably due to the increased number of epigraph variables) but achieves even lower optimality gaps at termination: less than 5\% for 10--50 nodes, 7--11\% for 70--300 nodes, and 20--30\% for 400-700 nodes.
Figures \ref{fig:max_problems_experiments_variable_nC.slim}-\ref{fig:max_problems_experiments_variable_nC.hybrid} in Appendix \ref{subsection:a.exp_results} display the optimality gap achieved by both our methods for each value of $|\mathcal{N}|$, $|\mathcal{R}|$, and $|\mathcal{K}|$, and shows that our methods are most sensitive to the number of commodities and nodes. 

\subsection{Benchmarking on the {\bf R} Instances }\label{ssec:exp.rinstances}
We now benchmark our methods on network design instances from the literature, the so-called \textbf{R} instances, originally introduced by \cite{crainic_simplex-based_2000} for deterministic network design problems. 
We use a total 54 instances (see details in Appendix \ref{ssec:a.instance}) with 10--20 nodes and 10--50 commodities. So, compared to our synthetic instances, the ratio $|\mathcal{K}| / |\mathcal{N}|$ is higher for these instances.
The number of scenarios generated varies from 10 to {1,000}. 
In line with our results in the previous section, we evaluate the performance of our single-cut and accelerated multi-cut stochastic cutting plane algorithm with a 10\% sampling rate and a time limit of 7,200 seconds for all methods.

Compared with the synthetic instances considered in the previous sections, each commodity in these instances has one origin and one destination. In addition, we start from a network without any edge, which is challenging for decomposition schemes because they do not have access to second-stage variables $\bs{x}^k$ to ensure primal feasibility of $\bs{z}$. 
In this setting, the separation oracle in the (deterministic or stochastic) cutting-plane algorithm generates either a feasibility or an optimality cut, depending on the feasibility of the incumbent solution $\bs{z}^t$, as described in Section \ref{ssec:cp.benders}. We also implemented the strategies presented in Section \ref{ssec:bounds.acceleration}. 

\begin{figure}
    \centering
    \begin{subfigure}{0.3\textwidth}
    \includegraphics[width=\textwidth]{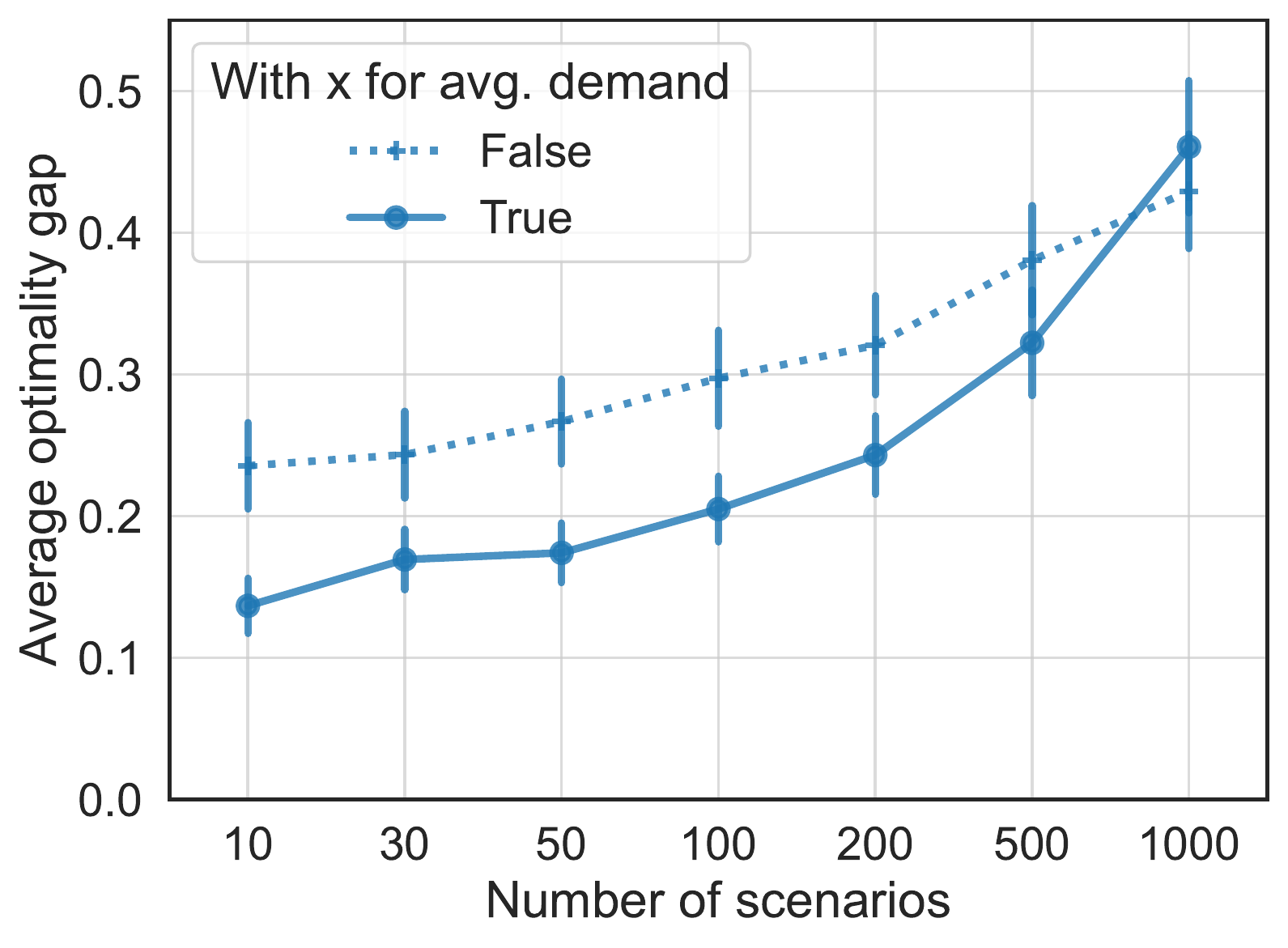}
    \subcaption{Deterministic}
    \end{subfigure}
    \begin{subfigure}{0.3\textwidth}
    \includegraphics[width=\textwidth]{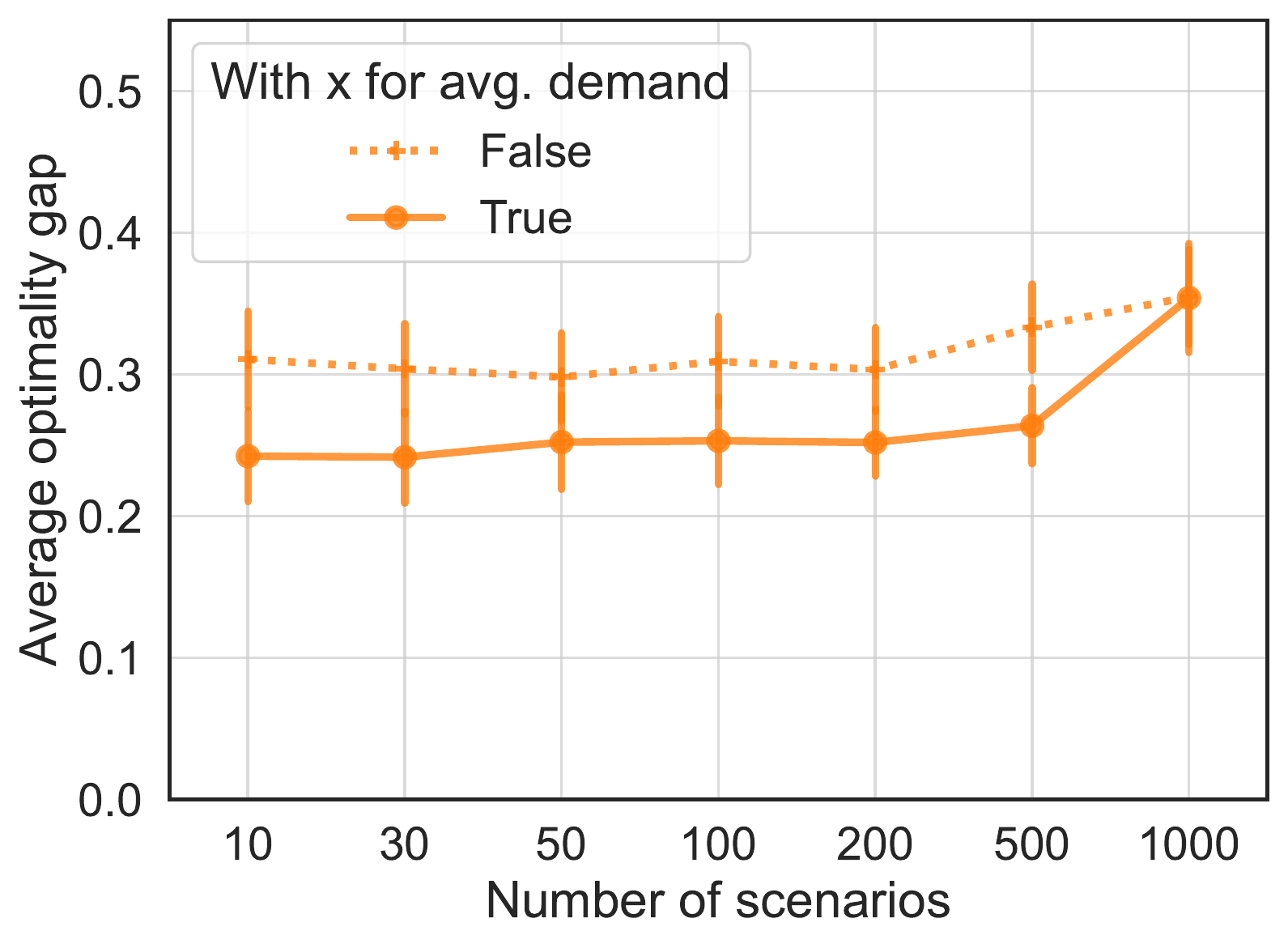}
    \subcaption{Stochastic Single}
    \end{subfigure}
    \begin{subfigure}{0.3\textwidth}
    \includegraphics[width=\textwidth]{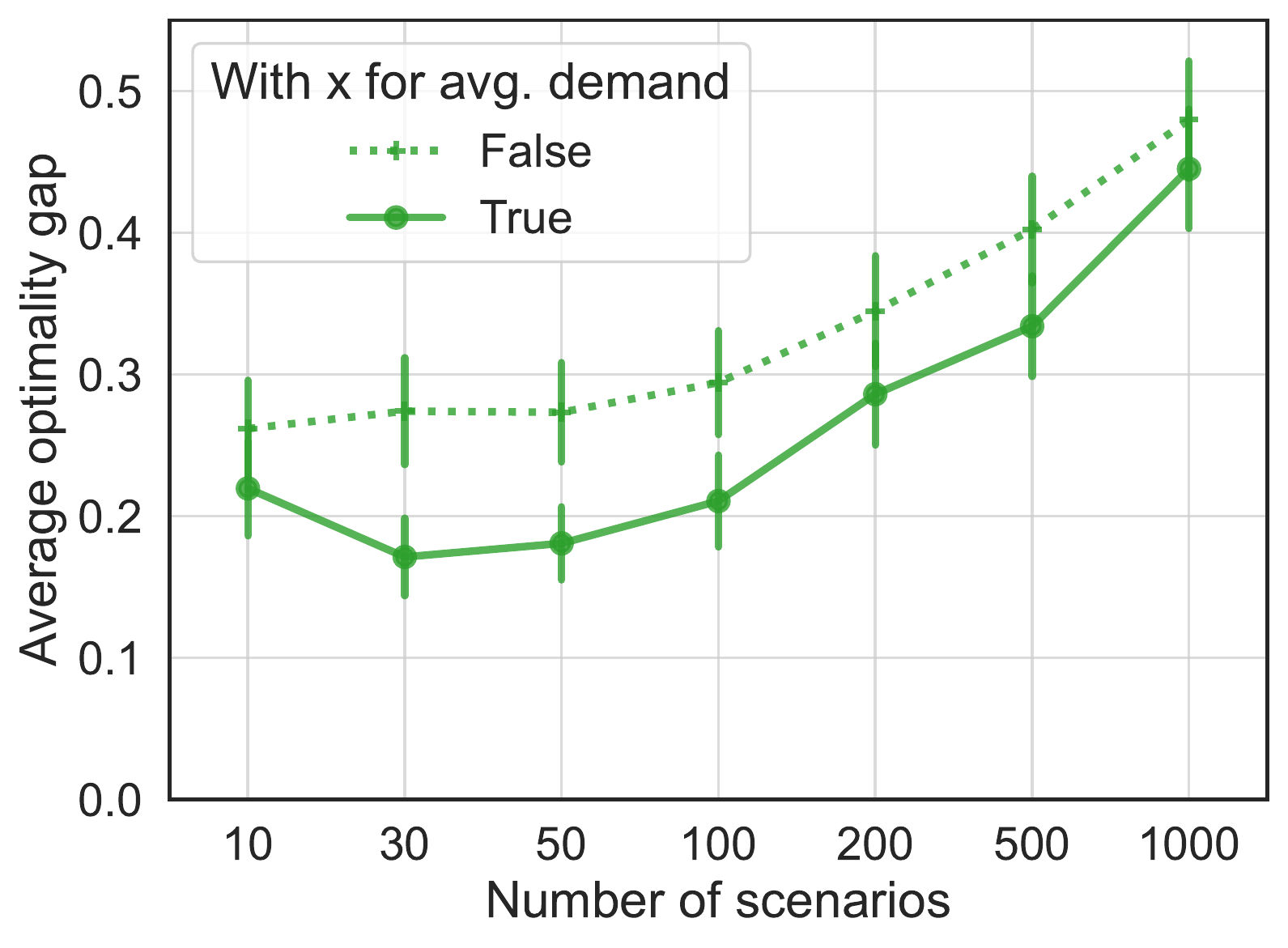}
    \subcaption{Stochastic Acc. Multi}
    \end{subfigure}
    \caption{Impact of implementing a cutting-plane algorithm with one second-stage variable $\bs{x}$ in the master problem (to ensure at least feasibility for the average demand) on the average optimality gap achieved on the R instances, as the number of scenarios $|\mathcal{R}|$ increases. Bars represent standard errors.}
    \label{fig:R.benders.xavg}
\end{figure}

Furthermore, a necessary condition for feasibility is to be feasible for one particular scenario. In particular, we consider adding one second-stage variable $\bs{x}$ to the master problem, to enforce feasibility with respect to the average demand. As displayed in Figure \ref{fig:R.benders.xavg}, we find that this simple strategy effectively reduces the optimality gap achieved by all decomposition schemes. However, the benefit shrinks as the total number of scenarios increases, which suggests that more than one second-stage variable might be needed to achieve the same gain in these instances.

\begin{figure}
    \centering
    \includegraphics[width=0.6\textwidth]{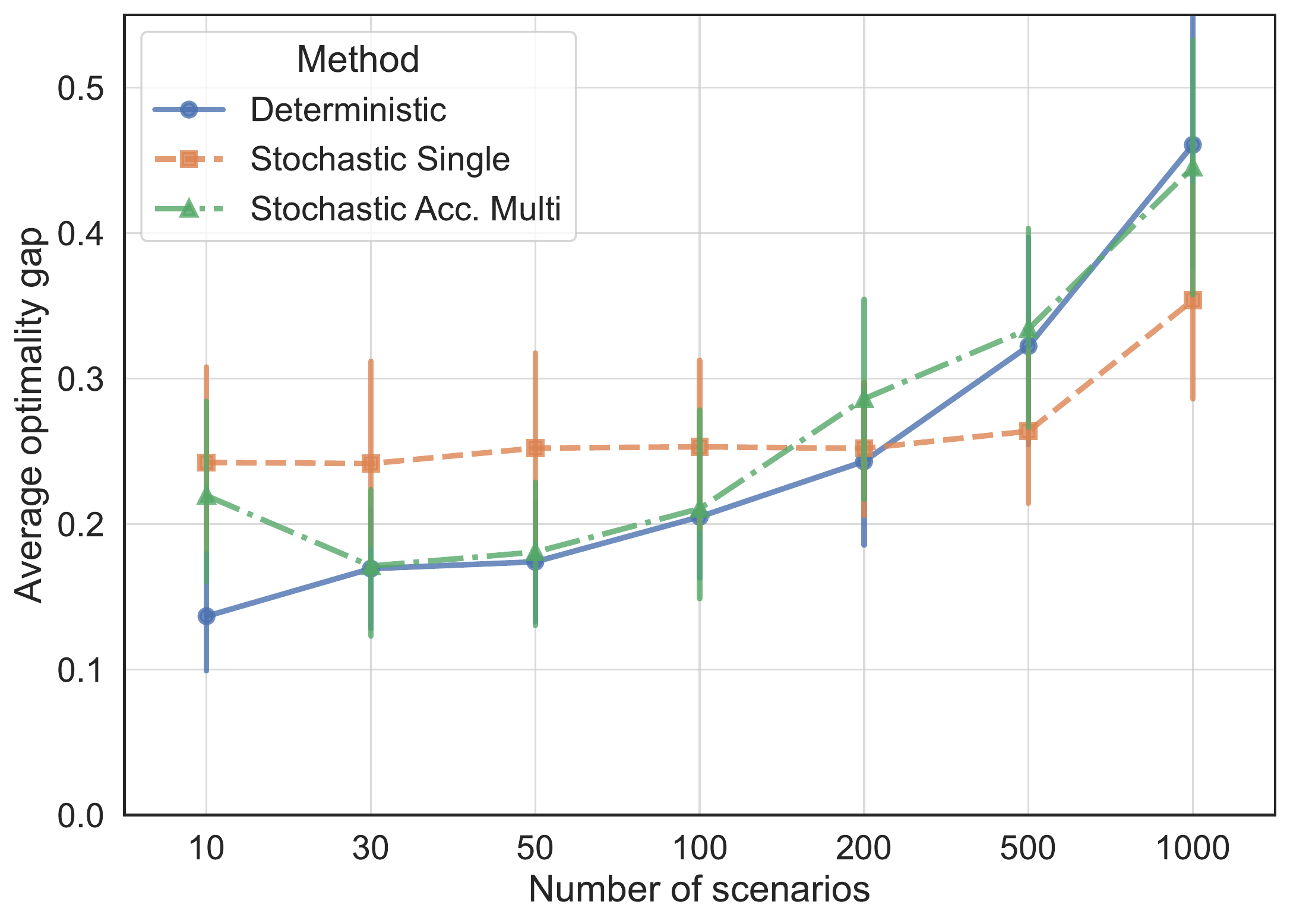}
    \caption{Average optimality gap achieved on the R instances by the deterministic and our stochastic (single-cut and accelerated multi-cut) cutting-plane algorithms, for different number of scenarios $|\mathcal{R}|$. Bars represent standard errors.}
    \label{fig:R.benders.comparison}
\end{figure}

Figure \ref{fig:R.benders.comparison} compares the optimality gap achieved at termination (with a 2-hour time limit) for the three cutting-plane algorithms. As expected, we observe that our stochastic cutting-plane algorithm achieves smaller optimality gaps than its deterministic counterpart as the total number of scenarios increases. Compared to the results on synthetic instances, however, we observe that the deterministic cutting-plane approach outperforms the single-cut stochastic one when the number of scenarios is smaller. We believe this behavior could be explained by the fact that the {\bf R} instances are small ($\mathcal{N}=10$--$20$ nodes), a regime where the different cutting-plane algorithms were achieving comparable optimality gaps in Table \ref{tab:overall_nodes}, and the fact the deterministic implementation (computes and) imposes more constraints per iteration, which might be more valuable when feasibility constraints are needed, especially in the first iterations of the algorithm.

In Appendix \ref{ssec:a.exp.rinstances}, we report the average gap achieve by the naive formulation solved with Gurobi (Figure \ref{fig:R.gurobiVSbenders}), as well as the distribution of optimality gaps achieved by each method (Figure \ref{fig:R.benders.distribution}).

\section{Conclusion}\label{section:conclusion}

We propose a stochastic Benders decomposition scheme which solves large-scale 
stochastic network design problems. Our approach mitigates the high computational cost of generating each cut by sampling a subset of the data at each iteration, while applying a dual-averaging technique to ensure that the cuts generated remain valid for the original problem. We also propose an outer loop technique to ensure the safe termination of our algorithm when the Benders decomposition scheme is implemented via lazy callbacks. To our knowledge, this is the first work synthesizing sampling with a single-tree approach for generating Benders cuts. We consider multi- and single-cut variants of our algorithm (and $k$-cut in Appendix) and discuss its implementation within a branch-and-cut solver. In numerical experiments, we demonstrate that our stochastic decomposition schemes obtain optimality gaps of $5$--$7\%$ on instances with $100$--$200$ nodes, compared to $16$--$26\%$ for deterministic Benders schemes. Moreover, we obtain bound gaps of around $30\%$ on instances with up to $700$ nodes and $50$ commodities, i.e., problem sizes an order of magnitude larger than any instances addressed by exact methods in the literature. Beyond network design, we believe our approach could be applied to other two-stage stochastic optimization problems addressed via sample average approximations.



{\scriptsize
\bibliographystyle{informs2014} 
\bibliography{Bibliography} 

\begin{thebibliography}{87}
\providecommand{\natexlab}[1]{#1}
\providecommand{\url}[1]{\texttt{#1}}
\providecommand{\urlprefix}{URL }

\bibitem[{Agrawal et~al.(1991)Agrawal, Klein, \protect\BIBand{}
  Ravi}]{agrawal1991trees}
Agrawal A, Klein P, Ravi R (1991) When trees collide: An approximation
  algorithm for the generalized {Steiner} problem on networks.
  \emph{Proceedings of ACM Symposium on Theory of Computing}, 134--144.

\bibitem[{Atamt{\"u}rk \protect\BIBand{}
  G{\"u}nl{\"u}k(2018)}]{atamturk2018note}
Atamt{\"u}rk A, G{\"u}nl{\"u}k O (2018) A note on capacity models for network
  design. \emph{Operations Research Letters} 46(4):414--417.

\bibitem[{Atamt{\"u}rk \protect\BIBand{}
  G{\"u}nl{\"u}k(2021)}]{atamturk2021multicommodity}
Atamt{\"u}rk A, G{\"u}nl{\"u}k O (2021) Multicommodity multifacility network
  design. \emph{Network Design with Applications to Transportation and
  Logistics}, 141--166 (Springer).

\bibitem[{Balakrishnan et~al.(1991)Balakrishnan, Magnanti, Shulman,
  \protect\BIBand{} Wong}]{balakrishnan_models_1991}
Balakrishnan A, Magnanti TL, Shulman A, Wong RT (1991) Models for planning
  capacity expansion in local access telecommunication networks. \emph{Annals
  of Operations Research} 33(4):237--284.

\bibitem[{Bardenet \protect\BIBand{}
  Maillard(2015)}]{bardenet2015concentration}
Bardenet R, Maillard OA (2015) Concentration inequalities for sampling without
  replacement. \emph{Bernoulli} 21(3):1361--1385.

\bibitem[{Barnhart et~al.(2003)Barnhart, Belobaba, \protect\BIBand{}
  Odoni}]{barnhart2003applications}
Barnhart C, Belobaba P, Odoni AR (2003) Applications of operations research in
  the air transport industry. \emph{Transportation Science} 37(4):368--391.

\bibitem[{Bayraksan \protect\BIBand{} Morton(2011)}]{bayraksan2011sequential}
Bayraksan G, Morton DP (2011) A sequential sampling procedure for stochastic
  programming. \emph{Operations Research} 59(4):898--913.

\bibitem[{Beale(1955)}]{beale1955minimizing}
Beale EM (1955) On minimizing a convex function subject to linear inequalities.
  \emph{Journal of the Royal Statistical Society: Series B (Methodological)}
  17(2):173--184.

\bibitem[{Bertsekas(1999)}]{bertsekas1999nonlinear}
Bertsekas DP (1999) \emph{Nonlinear Optimization} (Athena Scientific, Belmont).

\bibitem[{Bertsimas \protect\BIBand{}
  Cory-Wright(2022)}]{bertsimas_scalable_2022}
Bertsimas D, Cory-Wright R (2022) A scalable algorithm for sparse portfolio
  selection. \emph{INFORMS Journal on Computing} 34:1489--1511.

\bibitem[{Bertsimas et~al.(2021)Bertsimas, Cory-Wright, \protect\BIBand{}
  Pauphilet}]{bertsimas_unified_2021}
Bertsimas D, Cory-Wright R, Pauphilet J (2021) A unified approach to
  mixed-integer optimization problems with logical constraints. \emph{SIAM
  Journal on Optimization} 31(3):2340--2367.

\bibitem[{Bertsimas \protect\BIBand{} Li(2022)}]{bertsimas_stochastic_2021}
Bertsimas D, Li ML (2022) Stochastic cutting planes for data-driven
  optimization. \emph{INFORMS Journal on Computing} 34(5):2400--2409.

\bibitem[{Bertsimas \protect\BIBand{} Teo(1998)}]{bertsimas1998valid}
Bertsimas D, Teo CP (1998) From valid inequalities to heuristics: A unified
  view of primal-dual approximation algorithms in covering problems.
  \emph{Operations Research} 46(4):503--514.

\bibitem[{Bertsimas \protect\BIBand{}
  Tsitsiklis(1997)}]{bertsimas1997introduction}
Bertsimas D, Tsitsiklis JN (1997) \emph{Introduction to Linear Optimization},
  volume~6 (Athena Scientific Belmont, MA).

\bibitem[{Bezanson et~al.(2017)Bezanson, Edelman, Karpinski, \protect\BIBand{}
  Shah}]{bezanson2017julia}
Bezanson J, Edelman A, Karpinski S, Shah VB (2017) Julia: A fresh approach to
  numerical computing. \emph{SIAM Review} 59(1):65--98.

\bibitem[{Bienstock et~al.(1998)Bienstock, Chopra, G{\"u}nl{\"u}k,
  \protect\BIBand{} Tsai}]{bienstock1998minimum}
Bienstock D, Chopra S, G{\"u}nl{\"u}k O, Tsai CY (1998) Minimum cost capacity
  installation for multicommodity network flows. \emph{Mathematical
  Programming} 81(2):177--199.

\bibitem[{Binato et~al.(2001)Binato, Pereira, \protect\BIBand{}
  Granville}]{binato_new_2001}
Binato S, Pereira MVF, Granville S (2001) A new {B}enders decomposition
  approach to solve power transmission network design problems. \emph{IEEE
  Transactions on Power Systems} 16(2):235--240.

\bibitem[{Birge \protect\BIBand{} Louveaux(2011)}]{birge_introduction_2011}
Birge JR, Louveaux F (2011) \emph{Introduction to {Stochastic} {Programming}}.
  Springer {Series} in {Operations} {Research} and {Financial} {Engineering}
  (New York, NY: Springer New York).

\bibitem[{Birge \protect\BIBand{} Louveaux(1988)}]{birge_multicut_1988}
Birge JR, Louveaux FV (1988) A multicut algorithm for two-stage stochastic
  linear programs. \emph{European Journal of Operational Research}
  34(3):384--392.

\bibitem[{Bixby(2012)}]{bixby2012brief}
Bixby RE (2012) A brief history of linear and mixed-integer programming
  computation. \emph{Documenta Mathematica} 2012:107--121.

\bibitem[{Boland et~al.(2016)Boland, Fischetti, Monaci, \protect\BIBand{}
  Savelsbergh}]{boland_proximity_2016}
Boland N, Fischetti M, Monaci M, Savelsbergh M (2016) Proximity {B}enders: a
  decomposition heuristic for stochastic programs. \emph{Journal of Heuristics}
  22(2):181--198.

\bibitem[{Boyce et~al.(1973)Boyce, Farhi, \protect\BIBand{}
  Weischedel}]{boyce1973optimal}
Boyce DE, Farhi A, Weischedel R (1973) Optimal network problem: a
  branch-and-bound algorithm. \emph{Environment and Planning A} 5(4):519--533.

\bibitem[{Ceria \protect\BIBand{} Soares(1999)}]{ceria1999convex}
Ceria S, Soares J (1999) Convex programming for disjunctive convex
  optimization. \emph{Mathematical Programming} 86(3):595--614.

\bibitem[{Contreras et~al.(2011)Contreras, Cordeau, \protect\BIBand{}
  Laporte}]{contreras2011benders}
Contreras I, Cordeau JF, Laporte G (2011) {B}enders decomposition for
  large-scale uncapacitated hub location. \emph{Operations Research}
  59(6):1477--1490.

\bibitem[{Cornuejols et~al.(1980)Cornuejols, Nemhauser, \protect\BIBand{}
  Wolsey}]{cornuejols1980canonical}
Cornuejols G, Nemhauser GL, Wolsey LA (1980) A canonical representation of
  simple plant location problems and its applications. \emph{SIAM Journal on
  Algebraic Discrete Methods} 1(3):261--272.

\bibitem[{Costa(2005)}]{costa_survey_2005}
Costa AM (2005) A survey on {B}enders decomposition applied to fixed-charge
  network design problems. \emph{Computers \& Operations Research}
  32(6):1429--1450.

\bibitem[{Crainic et~al.(2000)Crainic, Gendreau, \protect\BIBand{}
  Farvolden}]{crainic_simplex-based_2000}
Crainic TG, Gendreau M, Farvolden JM (2000) A simplex-based {Tabu} search
  method for capacitated network design. \emph{INFORMS Journal on Computing}
  12(3):223--236.

\bibitem[{Crainic et~al.(2021{\natexlab{a}})Crainic, Gendreau,
  \protect\BIBand{} Gendron}]{crainic_network_2021}
Crainic TG, Gendreau M, Gendron B, eds. (2021{\natexlab{a}}) \emph{Network
  Design with Applications to Transportation and Logistics} (Cham: Springer
  International Publishing).

\bibitem[{Crainic et~al.(2021{\natexlab{b}})Crainic, Hewitt, Maggioni,
  \protect\BIBand{} Rei}]{crainic_partial_2021}
Crainic TG, Hewitt M, Maggioni F, Rei W (2021{\natexlab{b}}) Partial {B}enders
  decomposition: general methodology and application to stochastic network
  design. \emph{Transportation Science} 55(2):414--435.

\bibitem[{Crainic et~al.(2016)Crainic, Rei, Hewitt, \protect\BIBand{}
  Maggioni}]{crainic2016partial}
Crainic TG, Rei W, Hewitt M, Maggioni F (2016) \emph{Partial {Benders}
  {D}ecomposition {S}trategies for {T}wo-{S}tage {S}tochastic {I}nteger
  {P}rograms}, volume~37 (CIRRELT).

\bibitem[{Dantzig(1955)}]{dantzig1955linear}
Dantzig GB (1955) Linear programming under uncertainty. \emph{Management
  Science} 1(3-4):197--206.

\bibitem[{Dantzig \protect\BIBand{} Infanger(1993)}]{dantzig_multi-stage_1993}
Dantzig GB, Infanger G (1993) Multi-stage stochastic linear programs for
  portfolio optimization. \emph{Annals of Operations Research} 45(1):59--76.

\bibitem[{Davis et~al.(2020)Davis, Drusvyatskiy, Kakade, \protect\BIBand{}
  Lee}]{davis2020stochastic}
Davis D, Drusvyatskiy D, Kakade S, Lee JD (2020) Stochastic subgradient method
  converges on tame functions. \emph{Foundations of Computational Msathematics}
  20(1):119--154.

\bibitem[{de~Camargo et~al.(2008)de~Camargo, Miranda~Jr, \protect\BIBand{}
  Luna}]{de2008benders}
de~Camargo RS, Miranda~Jr G, Luna HP (2008) Benders decomposition for the
  uncapacitated multiple allocation hub location problem. \emph{Computers \&
  Operations Research} 35(4):1047--1064.

\bibitem[{De~Matos et~al.(2015)De~Matos, Philpott, \protect\BIBand{}
  Finardi}]{de2015improving}
De~Matos VL, Philpott AB, Finardi EC (2015) Improving the performance of
  stochastic dual dynamic programming. \emph{Journal of Computational and
  Applied Mathematics} 290:196--208.

\bibitem[{Dunning et~al.(2017)Dunning, Huchette, \protect\BIBand{}
  Lubin}]{DunningHuchetteLubin2017Jump}
Dunning I, Huchette J, Lubin M (2017) {JuMP}: A modeling language for
  mathematical optimization. \emph{SIAM Review} 59(2):295--320.

\bibitem[{Erickson et~al.(1987)Erickson, Monma, \protect\BIBand{}
  Veinott~Jr}]{erickson1987send}
Erickson RE, Monma CL, Veinott~Jr AF (1987) Send-and-split method for
  minimum-concave-cost network flows. \emph{Mathematics of Operations Research}
  12(4):634--664.

\bibitem[{F{\'a}bi{\'a}n(2000)}]{fabian2000bundle}
F{\'a}bi{\'a}n CI (2000) Bundle-type methods for inexact data. \emph{Central
  European Journal of Operations Research} 8(1):35--55.

\bibitem[{Fischetti et~al.(2016)Fischetti, Ljubić, \protect\BIBand{}
  Sinnl}]{fischetti_Benders_2016}
Fischetti M, Ljubić I, Sinnl M (2016) {Benders} decomposition without
  separability: {A} computational study for capacitated facility location
  problems. \emph{European Journal of Operational Research} 253(3):557--569.

\bibitem[{Fischetti et~al.(2017)Fischetti, Ljubić, \protect\BIBand{}
  Sinnl}]{fischetti_redesigning_2017}
Fischetti M, Ljubić I, Sinnl M (2017) Redesigning {{Benders}} {Decomposition}
  for {Large}-{Scale} {Facility} {Location}. \emph{Management Science}
  63(7):2146--2162.

\bibitem[{Florian et~al.(1976)Florian, Bushell, Ferland, Guerin,
  \protect\BIBand{} Nastansky}]{florian1976engine}
Florian M, Bushell G, Ferland J, Guerin G, Nastansky L (1976) The engine
  scheduling problem in a railway network. \emph{INFOR: Information Systems and
  Operational Research} 14(2):121--138.

\bibitem[{Garey \protect\BIBand{} Johnson(1977)}]{garey1977rectilinear}
Garey MR, Johnson DS (1977) The rectilinear {S}teiner tree problem is
  {NP}-complete. \emph{SIAM Journal on Applied Mathematics} 32(4):826--834.

\bibitem[{Gendron et~al.(1999)Gendron, Crainic, \protect\BIBand{}
  Frangioni}]{gendron1999multicommodity}
Gendron B, Crainic TG, Frangioni A (1999) Multicommodity capacitated network
  design. \emph{Telecommunications Network Planning}, 1--19 (Springer).

\bibitem[{Gendron et~al.(2018)Gendron, Hanafi, \protect\BIBand{}
  Todosijevi{\'c}}]{gendron2018matheuristics}
Gendron B, Hanafi S, Todosijevi{\'c} R (2018) Matheuristics based on iterative
  linear programming and slope scaling for multicommodity capacitated fixed
  charge network design. \emph{European Journal of Operational Research}
  268(1):70--81.

\bibitem[{Geoffrion(1972)}]{geoffrion1972generalized}
Geoffrion AM (1972) Generalized {B}enders decomposition. \emph{Journal of
  Optimization Theory and Applications} 10(4):237--260.

\bibitem[{Geoffrion \protect\BIBand{}
  Graves(1974)}]{geoffrion1974multicommodity}
Geoffrion AM, Graves GW (1974) Multicommodity distribution system design by
  {B}enders decomposition. \emph{Management Science} 20(5):822--844.

\bibitem[{Glover(1975)}]{glover1975improved}
Glover F (1975) Improved linear integer programming formulations of nonlinear
  integer problems. \emph{Management Science} 22(4):455--460.

\bibitem[{Goemans \protect\BIBand{} Bertsimas(1993)}]{goemans1993survivable}
Goemans MX, Bertsimas DJ (1993) Survivable networks, linear programming
  relaxations and the parsimonious property. \emph{Mathematical Programming}
  60(1):145--166.

\bibitem[{Grimmett \protect\BIBand{} Stirzaker(2020)}]{grimmett2020probability}
Grimmett G, Stirzaker D (2020) \emph{Probability and Random Processes} (Oxford
  university press).

\bibitem[{Guigues(2020)}]{guigues2020inexact}
Guigues V (2020) Inexact cuts in stochastic dual dynamic programming.
  \emph{SIAM Journal on Optimization} 30(1):407--438.

\bibitem[{G{\"u}nl{\"u}k(1999)}]{gunluk1999branch}
G{\"u}nl{\"u}k O (1999) A branch-and-cut algorithm for capacitated network
  design problems. \emph{Mathematical Programming} 86(1):17--39.

\bibitem[{G\"unl\"uk \protect\BIBand{}
  Linderoth(2009)}]{gunluk_perspective_2009}
G\"unl\"uk O, Linderoth J (2009) Perspective reformulations of mixed integer
  nonlinear programs with indicator variables. \emph{Mathematical Programming}
  183--205.

\bibitem[{{Gurobi Optimization, LLC}(2022)}]{gurobi}
{Gurobi Optimization, LLC} (2022) {Gurobi Optimizer Reference Manual}.

\bibitem[{Higle \protect\BIBand{} Sen(1991)}]{higle1991stochastic}
Higle JL, Sen S (1991) Stochastic decomposition: An algorithm for two-stage
  linear programs with recourse. \emph{Mathematics of Operations Research}
  16(3):650--669.

\bibitem[{Higle \protect\BIBand{} Sen(1996)}]{higle1996duality}
Higle JL, Sen S (1996) Duality and statistical tests of optimality for two
  stage stochastic programs. \emph{Mathematical Programming} 75(2):257--275.

\bibitem[{Infanger(1992)}]{infanger_monte_1992}
Infanger G (1992) Monte {Carlo} (importance) sampling within a {Benders}
  decomposition algorithm for stochastic linear programs. \emph{Annals of
  Operations Research} 39(1):69--95.

\bibitem[{Ke et~al.(2015)Ke, Ferrara, Radicchi, \protect\BIBand{}
  Flammini}]{ke2015defining}
Ke Q, Ferrara E, Radicchi F, Flammini A (2015) Defining and identifying
  sleeping beauties in science. \emph{Proceedings of the National Academy of
  Sciences} 112(24):7426--7431.

\bibitem[{Kingma \protect\BIBand{} Ba(2014)}]{kingma2014adam}
Kingma DP, Ba J (2014) Adam: A method for stochastic optimization. \emph{arXiv
  preprint arXiv:1412.6980} .

\bibitem[{Magnanti et~al.(1993)Magnanti, Mirchandani, \protect\BIBand{}
  Vachani}]{magnanti1993convex}
Magnanti TL, Mirchandani P, Vachani R (1993) The convex hull of two core
  capacitated network design problems. \emph{Mathematical Programming}
  60(1):233--250.

\bibitem[{Magnanti et~al.(1995)Magnanti, Mirchandani, \protect\BIBand{}
  Vachani}]{magnanti1995modeling}
Magnanti TL, Mirchandani P, Vachani R (1995) Modeling and solving the
  two-facility capacitated network loading problem. \emph{Operations Research}
  43(1):142--157.

\bibitem[{Magnanti \protect\BIBand{} Wong(1981)}]{magnanti1981accelerating}
Magnanti TL, Wong RT (1981) Accelerating {B}enders decomposition: Algorithmic
  enhancement and model selection criteria. \emph{Operations Research}
  29(3):464--484.

\bibitem[{Magnanti \protect\BIBand{} Wong(1984)}]{magnanti_network_1984}
Magnanti TL, Wong RT (1984) Network design and transportation planning: models
  and algorithms. \emph{Transportation Science} 18(1):1--55.

\bibitem[{Mak et~al.(1999)Mak, Morton, \protect\BIBand{} Wood}]{mak1999monte}
Mak WK, Morton DP, Wood RK (1999) Monte carlo bounding techniques for
  determining solution quality in stochastic programs. \emph{Operations
  Research Letters} 24(1-2):47--56.

\bibitem[{McDaniel \protect\BIBand{} Devine(1977)}]{mcdaniel1977modified}
McDaniel D, Devine M (1977) A modified {Benders'} partitioning algorithm for
  mixed integer programming. \emph{Management Science} 24(3):312--319.

\bibitem[{Morton(1998)}]{morton1998stopping}
Morton DP (1998) Stopping rules for a class of sampling-based stochastic
  programming algorithms. \emph{Operations Research} 46(5):710--718.

\bibitem[{Papadakos(2008)}]{papadakos2008practical}
Papadakos N (2008) Practical enhancements to the {M}agnanti--{W}ong method.
  \emph{Operations Research Letters} 36(4):444--449.

\bibitem[{Pereira \protect\BIBand{} Pinto(1991)}]{pereira1991multi}
Pereira MV, Pinto LM (1991) Multi-stage stochastic optimization applied to
  energy planning. \emph{Mathematical Programming} 52(1):359--375.

\bibitem[{Pishvaee et~al.(2014)Pishvaee, Razmi, \protect\BIBand{}
  Torabi}]{pishvaee_accelerated_2014}
Pishvaee MS, Razmi J, Torabi SA (2014) An accelerated {B}enders decomposition
  algorithm for sustainable supply chain network design under uncertainty: A
  case study of medical needle and syringe supply chain. \emph{Transportation
  Research Part E: Logistics and Transportation Review} 67:14--38.

\bibitem[{Rahmaniani et~al.(2018)Rahmaniani, Crainic, Gendreau,
  \protect\BIBand{} Rei}]{rahmaniani2018accelerating}
Rahmaniani R, Crainic TG, Gendreau M, Rei W (2018) Accelerating the benders
  decomposition method: Application to stochastic network design problems.
  \emph{SIAM Journal on Optimization} 28(1):875--903.

\bibitem[{Ram{\'\i}rez-Pico et~al.(2023)Ram{\'\i}rez-Pico, Ljubi{\'c},
  \protect\BIBand{} Moreno}]{ramirez2023benders}
Ram{\'\i}rez-Pico C, Ljubi{\'c} I, Moreno E (2023) Benders adaptive-cuts method
  for two-stage stochastic programs. \emph{Transportation Science}
  57(5):1252--1275.

\bibitem[{Rei et~al.(2009)Rei, Cordeau, Gendreau, \protect\BIBand{}
  Soriano}]{rei2009accelerating}
Rei W, Cordeau JF, Gendreau M, Soriano P (2009) Accelerating {B}enders
  decomposition by local branching. \emph{INFORMS Journal on Computing}
  21(2):333--345.

\bibitem[{Reuther et~al.(2018)Reuther, Kepner, Byun, Samsi, Arcand, Bestor,
  Bergeron, Gadepally, Houle, Hubbell, Jones, Klein, Milechin, Mullen, Prout,
  Rosa, Yee, \protect\BIBand{} Michaleas}]{reuther2018interactive}
Reuther A, Kepner J, Byun C, Samsi S, Arcand W, Bestor D, Bergeron B, Gadepally
  V, Houle M, Hubbell M, Jones M, Klein A, Milechin L, Mullen J, Prout A, Rosa
  A, Yee C, Michaleas P (2018) Interactive supercomputing on 40,000 cores for
  machine learning and data analysis. \emph{2018 IEEE High Performance extreme
  Computing Conference (HPEC)}, 1--6 (IEEE).

\bibitem[{Richardson(1976)}]{richardson1976optimization}
Richardson R (1976) An optimization approach to routing aircraft.
  \emph{Transportation Science} 10(1):52--71.

\bibitem[{Rodr{\'\i}guez-Mart{\'\i}n \protect\BIBand{}
  Salazar-Gonz{\'a}lez(2010)}]{rodriguez2010local}
Rodr{\'\i}guez-Mart{\'\i}n I, Salazar-Gonz{\'a}lez JJ (2010) A local branching
  heuristic for the capacitated fixed-charge network design problem.
  \emph{Computers \& Operations Research} 37(3):575--581.

\bibitem[{Santoso et~al.(2005)Santoso, Ahmed, Goetschalckx, \protect\BIBand{}
  Shapiro}]{santoso_stochastic_2005}
Santoso T, Ahmed S, Goetschalckx M, Shapiro A (2005) A stochastic programming
  approach for supply chain network design under uncertainty. \emph{European
  Journal of Operational Research} 167(1):96--115.

\bibitem[{Schmidt et~al.(2017)Schmidt, Le~Roux, \protect\BIBand{}
  Bach}]{schmidt2017minimizing}
Schmidt M, Le~Roux N, Bach F (2017) Minimizing finite sums with the stochastic
  average gradient. \emph{Mathematical Programming} 162(1):83--112.

\bibitem[{Shapiro et~al.(2021)Shapiro, Dentcheva, \protect\BIBand{}
  Ruszczynski}]{shapiro2021lectures}
Shapiro A, Dentcheva D, Ruszczynski A (2021) \emph{Lectures on Stochastic
  Programming: Modeling and Theory} (SIAM).

\bibitem[{Smith \protect\BIBand{} Winkler(2006)}]{smith2006optimizer}
Smith JE, Winkler RL (2006) The optimizer’s curse: Skepticism and
  postdecision surprise in decision analysis. \emph{Management Science}
  52(3):311--322.

\bibitem[{Stubbs \protect\BIBand{} Mehrotra(1999)}]{stubbs1999branch}
Stubbs RA, Mehrotra S (1999) A branch-and-cut method for 0-1 mixed convex
  programming. \emph{Mathematical Programming} 86(3):515--532.

\bibitem[{Trukhanov et~al.(2010)Trukhanov, Ntaimo, \protect\BIBand{}
  Schaefer}]{trukhanov2010adaptive}
Trukhanov S, Ntaimo L, Schaefer A (2010) Adaptive multicut aggregation for
  two-stage stochastic linear programs with recourse. \emph{European Journal of
  Operational Research} 206(2):395--406.

\bibitem[{Van~Roy \protect\BIBand{} Wolsey(1985)}]{van1985valid}
Van~Roy TJ, Wolsey LA (1985) Valid inequalities and separation for
  uncapacitated fixed charge networks. \emph{Operations Research Letters}
  4(3):105--112.

\bibitem[{Van~Slyke \protect\BIBand{} Wets(1969)}]{van1969shaped}
Van~Slyke RM, Wets R (1969) L-shaped linear programs with applications to
  optimal control and stochastic programming. \emph{SIAM Journal on Applied
  Mathematics} 17(4):638--663.

\bibitem[{Wei et~al.(2022)Wei, G{\'o}mez, \protect\BIBand{}
  K{\"u}{\c{c}}{\"u}kyavuz}]{wei2022ideal}
Wei L, G{\'o}mez A, K{\"u}{\c{c}}{\"u}kyavuz S (2022) Ideal formulations for
  constrained convex optimization problems with indicator variables.
  \emph{Mathematical Programming} 192(1):57--88.

\bibitem[{Wets(1966)}]{wets1966programming}
Wets RJB (1966) Programming under uncertainty: the equivalent convex program.
  \emph{SIAM Journal on Applied Mathematics} 14(1):89--105.

\bibitem[{Xie \protect\BIBand{} Deng(2020)}]{xie2020scalable}
Xie W, Deng X (2020) Scalable algorithms for the sparse ridge regression.
  \emph{SIAM Journal on Optimization} 30(4):3359--3386.

\bibitem[{You \protect\BIBand{} Grossmann(2013)}]{you2013multicut}
You F, Grossmann IE (2013) Multicut {B}enders decomposition algorithm for
  process supply chain planning under uncertainty. \emph{Annals of Operations
  Research} 210(1):191--211.

\bibitem[{Zakeri et~al.(2000)Zakeri, Philpott, \protect\BIBand{}
  Ryan}]{zakeri_inexact_2000}
Zakeri G, Philpott AB, Ryan DM (2000) Inexact cuts in {B}enders decomposition.
  \emph{SIAM Journal on Optimization} 10(3):643--657.

\end{thebibliography}
}
\newpage
\counterwithin{table}{section}
\counterwithin{figure}{section}
\counterwithin{theorem}{section}
\counterwithin{proposition}{section}

\setcounter{table}{0}
\setcounter{figure}{0}
\begin{APPENDICES}

\section{Justification for the Strongly Quadratic Penalty Term in Problem \eqref{eqn:prob_description}} \label{sec:a.regularization.theory}
In this section, we justify using the quadratic regularization term in Problem \eqref{eqn:prob_description}, from both a practical perspective and a theoretical one. We remark that the use of a regularization term in mixed-integer optimization is an increasingly popular modeling choice which has been discussed in detail in other works; we refer to \citet{bertsimas_unified_2021, bertsimas_scalable_2022} for a more detailed discussion of this matter.

From a practical perspective, a strongly quadratic term in the objective can model quadratic transportation costs or can be used to increase the robustness of the solution to parameter uncertainty. Indeed, as we show in Proposition \ref{prop:robust.equiv}, the second-stage cost in Problem \eqref{eqn:prob_description} is equivalent to a worst-case cost with an ellipsoidal uncertainty set around the second-stage transportation costs $f_{i,j}^k$. Furthermore, as advocated in \citet{bertsimas_unified_2021}, artificially introducing a quadratic penalty with $\gamma > 0$ is an efficient smoothing technique 
for MIO problems with logical constraints, approximating the nominal objective function (to arbitrary precision by taking $\gamma \rightarrow + \infty$) while improving computational tractability.
Formally, if $v(\gamma)$ denotes the optimal objective value of \eqref{eqn:prob_description}, one can easily show that 
$$  v(\gamma) - \dfrac{1}{2 \gamma}  \sum_{(i,j) \in \mathcal{E}} u_{i,j}^2  \: \leq \:  v(\infty) \: \leq \: v(\gamma). $$

Further, we should emphasize that the method we develop in this paper does not require $1 / \gamma > 0$ and applies in the case where there is no quadratic term in the objective as well ($\gamma \rightarrow \infty$ and $1 / \gamma = 0$), as discussed in \citet[Remark 2.6]{bertsimas_unified_2021}.  

We now provide some theoretical evidence to justify the presence of a smooth strongly convex term in the objective of \eqref{eqn:prob_description}. In particular, we show that adding this term can be interpreted as equivalent to considering a robust version of the linear objective.
\begin{proposition} \label{prop:robust.equiv}
Fix $\gamma > 0$ and $\bs{z} \in \{0,1\}^{\mathcal{E}}$. Let us denote $\mathcal{X}$ the set of feasible second-stage flow transportation variables $\bs{x}$, i.e., 
\begin{align*}
    \mathcal{X} := \left\lbrace \bs{x} \geq \bs{0} \: \left| 
    \begin{array}{ll}
		\bs{A}\bs{x}^{k,r} = \bs{d}^{k,r}, &\quad \forall k \in \mathcal{K}, r \in \mathcal{R} \\
		\sum_{k \in \mathcal{K}} x^{k,r}_{i,j} \leq u_{i,j},  
		&\quad \forall (i,j) \in \mathcal{E}, r \in \mathcal{R} \\
		x^{k,r}_{i,j} = 0 \textit{ if } z_{i,j} = 0, &\quad \forall (i,j) \in \mathcal{E}
    \end{array}
    \right. \right\rbrace.
\end{align*}

There exists a parameter value $\lambda \geq 0$ such that 
\begin{equation} \label{eqn:justify.gamma.penalty}
	\begin{aligned}
		\min_{\bs{x} \in \mathcal{X}} \quad & 
		\sum_{r \in \mathcal{R}} \sum_{(i,j) \in \mathcal{E}} \sum_{k \in \mathcal{K}} f_{i,j}^{k}x_{i,j}^{k,r} + \frac{1}{2 \gamma} \sum_{r \in \mathcal{R}} \sum_{(i,j) \in \mathcal{E}} \left( \sum_{k \in \mathcal{K}} x_{i,j}^{k,r}\right)^2
	\end{aligned}
\end{equation}
achieves the same optimal solution as
\begin{equation}\label{eqn:justify.gamma.robust}
	\begin{aligned}
		\min_{\bs{x} \in \mathcal{X}} \quad & \sum_{r \in \mathcal{R}} \sum_{(i,j) \in \mathcal{E}} \max \left\lbrace \sum_{k \in \mathcal{K}} \Tilde{f}_{i,j}^{k}x_{i,j}^{k,r} \: : \: \Tilde{f}_{i,j}^{k} = {f}_{i,j}^{k} + \zeta^{r}_{i,j},\ \sum_{r \in \mathcal{R}} \sum_{(i,j) \in \mathcal{E}} (\zeta^{r}_{i,j})^2 \leq \lambda^2 \right\rbrace .
        \end{aligned}
\end{equation} 
\end{proposition}
\proof{Proof of Proposition \ref{prop:robust.equiv}} 
Problems \eqref{eqn:justify.gamma.penalty} and \eqref{eqn:justify.gamma.robust} have the same feasible set. Consider a feasible solution $\bs{x}$. We have
\begin{align*}
    \sum_{r \in \mathcal{R}} \sum_{(i,j) \in \mathcal{E}} \sum_{k \in \mathcal{K}} \Tilde{f}_{i,j}^{k}x_{i,j}^{k,r} 
    &=   \sum_{r \in \mathcal{R}} \sum_{(i,j) \in \mathcal{E}} \sum_{k \in \mathcal{K}} {f}_{i,j}^{k}x_{i,j}^{k,r} +  \sum_{r \in \mathcal{R}} \sum_{(i,j) \in \mathcal{E}} \sum_{k \in \mathcal{K}} \Tilde{\zeta}_{i,j}^{r} x_{i,j}^{k,r} \\
    &=   \sum_{r \in \mathcal{R}} \sum_{(i,j) \in \mathcal{E}} \sum_{k \in \mathcal{K}} {f}_{i,j}^{k}x_{i,j}^{k,r} +  \sum_{r \in \mathcal{R}} \sum_{(i,j) \in \mathcal{E}} \Tilde{\zeta}_{i,j}^{r} \left( \sum_{k \in \mathcal{K}} x_{i,j}^{k,r} \right). 
\end{align*}
Hence, the worst-case value with respect to all vectors $\Tilde{\bs{\zeta}}$ such that $\displaystyle \sum_{r \in \mathcal{R}} \sum_{(i,j) \in \mathcal{E}} (\Tilde{\zeta}_{i,j}^{r})^2 \leq \lambda^2$ (i.e., the value of the inner maximization problem in \eqref{eqn:justify.gamma.robust}) is equal to
\begin{align*}
    \sum_{r \in \mathcal{R}} \sum_{(i,j) \in \mathcal{E}} \sum_{k \in \mathcal{K}} {f}_{i,j}^{k}x_{i,j}^{k,r} +  \lambda \sqrt{ \sum_{r \in \mathcal{R}} \sum_{(i,j) \in \mathcal{E}} \left( \sum_{k \in \mathcal{K}} x_{i,j}^{k,r} \right)^2 } =: \sum_{r \in \mathcal{R}} \sum_{(i,j) \in \mathcal{E}} \sum_{k \in \mathcal{K}} {f}_{i,j}^{k}x_{i,j}^{k,r} +  \lambda \sqrt{ q(\bs{x}) } . 
\end{align*}
Hence, minimizing the worst-case transportation cost in \eqref{eqn:justify.gamma.robust} is equivalent to minimizing the nominal cost, $\sum_{r} \sum_{(i,j)} \sum_{k} {f}_{i,j}^{k}x_{i,j}^{k,r}$  plus a penalty term, which is equal to the square root of the quadratic regularization term in \eqref{eqn:justify.gamma.penalty}, $q(\bs{x})$. To conclude the proof, we need to show that for any $\lambda > 0$ there exists a parameter value $\gamma > 0$ such that the penalties $\lambda \sqrt{q(\bs{x})}$ and $\tfrac{1}{2\gamma} q(\bs{x})$ lead to the same optimal solution.

For any $\lambda > 0$, by duality, there exists a constant $q_0$ such that Problem \eqref{eqn:justify.gamma.robust} is equivalent to 
\begin{equation} \label{eqn:constraint.reg} 
\begin{aligned}
		\min_{\bs{x} \in \mathcal{X}} \quad & \sum_{r \in \mathcal{R}} \sum_{(i,j) \in \mathcal{E}} \sum_{k \in \mathcal{K}} {f}_{i,j}^{k}x_{i,j}^{k,r} \quad  \mbox{s.t.} \quad \sqrt{q(\bs{x})} \leq \sqrt{q_0}.
\end{aligned}
\end{equation}
Since $t \mapsto \sqrt{t}$ is increasing over $\mathbb{R}_+$, constraint \eqref{eqn:constraint.reg} is equivalent to $q(\bs{x}) \leq q_0$ and the resulting problem is in turn equivalent to a problem of the same form as \eqref{eqn:justify.gamma.penalty}. \hfill \Halmos
\endproof
Proposition \ref{prop:robust.equiv} shows that the regularized objective in \eqref{eqn:prob_description} is equivalent to a robust linear objective with ellipsoidal uncertainty set, \emph{for each first-stage design decision $\bs{z}$}. Optimizing (over $\bs{z}$) for each objective can result in different solutions $(\bs{z}, \bs{x})$ though, because the first stage decision $\bs{z}$ is discrete so the value of $\lambda$ that makes the penalized formulation equivalent to the constrained one is $\bs{z}$-specific (the equivalence requires a duality type of argument, which does not hold in general when jointly optimizing for $(\bs{z},\bs{x})$ with $\bs{z}$ binary). However, after reformulating the logical constraints via algebraic linear/second-order cone constraints, we can show an equivalence results between the Boolean relaxations (proof omitted):
\begin{proposition} \label{prop:robust.equiv2}
Fix $\gamma > 0$. Let us denote $\mathcal{P}$ the set of feasible relaxed decision variables $(\bs{z}, \bs{x})$, i.e., 
\begin{align*}
    \mathcal{P} := \left\lbrace (\bs{z}, \bs{x}) \: \left| 
    \begin{array}{ll}
        \bs{z} \in [0,1]^{\mathcal{E}}, \\ 
        \sum_{(i,j) \in \mathcal{E}} z_{i,j} \leq c_0, \\
		\bs{A}\bs{x}^{k,r} = \bs{d}^{k,r}, &\quad \forall k \in \mathcal{K}, r \in \mathcal{R} \\
		\sum_{k \in \mathcal{K}} x^{k,r}_{i,j} \leq u_{i,j} z_{i,j},  
		&\quad \forall (i,j) \in \mathcal{E}, r \in \mathcal{R} \\
        \bs{x} \geq \bs{0}
    \end{array}
    \right. \right\rbrace.
\end{align*}

There exists a parameter value $\lambda \geq 0$ such that 
\begin{equation*}
	\begin{aligned}
		\min_{(\bs{z},\bs{x}) \in \mathcal{P}} \quad & \sum_{(i,j) \in \mathcal{E}} c_{i,j}z_{i,j} + 
		\dfrac{1}{|\mathcal{R}|} \sum_{r \in \mathcal{R}} \sum_{(i,j) \in \mathcal{E}} \left( \sum_{k \in \mathcal{K}} f_{i,j}^{k}x_{i,j}^{k,r} + \frac{1}{2 \gamma} \left( \sum_{k \in \mathcal{K}} x_{i,j}^{k,r}\right)^2 \right)
	\end{aligned}
\end{equation*}
achieves the same optimal solution as
\begin{equation*}
	\begin{aligned}
		\min_{(\bs{z},\bs{x}) \in \mathcal{P}} \quad & \sum_{(i,j) \in \mathcal{E}} c_{i,j}z_{i,j} + 
		\dfrac{1}{|\mathcal{R}|} \sum_{r \in \mathcal{R}} \sum_{(i,j) \in \mathcal{E}} \max \left\lbrace \sum_{k \in \mathcal{K}} \Tilde{f}_{i,j}^{k}x_{i,j}^{k,r} \: : \: \Tilde{f}_{i,j}^{k} = {f}_{i,j}^{k} + \zeta^{r}_{i,j},\ \sum_{r \in \mathcal{R}} \sum_{(i,j) \in \mathcal{E}} (\zeta^{r}_{i,j})^2 \leq \lambda^2 \right\rbrace .
        \end{aligned}
\end{equation*} 
\end{proposition}

\newpage
\section{Omitted Proofs}
\subsection{Proof of Proposition \ref{prop:f.dual}} \label{ssec:a.proof.fdual}
\proof{Proof of Proposition \ref{prop:f.dual}} The minimization problem defining  $f(\bs{z}; \bs{d})$ can be seen as the sum of two minimization problems
\begin{align*}
\min_{\boldsymbol{x}^k \in \mathbb{R}_+^{\mathcal{E}}, k \in \mathcal{K}} \:
 \sum_{k \in \mathcal{K}} \langle \bs{f}^{k}, \bs{x}^{k}\rangle  \mbox{ s.t. } 
\bs{Ax}^{k} = \bs{d}^{k}, \, \forall k \in \mathcal{K},
\end{align*}
and
\begin{align*}
\min_{\boldsymbol{y} \in \mathbb{R}^{\mathcal{E}}} \:  \frac{1}{2 \gamma} \sum_{{(i,j)} \in\mathcal{E} } {y}_{i,j}^2 \mbox{ s.t. } & y_{i,j} \leq u_{i,j}, \, \forall (i,j) \in \mathcal{E}, \\
& y_{i,j} = 0 \textit{ if } z_{i,j} = 0, \, \forall (i,j) \in \mathcal{E},
\end{align*}
coupled via the constraints $\sum_{k \in \mathcal{K}} x_{i,j}^k = y_{i,j}, \forall (i,j) \in \mathcal{E}$. Therefore, by associating a dual variable $\alpha_{i,j} \in \mathbb{R}$ with each coupling constraint, we rewrite $f(\bs{z}; \bs{d})$ as
\begin{align*}
\min_{\substack{\boldsymbol{x}^k \in \mathbb{R}_+^{\mathcal{E}}, k \in \mathcal{K}: \\ \bs{Ax}^{k} = \bs{d}^{k}, \, \forall k \in \mathcal{K}}} \quad
\min_{\substack{\boldsymbol{y} \in \mathbb{R}^{\mathcal{E}}: \\ y_{i,j} \leq u_{i,j}, \, \forall (i,j) \in \mathcal{E} \\ y_{i,j} = 0 \textit{ if } z_{i,j} = 0, \, \forall (i,j) \in \mathcal{E} }} \: 
\max_{\bs{\alpha} \in \mathbb{R}_+^{\mathcal{E}}} \: 
\sum_{k \in \mathcal{K}} \langle \bs{f}^{k} - \bs{\alpha}, \bs{x}^{k}\rangle  + \sum_{{(i,j)} \in\mathcal{E} } \left( \alpha_{i,j} {y}_{i,j} +  \tfrac{1}{2 \gamma}{y}_{i,j}^2 \right).
\end{align*}

By invoking standard results on saddle-point theorems \citep[see, e.g.,][]{bertsekas1999nonlinear}, the order of the minimization and maximization operators on the function $f(\bs{z}, \bs{d})$ can be exchanged\footnote{In general, we require that a constraint qualification holds to be able to exchange the order of minimization and maximization operators \citep[see, e.g.,][]{bertsekas1999nonlinear}. However, all constraints in Problem \eqref{eqn:f_z_definition_y} are linear and it has a convex quadratic objective. Therefore, we can exchange the order of the operators in an assumption-free manner.} without altering the objective value. Moreover, after exchanging these operators, we can compute the dual of each minimization problem separately. Indeed,
\begin{align*}
\min_{\substack{\boldsymbol{x}^k \in \mathbb{R}_+^{\mathcal{E}}, k \in \mathcal{K}: \\ \bs{Ax}^{k} = \bs{d}^{k}, \, \forall k \in \mathcal{K}}} \:
 \sum_{k \in \mathcal{K}} \langle \bs{f}^{k} - \bs{\alpha}, \bs{x}^{k}\rangle 
 \quad = \quad  
 \max_{\substack{\bs{p}^k \in \mathbb{R}^{\mathcal{N}}: \\ \bs{A}^\top \bs{p}^{k} \leq \bs{f}^{k} - \bs{\alpha}, \forall k \in \mathcal{K}}} \:
		\sum_{k \in \mathcal{K}} \langle \bs{p}^{k}, \bs{d}^{k} \rangle.
\end{align*}
Second, to dualize
\begin{align*}
\min_{\boldsymbol{y} \in \mathbb{R}^{\mathcal{E}}} \: \sum_{{(i,j)} \in\mathcal{E} } \left( \alpha_{i,j} y_{i,j} + \tfrac{1}{2 \gamma}{y}_{i,j}^2 \right)\mbox{ s.t. } &  y_{i,j} \leq u_{i,j}, \, \forall (i,j) \in \mathcal{E}, \\
& y_{i,j} = 0 \textit{ if } z_{i,j} = 0, \, \forall (i,j) \in \mathcal{E},
\end{align*}
let us first observe that we can omit the logical constraints by considering the change of variables $y_{i,j} = z_{i,j} w_{i,j}$ for $\boldsymbol{w} \in \mathbb{R}^{\mathcal{E}}$. Hence, we obtain
\begin{align*}
& \min_{\boldsymbol{w} \in \mathbb{R}^{\mathcal{E}}} \: \sum_{{(i,j)} \in\mathcal{E} } \left[ z_{i,j}\alpha_{i,j} w_{i,j} + \tfrac{1}{2 \gamma} z_{i,j}{w}_{i,j}^2 \right] \mbox{ s.t. } z_{i,j} w_{i,j} \leq z_{i,j} u_{i,j}, \, \forall (i,j) \in \mathcal{E} \\
= & \max_{\bs{\beta} \in \mathbb{R}_+^{\mathcal{E}}} \: 
\min_{\boldsymbol{y} \in \mathbb{R}^{\mathcal{E}}} \: 
 - \sum_{(i,j) \in \mathcal{E}} z_{i,j} \beta_{i,j} u_{i,j}
\, + \sum_{{(i,j)} \in\mathcal{E} } \left[ z_{i,j}(\alpha_{i,j} + \beta_{i,j}) w_{i,j} + \tfrac{1}{2 \gamma} z_{i,j}{w}_{i,j}^2 \right] \\
= & \max_{\bs{\beta} \in \mathbb{R}_+^{\mathcal{E}}} \: 
 - \sum_{(i,j) \in \mathcal{E}} z_{i,j} \beta_{i,j} u_{i,j}
\, - \dfrac{\gamma}{2} \sum_{{(i,j)} \in\mathcal{E} } z_{i,j} (\alpha_{i,j} + \beta_{i,j})^2.
\end{align*}

All together, we obtain the desired reformulation. \hfill \Halmos
\endproof

\subsection{Proof of Proposition \ref{prop:avg.approx}} \label{ssec:a.proof.approx}
In this section, we provide a proof of Proposition \ref{prop:avg.approx}. To clarify the presentation, we adopt a lighter set of notations: 

Fix $\bs{z}$. For any $r \in \mathcal{R}$, we denote  $\bs{\xi}^r := (\bs{\alpha}^r, \bs{\beta}^r, \bs{p}^r)$ the optimal dual solutions of \eqref{eqn:f.dual} for $\bs{d} = \bs{d}^{\cdot,r}$. For any subset $\mathcal{S} \subseteq \mathcal{R}$, let us denote  $\bar{\bs{\xi}}^{\mathcal{S}} := \tfrac{1}{|\mathcal{S}|} \sum_{r \in \mathcal{S}} \bs{\xi}^r$ the average of the optimal dual solutions $\bs{\xi}^r$ for $r \in \mathcal{S}$. For a random $\mathcal{S} \subseteq \mathcal{R}$ of fixed size $|\mathcal{S}|$, 
we will analyze the sub-optimality gap of $\bar{\bs{\xi}}^{\mathcal{S}}$, i.e., the quantity $q(\bs{z}, \bs{\xi}^r; \bs{d}^r) - q(\bs{z}, \bar{\bs{\xi}}^{\mathcal{S}}; \bs{d}^r) \, (\geq 0)$, for scenarios $r \in \mathcal{S}^c := \mathcal{R} \setminus \mathcal{S}$. 

\proof{Proof of Proposition \ref{prop:avg.approx}} 
Let us denote $M := \max_{r \in \mathcal{R}} \| \bs{\xi}^r \|_\infty$. Since $\| \bs{\xi}^r \|_\infty \leq M$, then $\left\| \bar{\bs{\xi}}^{\mathcal{S}} \right\|_\infty \leq M$ and there exists some constant $L > 0$ such that, for any $\mathcal{S}$ and any $r \notin \mathcal{S}$
\begin{align*}
    \left|  q(\bs{z}, \bar{\bs{\xi}}^{\mathcal{S}}, \bs{d}^r) - q(\bs{z}, {\bs{\xi}^r}, \bs{d}^r) \right| 
&\leq L \left\| \bar{\bs{\xi}}^{\mathcal{S}} - {\bs{\xi}}^r \right\|.
\end{align*}
We further decompose the right-hand side via a triangle inequality and sum the inequalities above across all $r \notin \mathcal{S}$ to obtain
\begin{align*}
    \sum_{r \in \mathcal{S}^c} \left|  q(\bs{z}, \bar{\bs{\xi}}^{\mathcal{S}}, \bs{d}^r) - q(\bs{z}, {\bs{\xi}^r}, \bs{d}^r) \right| &\leq  L |\mathcal{S}^c| \left\| \bar{\bs{\xi}}^{\mathcal{S}} - \bar{\bs{\xi}}^{\mathcal{R}} \right\| +  L \, \sum_{r \in \mathcal{S}^c} \left\| \bar{\bs{\xi}}^{\mathcal{R}} - {\bs{\xi}}^r \right\|.
\end{align*}

The first term corresponds to the difference between $\bar{\bs{\xi}}^{\mathcal{R}}$ and an unbiased estimate obtained via sampling without replacement. Denote $d$ the dimension of $\bs{\xi}$. Hence, since the $\bs{\xi}^r$ are uniformly bounded, by \citet[][corollary 2.5]{bardenet2015concentration}, there exists some universal constant $M_1$ such that for any $\delta > 0$, we have, with probability $1-\delta$ on the subset $\mathcal{S}$ of fixed size $|\mathcal{S}|$,
\begin{align}  \label{eqn:concentration.mean}
\left\| \bar{\bs{\xi}}^{\mathcal{S}} - \bar{\bs{\xi}}^{\mathcal{R}} \right\| \leq M_1 \sqrt{ d \left(\dfrac{1}{|\mathcal{S}|} - \dfrac{1}{|\mathcal{R}|} \right)\log(1/\delta)} , 
\end{align}

For the second term, we simply use the bound
\begin{align*}
    \sum_{r \in \mathcal{S}^c} \left\| \bar{\bs{\xi}}^{\mathcal{R}} - {\bs{\xi}}^r \right\| 
    \leq \sqrt{|\mathcal{S}^c|} \sqrt{ \dfrac{1}{{|\mathcal{S}^c|}} \sum_{r \in \mathcal{S}^c} \left\| \bar{\bs{\xi}}^{\mathcal{R}} - {\bs{\xi}}^r \right\|^2}.
\end{align*}
For interpretation, we denote $\nu^2 := \dfrac{1}{|\mathcal{R}|}\sum_{r \in \mathcal{R}} \left\| \bar{\bs{\xi}}^{\mathcal{R}} - {\bs{\xi}}^r \right\|^2$, which can be interpreted as the variance in optimal dual variables of our problem. Then, the term on the right-hand side of the inequality above can be viewed as a bootstrap estimator of $\nu$, which intuitively converges to $\nu$ as $\mathcal{S}^c \rightarrow \mathcal{R}$. To formalize this intuition, let us expand the squared norm term and apply the triangle inequality:
\begin{align*}
    \left| \dfrac{1}{{|\mathcal{S}^c|}} \sum_{r \in \mathcal{S}^c} \left\| \bar{\bs{\xi}}^{\mathcal{R}} - {\bs{\xi}}^r \right\|^2 - \nu^2 \right| 
    & = \left| \dfrac{1}{{|\mathcal{S}^c|}} \sum_{r \in \mathcal{S}^c} \left\| {\bs{\xi}}^r \right\|^2 - 2 \langle \bar{\bs{\xi}}^{\mathcal{R}}, \bar{\bs{\xi}}^{\mathcal{S}^c}\rangle - \dfrac{1}{{|\mathcal{R}|}} \sum_{r \in \mathcal{R}} \left\| {\bs{\xi}}^r \right\|^2  + 2 \langle \bar{\bs{\xi}}^{\mathcal{R}}, \bar{\bs{\xi}}^{\mathcal{R}}\rangle \right| \\ 
    &\leq \left| \dfrac{1}{{|\mathcal{S}^c|}} \sum_{r \in \mathcal{S}^c} \left\| {\bs{\xi}}^r \right\|^2  - \dfrac{1}{{|\mathcal{R}|}} \sum_{r \in \mathcal{R}} \left\| {\bs{\xi}}^r \right\|^2  \right| 
    + 2 \left| \langle \bar{\bs{\xi}}^{\mathcal{R}},  \bar{\bs{\xi}}^{\mathcal{R}} - \bar{\bs{\xi}}^{\mathcal{S}^c} \rangle \right| \\ 
     &\leq \left| \dfrac{1}{{|\mathcal{S}^c|}} \sum_{r \in \mathcal{S}^c} \left\| {\bs{\xi}}^r \right\|^2  - \dfrac{1}{{|\mathcal{R}|}} \sum_{r \in \mathcal{R}} \left\| {\bs{\xi}}^r \right\|^2  \right| 
    + 2  M \left\|  \bar{\bs{\xi}}^{\mathcal{R}} - \bar{\bs{\xi}}^{\mathcal{S}^c} \right\|
\end{align*}
By \citet[][corollary 2.5]{bardenet2015concentration} again, there exists $M_2 > 0$ such that, with probability $1-\delta$, 
\begin{align*}
    \left| \dfrac{1}{{|\mathcal{S}^c|}} \sum_{r \in \mathcal{S}^c} \left\| {\bs{\xi}}^r \right\|^2  - \dfrac{1}{{|\mathcal{R}|}} \sum_{r \in \mathcal{R}} \left\| {\bs{\xi}}^r \right\|^2  \right| 
    \leq M_2 \sqrt{\left(\dfrac{1}{|\mathcal{S}^c|} - \dfrac{1}{|\mathcal{R}|} \right)\log(1/\delta)}, 
\end{align*}
and $\left\|  \bar{\bs{\xi}}^{\mathcal{R}} - \bar{\bs{\xi}}^{\mathcal{S}^c} \right\|$ satisfies a similar inequality as \eqref{eqn:concentration.mean}. All together, with probability $1-2\delta$,
\begin{align*}
    \left| \dfrac{1}{{|\mathcal{S}^c|}} \sum_{r \in \mathcal{S}^c} \left\| \bar{\bs{\xi}}^{\mathcal{R}} - {\bs{\xi}}^r \right\|^2 - \nu^2 \right| 
     &\leq M_2 \sqrt{\left(\dfrac{1}{|\mathcal{S}^c|} - \dfrac{1}{|\mathcal{R}|} \right)\log(1/\delta)}
    + M_1 \sqrt{d \left(\dfrac{1}{|\mathcal{S}^c|} - \dfrac{1}{|\mathcal{R}|} \right)\log(1/\delta)},
\end{align*}
yielding 
\begin{align} \label{eqn:concentration.var}
    \sum_{r \in \mathcal{S}^c} \left\| \bar{\bs{\xi}}^{\mathcal{R}} - {\bs{\xi}}^r \right\| 
    \leq \sqrt{|\mathcal{S}^c|} \nu  + \sqrt{|\mathcal{S}^c|} M_3 \left(\dfrac{1}{|\mathcal{S}^c|} - \dfrac{1}{|\mathcal{R}|} \right)^{1/4} \left(\log(1/\delta)\right)^{1/4},
\end{align}
with $M_3 := M_2 + M_1\sqrt{d}$.

Combining \eqref{eqn:concentration.mean} and \eqref{eqn:concentration.var}, we obtain that, with probability $1-3\delta$ over the sample $\mathcal{S}$,
\begin{align*}
    \sum_{r \in \mathcal{S}^c} \left|  q(\bs{z}, \bar{\bs{\xi}}^{\mathcal{S}}, \bs{d}^r) - q(\bs{z}, {\bs{\xi}^r}, \bs{d}^r) \right| 
    &\leq  L \sqrt{|\mathcal{S}^c|} \nu  + E
\end{align*}
where $E$ is a bootstrap error term equal to
\begin{align*}
    E &=  L |\mathcal{S}^c| M_1 \sqrt{ d \left(\dfrac{1}{|\mathcal{S}|} - \dfrac{1}{|\mathcal{R}|} \right)\log(1/\delta)} + \sqrt{|\mathcal{S}^c|} M_3 L \left(\dfrac{1}{|\mathcal{S}^c|} - \dfrac{1}{|\mathcal{R}|} \right)^{1/4} \left(\log(1/\delta)\right)^{1/4} \\ 
    &\leq \sqrt{|\mathcal{S}^c|} M_3 L \left[ \sqrt{|\mathcal{R}|} \left(\dfrac{1}{|\mathcal{S}|} - \dfrac{1}{|\mathcal{R}|} \right)^{1/2} + \left(\dfrac{1}{|\mathcal{S}^c|} - \dfrac{1}{|\mathcal{R}|} \right)^{1/4} \right] \sqrt{\log(1/\delta)},
\end{align*}
because $M_3 \geq M_1 \sqrt{d}$ and for $\delta$ such that $\log(1/\delta) > 1$.

To conclude the proof, let us observe that  $M_3 = M_2 + M_1\sqrt{d}$ and $d = 2|\mathcal{E}| + |\mathcal{N}| \times |\mathcal{K}|$. 
 \hfill \Halmos
\endproof

\newpage  
\section{A $k$-cut Implementation of Benders Decomposition} \label{sec:a.kcut}
To successfully combine the best aspects of single and multi-cut approaches, a $k$-cut approach was proposed by \cite{trukhanov2010adaptive, contreras2011benders}. They observed that scenarios can often be partitioned into subsets (or clusters) that are very similar to one another. Moreover, aggregating the cuts in each partition successfully compresses information about the second-stage cost surface and, on a per-iteration basis, is almost as fast as a single-cut approach. Accordingly, let $\cup_{c \in [k]} \mathcal{S}_c$ be a partition of $\mathcal{R}$. Then, at each iteration, the $k$-cut approach solves the MIO:
\begin{equation} \label{eqn:cp.kcut}
	\min _{\substack{\boldsymbol{z} \in \mathcal{Z}\\ \eta_c \in \mathbb{R}, c\in [k]}} \langle \bs{c}, \boldsymbol{z} \rangle + \dfrac{1}{|\mathcal{R}|}\sum_{c \in [k]} \eta_c \mbox{ s.t. } \eta_c \geq \sum_{r \in \mathcal{S}_c}f(\bs{z}^t; \bs{d}^{\cdot,r}) + \langle \nabla f(\bs{z}^t; \bs{d}^{\cdot,r}), \bs{z} - \bs{z}^t\rangle, \, \forall t \in [T], \, \forall c \in [k],
\end{equation}
and constructs each cut similarly to the single and multi-cut approaches. At each iteration, the $k$-cut approach adds $k$ linear constraints (one per cluster $c \in [k]$). If $k=1$ (resp. $|\mathcal{R}|$), we recover the single-cut (resp. multi-cut) algorithm. 

\subsection{Stochastic Variant of $k$-cut Benders' Decomposition}
A stochastic variant of the $k$-cut approach can be developed analogously to the stochastic single-cut approach in Section \ref{ssec:cp.stochastic}, by applying our method for stochastic single-cut to each cluster $c \in [k]$. 

Namely, we partition the set of scenarios $\mathcal{R}$ into $k$ sets $\mathcal{S}_c: c \in [k]$, and impose valid constraints of the form \eqref{eqn:scp.single.dual} to each epigraph variable $\eta_c$: 
 \begin{align} \label{eqn:scp.kcut.dual}
\eta_c \geq \sum_{r \in \mathcal{S}_c} q(\bs{z}^t, \bs{\alpha}^r, \bs{\beta}^r, \bs{p}^r; \bs{d}^{\cdot,r}) + \sum_{r \in \mathcal{S}_c} \left\langle \nabla_{\bs{z}} q(\bs{z}^t, \bs{\alpha}^r, \bs{\beta}^r, \bs{p}^r; \bs{d}^{\cdot,r}), \bs{z} - \bs{z}^t \right\rangle, \ \forall c \in [k].
\end{align}
Then, at each iteration, we sample and solve \eqref{eqn:f.dual} for a random subset $\mathcal{R}_{t,c} \subseteq \mathcal{S}_c$ of scenarios in each cluster and set $(\bs{\alpha}^r, \bs{\beta}^r,\bs{p}^r) = (\bar{\bs{\alpha}}^{\mathcal{R}_{r,c}}, \bar{\bs{\beta}}^{\mathcal{R}_{r,c}}, \bar{\bs{p}}^{\mathcal{R}_{r,c}})$ for $r \in \mathcal{S}_c \setminus \mathcal{R}_{t,c}$. From Proposition \ref{prop:avg.approx} applied to each cluster separately, we obtain that the approximation error for cluster $c$ is bounded, with high probability, by a term that depends on the variance in dual optimal variables within cluster $c$, $\nu_c^2$, plus a bootstrap estimation error term. Hence, if the clustering successfully reduces total weighted variance $\sum_{c \in [k]} \sqrt{|\mathcal{S}_c|} \nu_c$, a $k$-cut approach could improve the lower bound obtained by single-cut, while using the same number of samples per iteration.

In practice (and in our implementation), it is not feasible to cluster the set of scenarios $\mathcal{R}$ based on their associated optimal dual solutions $\bs{\alpha}^r$ \emph{at the incumbent solution}, because the clustering cannot change throughout the algorithm (it has to be independent from the incumbent). Intuitively, however, the optimization problem defining $\bs{\alpha}^r$, \eqref{eqn:f.dual}, is smooth and parametrized by $\bs{d}^{\cdot,r}$. From sensitivity analysis, the clusters obtained by applying the $k$-means algorithm on the demand vectors $\bs{d}^{\cdot,r}$ or on the optimal dual variables $\bs{\alpha}^r$ at some initial vector $\bs{z}^0$ should lead to relatively homogeneous clusters in terms of optimal dual variables $\bs{\alpha}^r$ throughout the algorithm. In our implementation, we use the latter clustering (on the optimal dual variables $\bs{\alpha}^r$ computed at the root node for the initial solution $\bs{z}^0$).

\subsection{Numerical Performance}
In this section, we report (Table \ref{tab:overall_nodes.kcut}) the performance of the stochatstic $k$-cut approach on the same instances as the one used in Section \ref{ssec:exp.scalability} (see Table \ref{tab:overall_nodes}). 
As for the other cutting-plane methods, we warmstart the algorithm with cuts obtained from solving the perspective relaxation with a single-cut stochastic cutting-plane algorithm and applying these cuts at the root node. 
Overall, we observe that the $k$-cut implementation improves upon the single-cut implementation in terms of computational time and optimality gap, although the benefit is less acute as the size of the instance grows. In terms of optimality, it achieves optimality gaps that are similar to those of the accelerated multi-cut stochastic cutting plane for the small instances but that deteriorate more as the number of nodes in the network increases. 
\begin{table}
	\centering \footnotesize 
	\caption{Runtime (in seconds) and final optimality gap (in \%) for stochastic $k$-cut approach, averaged over instances with the same number of nodes $|\mathcal{N}|$.}
	\label{tab:overall_nodes.kcut}
    \begin{tabular}{lrc}
    \toprule
     & \multicolumn{2}{c}{ Stochastic $k$-Cut} \\ 
     \cmidrule(lr){2-3}
    $|\mathcal{N}|$ &    Runtime &  Gap \\ 
    \midrule 
     10 & 86.31 & 0.39 \\
     30 & 4030.85 &	4.17 \\
     50 & 4462.14 &	2.58 \\
     70 & 4851.60 &	7.59 \\
    \midrule
    100 & 4946.72 &	6.58 \\
    150 & 5103.10 &	10.96 \\
    200 & 6028.52 &	10.54 \\
    \midrule
    300 & 4734.46  &  18.88 \\
    500 & 5456.29  &  25.83 \\
    700 & 5551.47  &  37.24 \\
    \bottomrule
    \end{tabular}
\end{table}

\newpage
\section{Computational Effect of Regularizer $\bs{\gamma}$} \label{sec:a.gamma}
The regularizing constant $\gamma$ plays a crucial role in the performance of decomposition algorithms like Benders' decomposition; see, e.g., \citet{bertsimas_unified_2021}. An appropriate value of $\gamma$ is essential for achieving optimal convergence and solution quality. When set too high, the regularizing term has minimal impact on the objective function, making the problem more challenging to solve. Conversely, when set too low, the regularizing term dominates the objective function, resulting in easier but less accurate solutions.

To illustrate this, Figure \ref{fig:gamma_experiments} presents the results of our experiments on \textbf{R} instances with 160 scenarios, where we vary the value of $\gamma$. It displays the average runtime and the objective value for each value of $\gamma$ to provide insights on how the choice of $\gamma$ affects the solution quality and computational performance of the algorithm. In the experimental results presented in Section \ref{section:experiments}, we selected an appropriate value for the regularizing constant $\gamma$ in order to strike a balance between the two extremes of the spectrum, as depicted in Figure \ref{fig:gamma_experiments}.

\begin{figure}[ht]
	\centering
	\resizebox{0.9\columnwidth}{!}{%
		\includegraphics[scale=0.2]{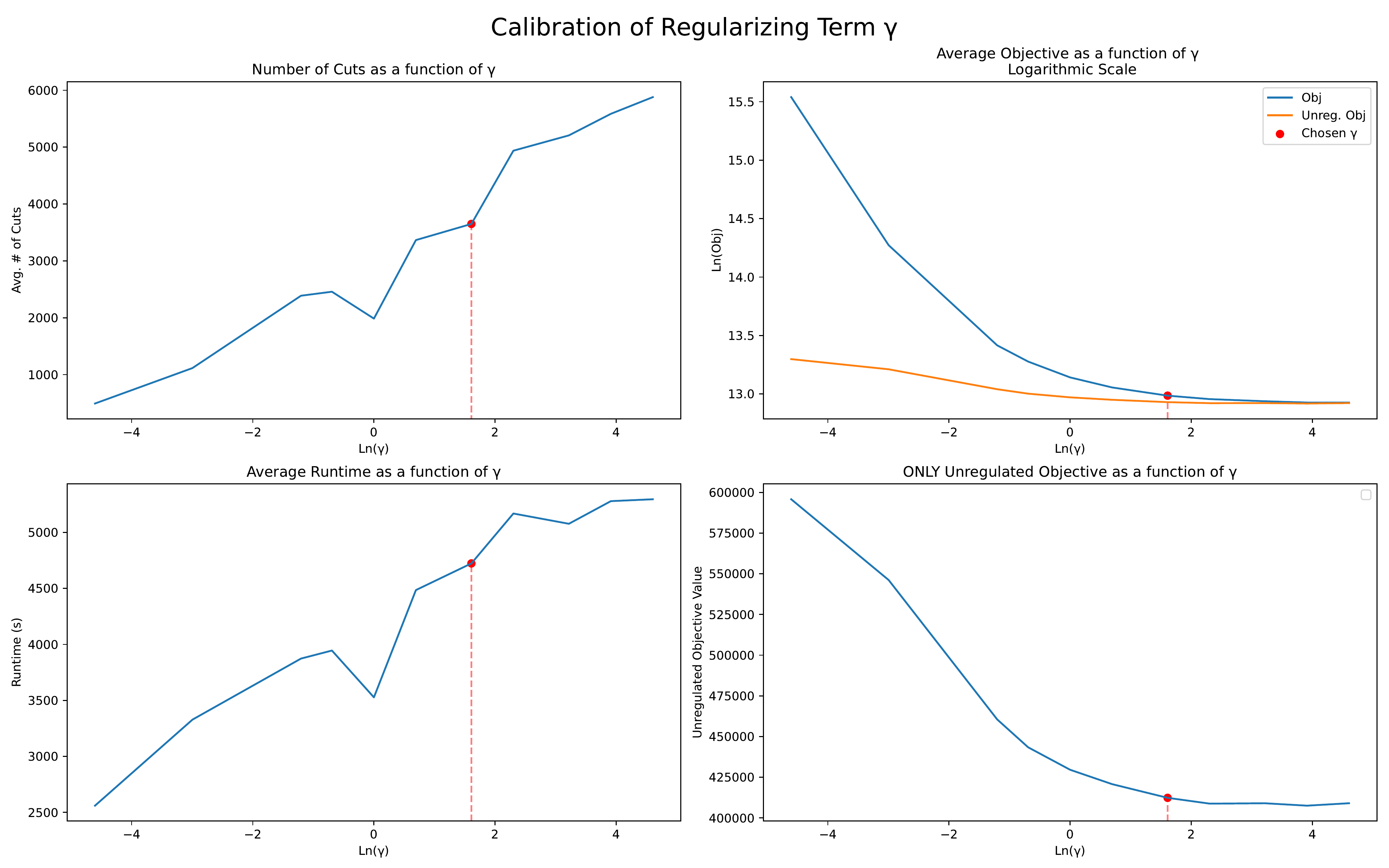}
	}
	\captionsetup{justification=centering,margin=5cm}
	\caption{Effect of Regularizer on the Algorithm's Performance.}
	\scriptsize
	\label{fig:gamma_experiments}
\end{figure}

\FloatBarrier \newpage
\section{Additional Numerical Results} \label{appendix:results}
In this section, we present additional numerical results that complement the results in Section \ref{section:experiments}.

\subsection{Instance Generation} \label{ssec:a.instance}
\subsubsection{Synthetic instances}
We generate instances according to a methodology inspired by that of \citet{gunluk_perspective_2009}. We construct a random graph by uniformly positioning $|\mathcal{N}|$ nodes over the unit square $[0,1]^2$ and randomly sampling edges to construct a set of {feasible} edges $\mathcal{E}_0$: We iteratively sample edges from the $k$-nearest neighbors graph (with $k=6$) until we obtain a connected graph to ensure the feasibility of our instances. 
The construction cost for each edge, $c_{i,j}$ is drawn uniformly from $\mathcal{U}(1,4)$.
Each commodity $k \in \mathcal{K}$ corresponds to an all-to-one shortest path problem with a single destination node $i_k \in \mathcal{N}$. For commodity $k$, we independently sample demands from all nodes $i' \in \mathcal{N}$, $d^k_{i'}$, uniformly between 5 to 20. 
We set $d^k_{i_k} := - \sum_{i' \neq i_k} d^k_{i'}$. We generate $\mathcal{R}$ demand scenarios for each commodity accordingly. 
This process is repeated for every scenario $r \in \mathcal{R}$.
Flow circulation costs, $f_{i,j}^k$, are proportional to the edge length (by a factor 10). 
The capacity of each arc is scaled based on the maximum cumulative demand across all scenarios: $B_{i,j} := \sum_{k \in \mathcal{K}} \sum_{(i,j) \in \mathcal{E}} \displaystyle \max_{r \in \mathcal{R}} d_{ij}^{k,r}$. Formally, we sample the capacity for arc $(i,j)$ according to $u_{ij} \sim \mathcal{U}(1,4) \cdot {B_{i,j}} / |\mathcal{E}_0| $. We fix the cardinality constraint to $c_0 = 2|\mathcal{E}_0|$. 

\subsubsection{{\bf R} instances}
As in \citet{crainic_partial_2021,crainic2016partial,boland_proximity_2016}, we use the {\bf R} instances from classes 4 to 10. 
Each class corresponds to a particular network with its set of nodes, arcs, and commodities. Within each class, the library contains nine instances, associated with different set of arc capacities, edge construction costs, and flow transportation costs. Precisely, each instance within each class is associated with a `class minor' ranging from 1 to 9 and corresponding to increasing ratios of fixed to variable cost and of total demand to total capacity (e.g., instance R4.1 or R10.9).
Each class also contains one vector of nominal demands for each commodity (the same nomial for all instances within each class). However, we need samples of demand \emph{scenarios}. Instead of implementing an ad-hoc sampling scheme, we use the same scenarios as those generated by \citep{rahmaniani2018accelerating} and available at \url{https://github.com/Ragheb2464/R-Instances}. For each class, these files contain 1,000 scenarios for different level of correlation between commodity demand.

Overall, each instance is characterized by a class number (e.g., R4), which determines the network, a class \emph{minor} (e.g., R4.1), and a correlation between demands. In our experiments, we consider 7 classes, 3 class minors (1, 3, and 9), and 2 correlation values (0.0 and 0.8).

\subsection{Computing requirements for each experiment}\label{ssec:a.exp.memories}
This section provides a breakdown of the computational resources allocated for the experiments described in \ref{ssec:exp.implementation}. Each CPU core of the MIT Supercloud Cluster corresponds to 4GB of allocated RAM.

\subsubsection{Synthetic instances} For the experiments on synthetic instances, we request a number of CPU cores/memory that is increasing with the number of nodes in the network, $|\mathcal{N}|$, as described in Table \ref{tab:memory_synthetic}. 
\begin{table}[H]
    \centering
    \caption{{Number of CPU Cores and Memory (GB) allocated for the experiments on synthetic instances, as a function of the number of nodes in the network, $|\mathcal{N}|$.}}
    \label{tab:memory_synthetic}
    \begin{tabular}{rrr}
        \textbf{$|\mathcal{N}|$} & \textbf{Memory (GB)} & \textbf{\# Cores} \\
        \midrule
        10   & 4  & 1  \\
        30   & 8  & 2  \\
        50   & 16 & 4  \\
        70   & 20 & 5  \\
        100  & 28 & 7  \\
        150  & 40 & 10 \\
        200  & 52 & 13 \\
        300  & 76 & 19 \\
        500  & 128& 32 \\
        700  & 176& 44 \\
        \hline
    \end{tabular}
\end{table}

\subsubsection{{\bf R} instances}
For the experiments on the {\bf R} instances, we request a number of CPU cores/memory that is increasing with the number of scenarios, $|\mathcal{R}|$, as summarized in Table \ref{tab:memory_synthetic}. 
\begin{table}[H]
    \centering
    \caption{{Number of CPU Cores and Memory (GB) allocated for the experiments on the {\bf R} instances, as a function of the number of scenarios, $|\mathcal{R}|$.}}
    \label{tab:memory_r_instances}
    \begin{tabular}{rrr}
        \textbf{$|\mathcal{R}|$} & \textbf{Memory (GB)} & \textbf{\# Cores} \\
        \midrule
        10  & 20 & 5  \\
        30  & 20 & 5  \\
        50  & 20 & 5  \\
        100 & 20 & 5  \\
        200 & 28 & 7  \\
        500 & 64 & 16 \\
        1000& 64 & 16 \\
        \end{tabular}
\end{table}

\subsection{Comparison of Different Stochastic Cutting-Plane Algorithms} \label{ssec:a.exp.warmstart}
In Section \ref{ssec:exp.warmstart}, we benchmark the performance of different variants of the stochastic cutting plane algorithm (namely the multi-, single-, and $k$-cut algorithms) with different warm-starting strategies at the root node. 
Recall that we terminate our algorithm after 7,200 seconds or as soon as it achieves an optimality gap of with confidence level $\alpha=0.90$.

Accordingly, the average computational time reported in Table \ref{tab:root_cuts_runtime} are capped at 7,200 seconds whenever the algorithm does not converge within this time limit. To appreciate this censoring issue, Table \ref{tab:root_cuts_conv} presents the fraction of instances solved to $\epsilon$-optimality for each combination of algorithm and warm-start strategy. 
\begin{table}[ht]
	\centering
	\caption{Percentage (in \%) of instances for which the algorithm converged within the time limit (7,200 seconds), for the multi-, single-, and accelerated multi-cut stochastic cutting plane algorithms, with different warm-start strategies at the root node (none, single, or multi). Metrics are averaged across instances with the same number of nodes $|\mathcal{N}|$.}
	\label{tab:root_cuts_conv}
     \begin{tabular}{l|rrr|rrr|rrr}
        \toprule
         & \multicolumn{3}{c}{Multi-Cut}& \multicolumn{3}{c}{Single-Cut} & \multicolumn{3}{c}{Accelerated Multi-Cut} \\
       $|\mathcal{N}|$ &  None & Multi & Single & None & Multi & Single & None & Multi & Single \\ 
       \midrule
         10 &    100     &          100.00  &            80.00 &       100.00 &           100.00 &            100.00 &        100.00 &           100.00 &              80.00    \\
         30 &    33.33   &          40.00   &            26.67 &        26.67 &            40.00 &             20.00 &        33.33  &            46.67 &              33.33    \\
         50 &    26.67   &          53.33   &            33.33 &        26.67 &            40.00 &             33.33 &        20.00  &            60.00 &              40.00   \\
         70 &    20.00   &          45.00   &            45.00 &        25.00 &            35.00 &             25.00 &        20.00  &            45.00 &              45.00 \\
         \midrule     
        100 &    15.00   &          45.00   &            50.00 &        25.00 &            30.00 &             40.00 &        20.00  &            50.00 &              45.00  \\
        150 &    10.00   &          45.00   &            45.00 &        25.00 &            25.00 &             35.00 &        15.00  &            40.00 &              50.00  \\
        200 &    15.00   &          30.00   &            50.00 &        20.00 &            35.00 &             45.00 &        10.00  &            50.00 &              50.00  \\
        \bottomrule
    \end{tabular}   
\end{table}

For the deterministic cutting-plane method (Benders' decomposition), numerous acceleration strategies have been proposed in the literature, e.g., based on valid inequalities (VI) or Magnanti-Wong Pareto-optimal cuts \citep[MW;][]{magnanti_network_1984}. We implemented different variants of the scheme based on this strategy and report their performance in Table \ref{tab:overall.deterministic.variants}. Consistent with \cite{papadakos2008practical}, we find that Magnanti-Wong cuts often do more harm than good, and accordingly we do not include them in our implementation in the main text.
\begin{table}
	\centering
	\caption{Runtime (in seconds) and final optimality gap (in \%) for different variants of the deterministic cutting-plane algorithm, averaged over instances with the same number of nodes $|\mathcal{N}|$. In addition to the naive implementation of the algorithm (with lazy callbacks), we consider adding valid inequalities (VI), Magnanti-Wong Pareto dominating cuts (MW), or both (VI+MW). }
	\label{tab:overall.deterministic.variants}
    \begin{tabular}{lrcrcrcrc}
    \toprule 
     Speed-up & \multicolumn{2}{c}{{None}} & \multicolumn{2}{c}{VI} & \multicolumn{2}{c}{MW}  & \multicolumn{2}{c}{MW + VI}\\ 
     \cmidrule(lr){2-3} \cmidrule(lr){4-5} \cmidrule(lr){6-7} \cmidrule(lr){8-9} 
     $|\mathcal{N}|$ & Runtime & Gap & Runtime & Gap & Runtime & Gap & Runtime & Gap \\ 
    \midrule   
    10 &           247.79  & 0.02          &    326.13	& 0.02                     &   454.61	& 0.02                          &      510.20	& 0.02  \\
    30 &           7163.94 & 6.22          &    7181.16	& 14.86                    &   7200.00	& 7.72                          &      7178.66	& 21.95 \\
    50 &           7200.00 & 4.87          &    7200.00	& 18.47                    &   7200.00	& 6.29                          &      7200.00	& 21.53 \\
    70 &           7200.00 & 11.85         &    7186.38	& 41.29                    &   7200.00	& 12.93                         &      7200.00	& 44.50 \\
    \midrule
    100 &          7165.78 & 16.37         &    7200.00	& 49.15                    &   7200.00	& 16.47                         &      7200.00 &  51.00 \\
    150 &          7186.61 & 23.49         &    7196.71	& 59.96                    &   7200.00	& 26.59                         &      7200.00 &  60.97 \\
    200 &          6853.71 & 26.68         &    7138.72	& 60.39                    &   7039.69	& 28.77                         &      7188.38 &  66.20 \\
    \midrule
    300 &          6237.87 & 23.11         &    6832.63 & 63.33                    &   6818.90	& 24.39                         &      7194.16 &  71.60 \\
    500 &          6441.49 & 49.09         &    6510.74 & 81.27                    &   6874.99	& 51.11                         &      7192.51 &  75.37 \\
    700 &          6499.08 & 53.39         &    6947.16 & 87.48                    &   7012.66	& 69.53                         &      6918.72 &  87.49 \\
    \bottomrule
\end{tabular}
\end{table}

\FloatBarrier
\subsection{Benchmarking Scalability on Synthetic Instances}\label{subsection:a.exp_results}
To verify the correctness of our implementation, we use the smallest instances to verify that all methods terminate with the same optimal solution. To this end, Table \ref{tab:small_scale} reports the optimality gap (in \%) and final objective value for each algorithm, averaged over instances with the same number of nodes $|\mathcal{N}|$ and for which Gurobi converged to within $5\%$ of optimality.
\begin{table}[ht]
	\centering
	\caption{Optimality gap (in \%) and final objective value for each algorithm, averaged over synthetic instances with the same number of nodes $|\mathcal{N}|$, where Gurobi converged to within $5\%$ of optimality.}
	\label{tab:small_scale}
    \begin{tabular}{lrcrcrc}
    \toprule
     & \multicolumn{2}{c}{\texttt{Gurobi} with \eqref{eqn:prob_description}}& \multicolumn{2}{c}{ Deterministic} & \multicolumn{2}{c}{Stochastic} \\ 
     \cmidrule(lr){2-3} \cmidrule(lr){4-5} \cmidrule(lr){6-7} $|\mathcal{N}|$ &    Gap &  Objective & Gap & Objective &  Gap &  Objective   \\
    \midrule
    10 &      0.00   &  10,502.04     &        0.02   &    10,502.04   &      0.23   &      10,509.23  \\
    30 &      42.68  &  1,206,452.87   &        6.22   &   441,821.98   &      4.30   &     404,461.28  \\
    50 &      67.71  &  3,994,890.55   &        4.87   &   1,109,686.50  &      3.72   &    1,025,044.43  \\
    70 &      77.56  &  19,248,454.24  &        11.85  &   5,085,529.77  &      12.12  &    5,055,830.36  \\   
    \bottomrule
\end{tabular}
\end{table}

Figure \ref{fig:max_problems_experiments_variable_nC.slim} (resp. \ref{fig:max_problems_experiments_variable_nC.hybrid}) illustrates the scalability of our single-cut method (resp. accelerated multi-cut) with respect to the number of scenarios and commodities and depicts the optimality gap achieved depending on the total number of nodes $|\mathcal{N}|$ (horizontal axis), the number of scenarios $|\mathcal{R}|$, and the number of commodities $|\mathcal{K}|$. We observe that the complexity of the problem (measured in terms of the final optimality gap) increases with all three problem dimensions, with the number of commodities and nodes having the most noticeable impact. 

\begin{figure}
	\centering
	\includegraphics[width=0.8\textwidth]{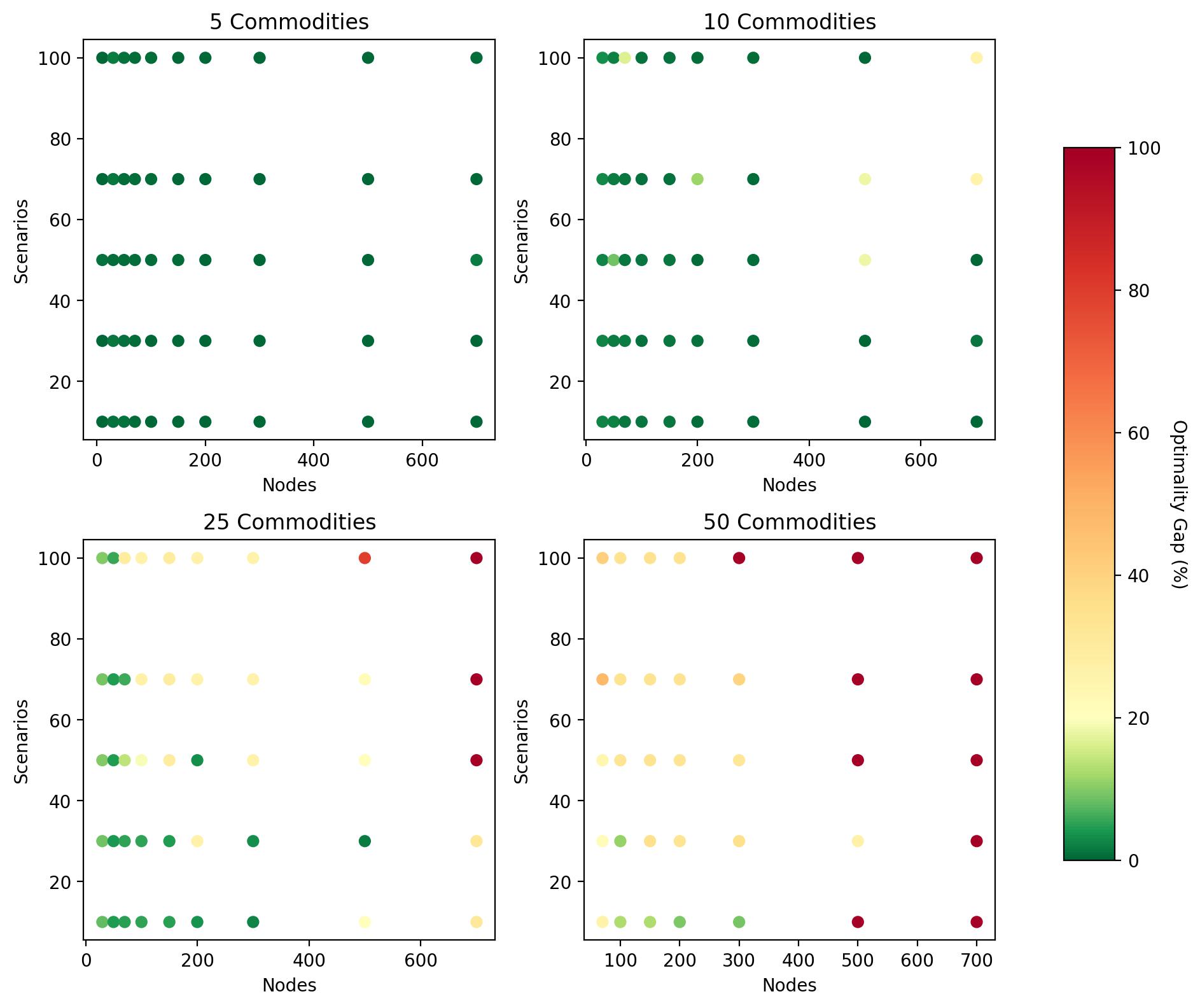}
	\caption{Optimality gaps achieved by the single-cut stochastic cutting plane algorithm on all synthetic instances. For each combination of number of nodes $|\mathcal{N}|$, number of commodities $|\mathcal{K}|$, and number of scenarios $|\mathcal{R}|$, results are averaged across 3 random instances.}
	\label{fig:max_problems_experiments_variable_nC.slim}
\end{figure}
\begin{figure}
	\centering
	\includegraphics[width=0.8\textwidth]{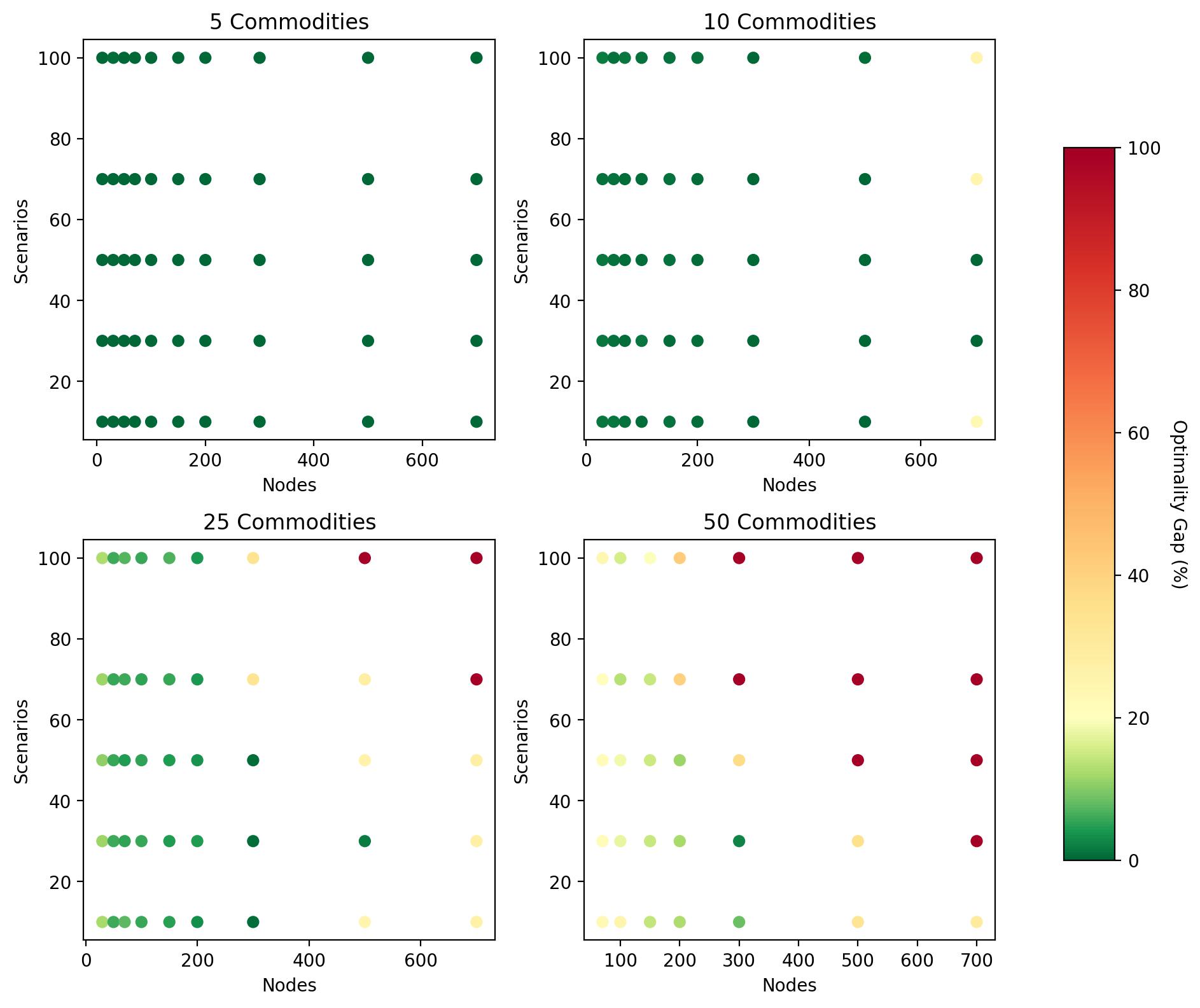}
	\caption{Optimality gaps achieved by the accelerated multi-cut stochastic cutting plane algorithm on all synthetic instances. For each combination of number of nodes $|\mathcal{N}|$, number of commodities $|\mathcal{K}|$, and number of scenarios $|\mathcal{R}|$, results are averaged across 3 random instances.}
	\label{fig:max_problems_experiments_variable_nC.hybrid}
\end{figure}

\FloatBarrier
\subsection{Benchmarking on the Instances from \cite{crainic_simplex-based_2000} }\label{ssec:a.exp.rinstances}
 
In Section \ref{ssec:exp.rinstances}, we compare the performance of different cutting-plane algorithms for solving stochastic network design instances from the {\bf R} instances. One alternative is to solve the perspective reformulation \eqref{eqn:prob_description} with a mixed-integer second-order cone solver like Gurobi. As displayed in Figure \ref{fig:R.gurobiVSbenders}, however, it performs worse than, e.g., a deterministic Benders decomposition scheme.
\begin{figure}
    \centering
    \includegraphics[width=0.5\textwidth]{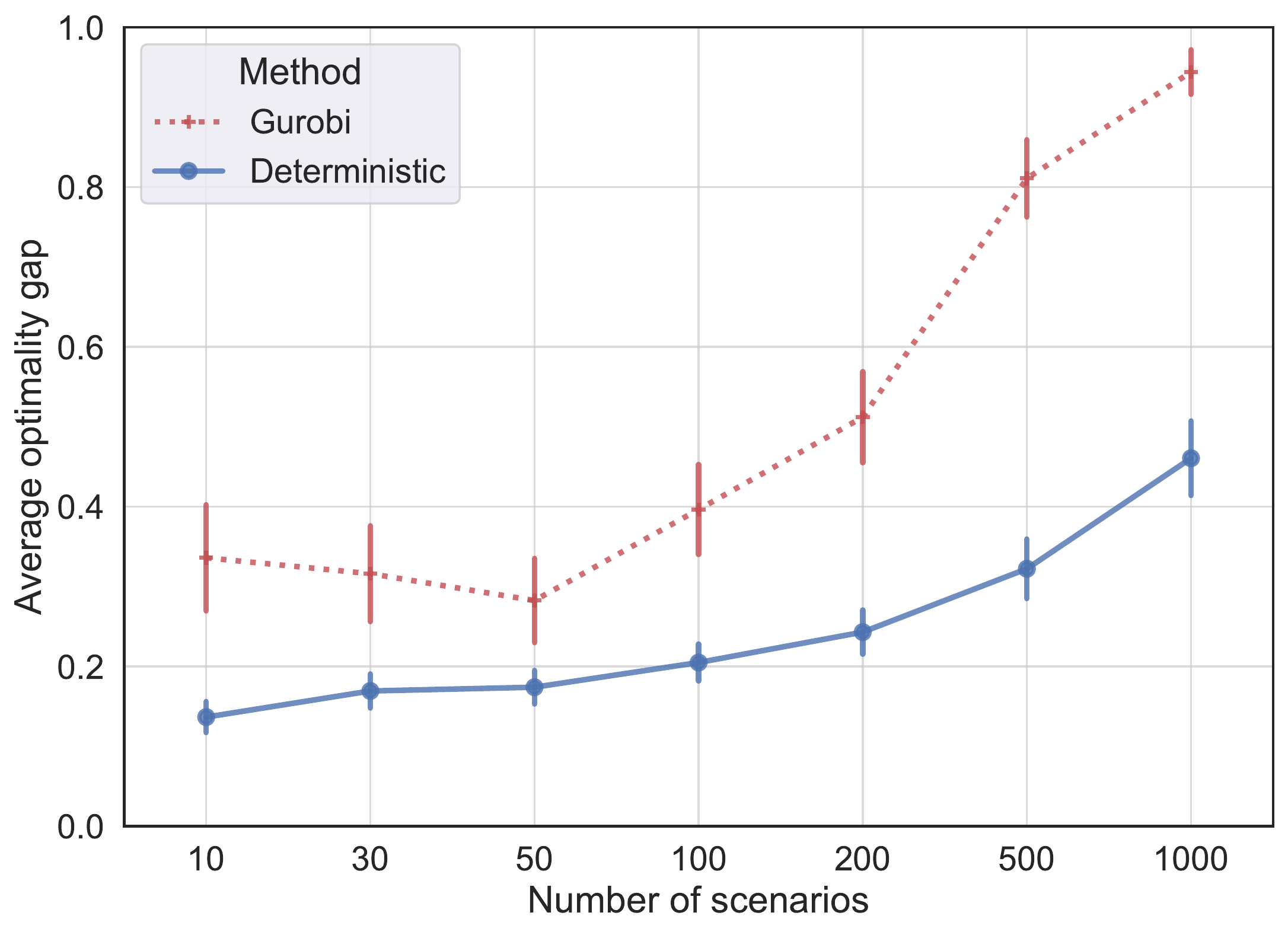}
    \caption{Average optimality gap on the R instances, achieved by the Gurobi on the formulation \eqref{eqn:prob_description} compared with the deterministic cutting-plane algorithms, for different number of scenarios. Bars represent standard errors.}
    \label{fig:R.gurobiVSbenders}
\end{figure}

Figure \ref{fig:R.benders.distribution} reports the distribution of the optimality gaps achieved by each method, over all {\bf R} instances. We find that the ratio of fixed to variable cost and of total demand to total capacity, as controlled by the `class minor' of each instance, is the main driver of the instance complexity, with a higher minor (i.e., higher ratios) resulting in larger optimality gaps (i.e., harder instances). 
\begin{figure}
    \centering
    \includegraphics[width=0.9\textwidth]{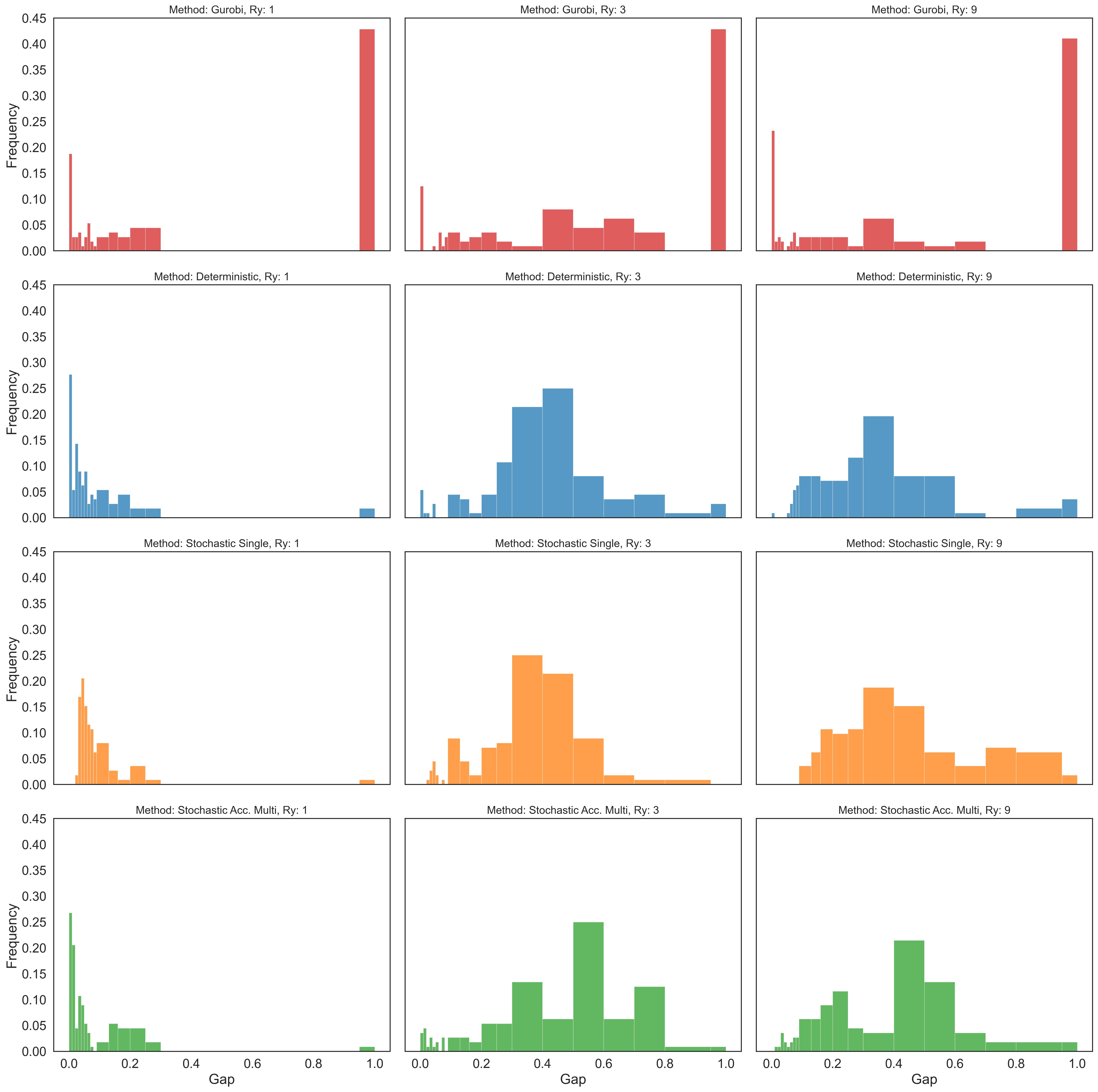}
    \caption{Distribution of the optimality gap achieved on the R instances by the perspective reformulation, the deterministic Benders decomposition, and our stochastic (single-cut and accelerated multi-cut) cutting-plane algorithms. Results are grouped according to the instance class minor (see definition in Section \ref{ssec:a.instance}).}
    \label{fig:R.benders.distribution}
\end{figure}

\end{APPENDICES}

\end{document}